\documentclass[trsc,nonblindrev]{informs3} 

\OneAndAHalfSpacedXI 

\usepackage[linesnumbered,ruled,vlined]{algorithm2e}

\SetCommentSty{mycommfont}
\usepackage{listings}
\usepackage{pdflscape}
\usepackage{xcolor}
\usepackage{multirow}
\usepackage{soul}

\usepackage{natbib}
 \bibpunct[, ]{(}{)}{,}{a}{}{,}%
 \def\bibfont{\small}%
\TheoremsNumberedThrough     

\EquationsNumberedThrough    

\MANUSCRIPTNO{} 

\begin{document}


\RUNAUTHOR{Martin-Iradi, Pacino, and Ropke}

\RUNTITLE{The multi-port berth allocation problem with speed optimization}

\TITLE{The Multi-Port Berth Allocation Problem with Speed Optimization: Exact Methods and a Cooperative Game Analysis}

\ARTICLEAUTHORS{%
\AUTHOR{Bernardo Martin-Iradi}
\AFF{DTU Management, Technical University of Denmark, Akademivej Building 358, 2800 Kgs. Lyngby, Denmark, \EMAIL{bmair@dtu.dk} \URL{}}
\AUTHOR{Dario Pacino}
\AFF{DTU Management, Technical University of Denmark, Akademivej Building 358, 2800 Kgs. Lyngby, Denmark, \EMAIL{darpa@dtu.dk} \URL{}}
\AUTHOR{Stefan Ropke}
\AFF{DTU Management, Technical University of Denmark, Akademivej Building 358, 2800 Kgs. Lyngby, Denmark, \EMAIL{ropke@dtu.dk} \URL{}}
} 

\ABSTRACT{%
We consider a variant of the berth allocation problem \textemdash i.e., the multi-port berth allocation problem\textemdash aimed at assigning berthing times and positions to vessels in container terminals. 
{This variant involves optimizing vessel travel speeds between multiple ports, thereby}
exploiting the potentials of a collaboration between carriers (shipping lines) and terminal operators.
Using a graph representation of the problem, we reformulate an existing mixed-integer problem into a generalized set partitioning problem, in which each variable refers to a sequence of feasible berths in the ports that the vessel visits.
By integrating column generation and cut separation in a branch-and-cut-and-price procedure, our proposed method is able to outperform commercial solvers in a set of benchmark instances and adapt better to larger instances.
In addition, we apply cooperative game theory methods to efficiently distribute the savings resulting from a potential collaboration and show that both carriers and terminal operators would benefit from collaborating. 
}%


\KEYWORDS{Transportation, Exact methods, Container terminal, Berth allocation problem, Speed Optimization, Cooperative game theory}

\maketitle

%


\section{Introduction}
\label{sec1}
The International Maritime Organization (IMO), in its fourth climate report \citep{IMO2020}, reflects on the increase in shipping's $CO_2$ emissions in the recent years. In the period 2012-2018 the shipping's total emissions have increased by 9.6\%. This alarming trend highlights the need for pursuing the strategies that the IMO adopted in 2018 for reducing greenhouse gas (GHG) emissions from ships \citep{IMO2018}.
The aim is to reduce total emissions from shipping by 50\% in 2050, and to reduce the average carbon intensity by 40\% in 2030 and 70\% in 2050, compared to 2008. 
Yet world maritime trade keeps growing at an annual average of 3\% reaching a record high of 11 billion tons of total volume in 2018 \textemdash a number that translates into almost 800 million twenty-foot equivalent units (TEUs) handled in container ports worldwide \citep{UNCTAD19}.
Given that trade volume has steadily increased since then, the need for more efficient and sustainable operations in maritime transport logistics is essential \citep{bekta2019a}.

From the terminal viewpoint, the growth in container trade involves more or larger vessels arriving at ports, in need of berthing. One solution to satisfying the increasing demand is to extend the existing quay. The problem is that doing so usually requires an expensive investment and sometimes may not even be physically feasible. An alternative strategy is to improve the efficiency of existing resources through optimization techniques that do not entail costly investment.

The berth planning of a terminal can be modelled mathematically as the Berth Allocation Problem (BAP). In the BAP, the aim is to assign incoming ships to berthing positions along the terminal.
\cite{steenken2004a} define this problem as highly critical within container container terminal planning logistics, due to the scarcity of berthing space.
\begin{figure}[th]
    \centering
    \includegraphics[width=0.8\textwidth]{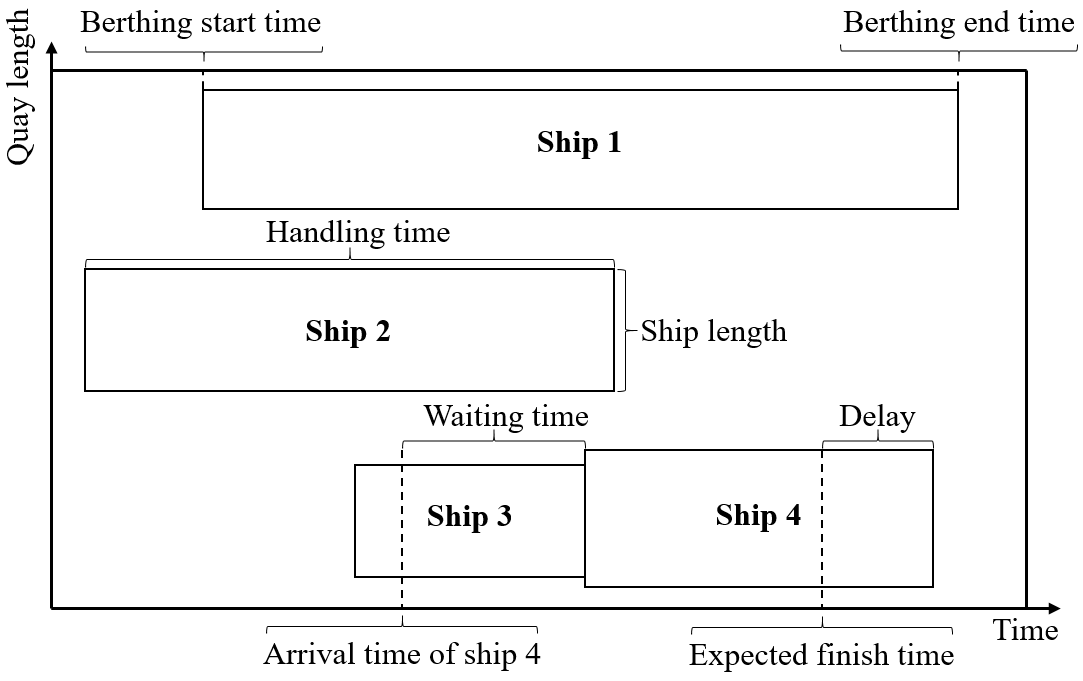}
    \caption{Example solution of the BAP for a port terminal with four vessels.}
    \label{fig:BAPex}
\end{figure}
Figure \ref{fig:BAPex} illustrates the problem in a two-dimensional diagram where one dimension is space (quay length), and the other one is time (the planning horizon). 
We depict each ship as a rectangle whose dimensions are the ship length and handling time, the time the ship spends at the berth (i.e., during unloading and loading)
Each ship usually has a fixed time window defined by its expected start and finish time.
Although ships can arrive before their berthing time, they will need to wait at the port. 
Similarly, ships can be allowed to exceed the expected finish time incurring in a delay.
We denote the entire time that the ship spends at the port (i.e waiting time plus berthing period) as the "service time."
Any non-overlapping positioning of the ship rectangles within the decision space defines a feasible solution for the BAP.

We can classify the BAP variants according to how the berths are distributed along the quay. In the \textit{discrete} BAP, we divide the quay into a discrete set of berths, with only one ship allowed to one berth at a time, whereas in the \textit{continuous} BAP, the ships can berth anywhere along the quay as long as they maintain {a safe distance} from the other ships. Moreover, the BAP can be either \textit{static} or \textit{dynamic}. In the \textit{static} variant, we assume that all ships are at the port when the berth planning is done, whereas in the \textit{dynamic} version, we assume that ships can arrive while the planning is in process.

{The efficient planning of a terminal requires the vessels to abide by their schedules. Thus, efficient vessel scheduling is also a critical aspect, not only for the carriers, but also for the terminal operators. 
The design of vessel schedules can be modeled mathematically as the Vessel Scheduling Problem (VSP).
The VSP aims at determining the sailing speeds between consecutive ports in the route (i.e., voyage legs) in order to optimize the vessels' fuel consumption and turnaround time at port and the number of vessels required to operate the route with a given frequency.}

Both the increasing volume of container trade and the up-sizing of the vessels have led to increased competition among container terminals.
each vying to become the port of call for more vessels
\citep{notteboom2017a}. As a result, most terminals are reticent to share information with other terminals and prefer to plan their operations independently. 
Terminals commonly plan berth allocation based on ship schedules. Nevertheless, these schedules are subject to a level of uncertainty,
because different types of disruptions\textemdash such as weather conditions or technical problems at the terminal\textemdash can alter the schedules and result in delays. 
When each terminal does its planning independently, a delay in one terminal can potentially be propagated through the shipping service to other ports \citep{notteboom2009a} or incur higher fuel costs for the carriers (shipping lines) if they need to increase the vessel's speed to make up for lost time. For example, a vessel stopping in ports A and B may encounter a congested terminal when arriving at port A and become delayed. The carrier can then order the vessel to either speed up to arrive at port B on time, entailing higher fuel consumption, or arrive late at port B, forcing the terminal to modify its berthing plan.

A potential solution to avoid this type of {scenario} is to establish some form of collaboration between players in the maritime industry. 
Collaboration can be established not only between same type of stakeholders (i.e., between multiple carriers) but also between more players (i.e., carriers and terminal operators). 
{The \cite{shippingcouncil2015} encourages terminals to establish collaborative agreements with carriers as one of the main 
ways of reducing port congestion and improving planning eﬃciency.}

{A certain degree of collaboration is assumed in the VSP, however, the problem does not explicitly consider the berth allocation at the terminal and this
can lead to a significant increase in service time. Therefore, integrating the BAP together with the scheduling of the vessels becomes relevant.} 
Sharing information allows planners to simultaneously plan the berthing at the terminals and be able to minimize disruptions and reduce costs and emissions. 
{Recent studies show that collaboration between carriers and terminals can 
lead to significant benefits for both
\citep{dulebenets2019a}.}
{This is the goal of the Multi-Port Berth Allocation Problem (MPBAP), first introduced by \cite{venturini2017a}, which
simultaneously plans the berth allocation of multiple ports taking into account the vessels' speed.
} 

{The MPBAP can either be applied either when one company controls both vessels and terminals, or by a third-party service provider which works as an orchestrator. An example of the former is \emph{Maersk}, owning both the carrier \emph{Maersk Line} \citep{maersk} and the terminal operator \emph{APM Terminals} \citep{apmt}. }

{
At present, there are companies in the market that offer optimization-based planning software separately to carriers and terminal operators \citep{portchain, navis, sealytix, tgi, rbs}. Such companies already have access to all the necessary data for the MPBAP, which makes them excellent candidates to orchestrate the collaboration. Since both carriers and terminals are already sharing data with those companies, trust issues should be minimal, but customers should of course be free to decline that their data is used in a joint optimization problem. The amount of flexibility that terminals and carriers are willing to commit to the collaboration, can easily be modeled with the time windows, making the MPBAP if not an operational tool, at least a tool to identify the potential savings.
}

{
To make the service attractive to customers, the software company needs to show that the collaboration is beneficial for all involved parties. Therefore, we apply cooperative game theory to demonstrate that the total costs in the MPBAP solution can be shared in a favorable way. Using this service only requires that participating carriers and terminal operators allow the third party to jointly use their data but does not entail sharing additional data or the disclosure of the customer’s data to other customers. Once the operations conclude, the third party would be in charge of returning the savings according to the initial calculations. 
}

{
Similar collaboration mechanisms have also been studied in the road transportation sector. \cite{ergun2007a} study collaborative logistics in truck transportation where part of the carriers’ savings are returned to the shippers. \cite{oezener2011a} also propose collaborative models where players receive more favorable rates in return.
In fact, they indicate that a centralized decision-maker with complete information about all participants would be ideal for collaborative models to work. However, in their study, \cite{oezener2011a} suggest that, due to lack of trust, players may not be willing to share additional information and therefore, they explore different mechanisms. Fortunately, this lack of trust is minimized in our case as the players already share the required information with the third party.
}

In studying the MPBAP, this paper makes four contributions. First, we present two new formulations for the MPBAP, based on a graph representation. Second, we propose exact methods based on column generation, together with branching, cutting, and symmetry-breaking enhancements. Third, we demonstrate the quality of our method by comparing it to a commercial solver and testing it through both a set of benchmark instances from a previous study and a new set of harder instances. Fourth, to demonstrate the benefits for both carriers and terminal operators in a scenario of a joint grand coalition, we apply {cost} allocation methods from cooperative game theory.

The structure of this paper is as follows. {Section \ref{sec:lr} reviews the state-of-the-art studies on berth allocation, speed optimization and collaboration on the shipping industry.} Section \ref{sec:ProbDes} describes the MPBAP by presenting two mathematical formulations, together with the one from \cite{venturini2017a}. Section \ref{sec:SolM} gives our solution method, and Section \ref{sec:cgt} introduces and discusses the cooperative game methods used for effectively distributing the {costs} of a coalition. Section \ref{sec:results} compares the models’ performance through extensive computational experiments and analyzes the cooperative game theory results. Section \ref{sec:conclusion} concludes by briefly discussing both the findings and possible future research directions.

\section{Literature review} \label{sec:lr}
This section has been divided into three. First, we describe the main studies related to the BAP. Secondly, we cover the literature concerning speed optimization, and the last part focuses on collaboration studies within the container shipping industry and literature where cooperative game theory has been applied to it.

\subsection{BAP literature}

The berth allocation problem {is known to be NP-hard \citep{lim1998a, hansen2015a}} and has received significant attention in the last two decades. \cite{carlo2014a} and \cite{bierwirth2015a} presented detailed literature surveys on the seaside operations of container terminals such as the BAP where they emphasized the raising interest on this particular problem in the last years.
\cite{imai2005a} conducted the first study considering a continuous BAP and  
\cite{cordeau2005a} studied both the discrete and continuous version of the problem and solved them through heuristic methods. 
\cite{guan2005a} presented a tree search exact method that performed better than commercial solvers and an efficient composite heuristic method.
\cite{du2015a} extended the problem to also include the effect of tides and adopted the \textit{virtual arrival} policy that is currently used in many terminals worldwide.
\cite{cheong2010a} considered priorities for each of the vessels. The BAP is optimized using an evolutionary algorithm that minimizes the make-span, the waiting time and the deviation from a reference schedule. 
\cite{buhrkal2011a} compared three different methods for the discrete BAP and showed that a generalized set-partitioning model outperforms the rest. 
\cite{saadaoui2015column} reformulated the problem into a set packing problem where variables refer to assignments of ships to berthing positions and solved it using delayed column generation. In our paper, we combine the applicability of column generation procedures using a generalized set partitioning problem formulation.
Regarding the discretization of the quay, \cite{kordi2016a} presented a hybrid variant of the BAP where ships can only berth in a given set of positions.
\cite{lalla-ruiz2016b} studied how the tides can limit the time available for ships to berth given their draft and the water depth and solved this variant of the BAP using a generalized set partitioning problem formulation.
{The multi-port version of the BAP studied in this paper was first defined by \cite{venturini2017a}. The mixed integer problem formulation they presented is used as a reference for the ones considered in this paper.} 
\cite{kramer2019a} proposed two new formulations for the discrete BAP: a time-indexed formulation and an arc-flow formulation that seem to perform better than the methods from \cite{buhrkal2011a}.
\cite{corry2019a} proposed a mixed integer problem formulation for the BAP with channel-constrained ports where the sequencing of channel movements is also optimized.

\subsection{Speed optimization literature}

The relation between vessel speed and fuel consumption is non-linear. Since fuel emissions are directly proportional to the fuel burnt, optimizing sailing speed becomes relevant from the carrier and environmental perspective. 
The policies of the IMO in the last years have raised debate on which measures to implement regarding speed optimization, speed reduction or slow steaming. In that aspect, multiple studies have been done analyzing the aspects and impacts of the different measures.
Based on the scenario of slow steaming, \cite{kontovas2011a} investigated a berthing policy that aims at reducing the waiting time at port.
\cite{psaraftis2013a}, \cite{wang2013a}, \cite{psaraftis2015b} and \cite{psaraftis2015a} presented taxonomies and surveys on speed models in the maritime transportation sector where the impacts and main trade-offs of slow steaming are analyzed and decision models proposed.

{
The VSP has speed optimization as its core concept and
the interest in this problem has continued increasing in the last decade \citep{dulebenets2019a}.}
{\cite{fagerholt2001a} presented a mathematical model for the VSP and solved it using a method based on the set partitioning formulation.}
{\cite{wang2014a} extended the VSP to also consider cargo allocation and indicated that carriers should consider the cargo costs  arising from additional waiting time at port.}
{\cite{dulebenets2018b} proposed a multi-objective model considering the route service costs. The results indicated that negotiating the port calls and handling rates with the terminal operator could lead to significant savings.}
{To some extent, the VSP can be seen as a collaborative problem, however, most of the studies focus on the interests of the carrier.}
{A variant of the VSP where shipping line companies and terminal operators collaborate has also been studied recently. This variant assumes that the terminal operator can offer multiple time windows or handling rates to the carrier, instead of the fixed ones considered in the generic VSP.}
{For instance, the MPBAP presented in \cite{venturini2017a} can be included in this problem category where the berth allocation planning is also considered.}
{\cite{dulebenets2019b} presented a mathematical model for the collaborative VSP where terminals offer both multiple port service time windows and handling rates. The results showed the benefits of the collaborative agreement on the liner shipping operations.}

{Environmental aspects have also been addressed in this type of problems.}
\cite{fagerholt2010a} minimized the fuel consumption by optimizing speeds along a shipping route. By discretizing the arrival times at each port, the cubic function relating speed and fuel emissions can be linearized and the problem solved as a shortest path problem. \cite{fagerholt2015a} and \cite{zhen2020a} extended the route and speed optimization study by also considering emission control areas (ECAs). \cite{fagerholt2015a} aimed at minimizing the fuel consumption whereas \cite{zhen2020a} also considered $SO_2$ emissions. Both studies showed that carriers tend to use slow steaming within ECAs or directly avoid sailing through these areas.
\cite{reinhardt2016a} optimized a liner shipping network by adjusting berthing times with the objective of minimizing fuel consumption.
The speed and routing of multiple vessels is optimized in \cite{wen2017a} under a unified objective that minimizes transit times, total costs and fuel emissions. They implemented  a branch-and-price heuristic and a constraint programming model which is tested in a subset of the Mediterranean ports. \cite{du2011a}, \cite{du2015a} and \cite{sun2018a} integrated speed optimization with the BAP by considering that ships still need to sail a certain distance to arrive at port. The second-order cone programming transformation used by \cite{du2011a} to approximate the relation between sailing speed and fuel consumption is improved by quadratic outer approximations in \cite{wang2013b}.

\subsection{Collaboration in the shipping industry}

{The MPBAP introduced by \cite{venturini2017a} can be seen as a problem with a high degree of collaboration and}
the study of different collaborative forms in the container shipping industry has gained interest in recent years.
\cite{song2003a} studied competition and co-operation in ports and coined the term \textit{co-opetition}.
\cite{wang2015a} presented two collaborative methods between shipping line companies and port operators where the aim is to create a win-win situation by balancing the priorities of both parties and encouraging them to share true information.
\cite{lalla-ruiz2016a} proposed a cooperative search for the discrete BAP based on a grouping strategy. Individuals are organized into groups where they can only share information with other individuals from the same group.
\cite{notteboom2017a} investigated alliance formations in container shipping by studying their strategies when choosing ports.
\cite{dulebenets2018a} presented the collaborative berth allocation problem (CBAP), which is a variation of the BAP that also allows to divert vessels to another terminal when there is a peak demand, and solved it using a memetic algorithm.
Collaboration is also studied by integrating berth allocation with other scheduling problems such as ship routing. \cite{pang2014a} studied such integration for a feeder company operating both vessels and container terminals. {This study also considered transhipments of containers but did not cover speed optimization.}

Game theory has also been widely applied in the container shipping industry  \citep{pujats2020a}. In our paper, the focus is on cooperative game theory where the target is on distributing the profits or savings among players.
The studies vary depending on which are the players considered (carriers, terminal operators or both) in the cooperation.
\cite{song2002a} applied cooperative game theory to depict a conceptual framework for liner shipping alliances showing that the core theory is applicable to the liner shipping market.
\cite{saeed2010a} presented a two-stage cooperative game for container terminals within the Karachi Port in Pakistan. The results indicated that a \textit{grand coalition} among all players gives the best payoff for all terminals.
The work by \cite{krajewska2008a} showed, by means of cooperative game theory, that collaboration among road freight carriers is practical and cost-effective for all players. 
\cite{wen2019a} studied the benefits of horizontal cooperation in a shipping pool by not only maximizing the pool profit but also allocating the profits fairly among participants. The profit sharing framework from \cite{krajewska2008a} and some of the profit allocation methods presented in \cite{wen2019a} have been used in this study and, to the best of our knowledge, it is the first time cooperative game theory is applied to the MPBAP.

\subsection{Research gap}
{While the BAP and VSP have been extensively studied in the literature, with an increasing interest in the last decade, very few papers address the potentials of integrating the two problems and only one paper has been found to address the MPBAP. Only an MIP formulation for the problem has been proposed, which shows good performance for small instances but struggles when the size of the instances increases. Therefore, there is a need for a more efficient solution method that scales better to larger instances.
Furthermore, the MPBAP implies collaboration between different parties in the shipping industry and an analysis of the model's applicability in real life is lacking in the literature. Thus, we believe that assessing the stakeholders' incentives to enter into such collaboration is relevant.}

\section{Problem description}\label{sec:ProbDes}

\begin{figure}[th]
    \centering
    \includegraphics[width=\textwidth]{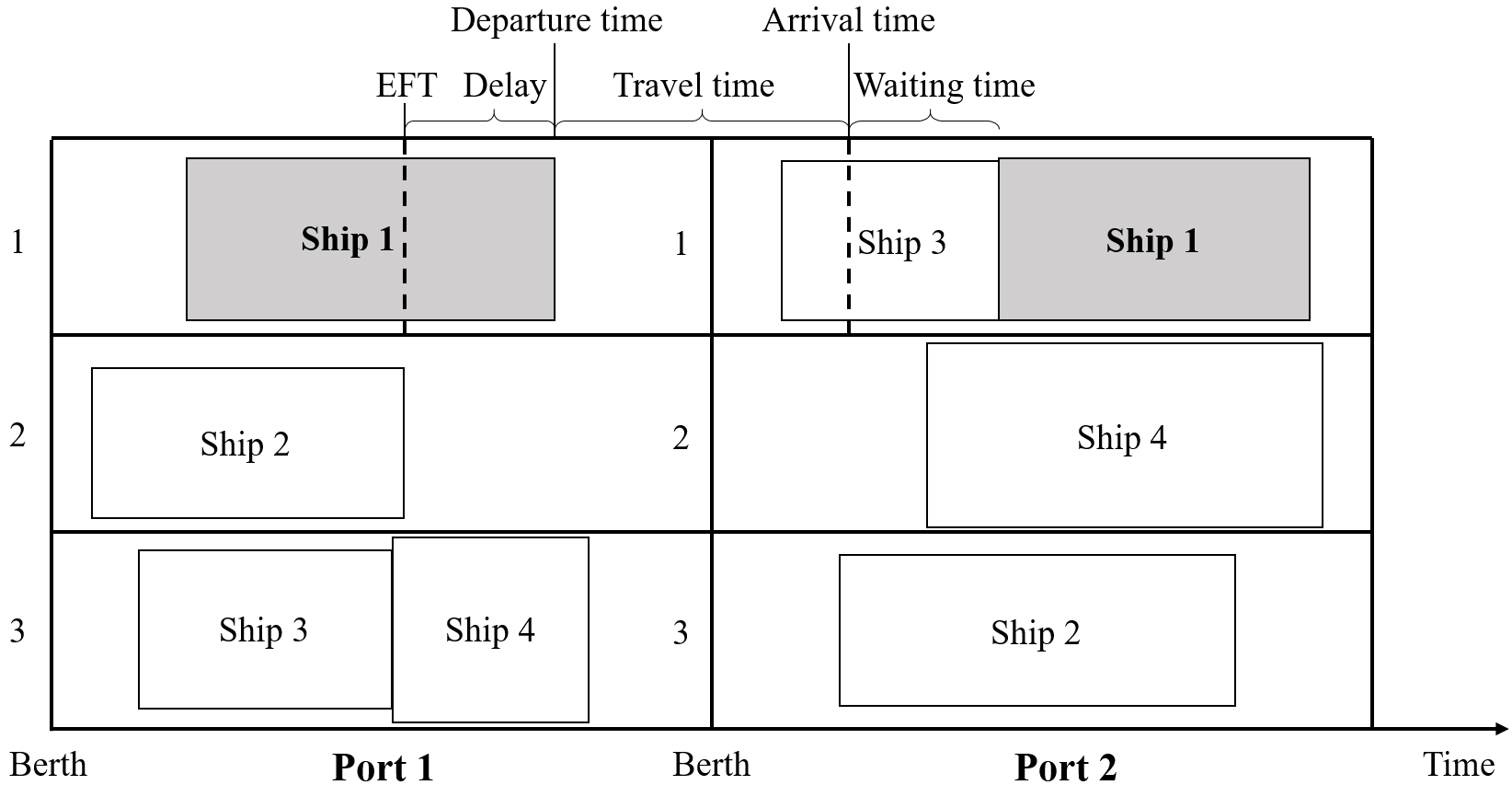}
    \caption{Example solution of the MPBAP for four vessels, two port terminals and three berths per port. The traveling timeline for ship 1 (in gray) is defined at the top, where $EFT$ denotes the \textit{expected finish time} at port 1.}
    \label{fig:CBAPex}
\end{figure}

The MPBAP 
can be seen as 
{a partial integration between the BAP and the VSP.}
Particularly, this study is based on the discrete and dynamic BAP and it is extended to cover multiple ports where the sailing speed between ports is optimized. 
This can be seen as a collaborative approach where information is shared among shipping line and terminal companies.
The main addition of the MPBAP compared to the BAP is the optimization of the sailing speed between ports and the simultaneous planning of multiple ports. 
{Figure \ref{fig:CBAPex} shows a solution example to a problem with four ships and two ports, each having three berthing positions.
As shown for ship 1, the travel time, which depends on the chosen sailing speed, determines the arrival time to the next port and this can constrain the available berthing time window further.}
The MPBAP aims at minimizing the total costs for both the carriers and terminal operators. 

\subsection{Fuel consumption model}\label{sec:fuelCon}
One of the main costs for a carrier is the fuel. The fuel consumption is directly linked to the sailing speed but not in a linear way. Thus, we need an accurate model that links the sailing speed with the fuel consumption realistically.

Many studies approximate the fuel consumption as a cubic function of the speed (e.g., \cite{meng2011a}, \cite{wang2012a}, \cite{reinhardt2016a}),
{
\begin{equation}\label{eq:fuelC}
    F(i,\delta) = \bigg(\frac{\delta}{\delta_i}\bigg)^3 \Gamma_i
\end{equation}
}
where equation (\ref{eq:fuelC}) measures the fuel consumption {$F(i,\delta)$ in tons/hour for a given ship $i$. $\delta_i$ is the design speed of vessel $i$ and $\delta$ is the sailing speed, both measured in knots (i.e., nautical miles per hour). Finally, $\Gamma_i$ is the fuel consumption in tons/hour for vessel $i$ at the design speed.}
This approximation is fairly accurate for container ships of limited size and for a range  of sailing speeds that are not significantly slow. In our study, we optimize the sailing speed between ports, where we expect speeds similar to the design speed ({$\delta_i$}) of the vessel and we do not consider the fuel consumption derived from entering or leaving a port where near-zero speeds are used. In order to avoid non-linearity in the mathematical formulation of the problem, we apply a discretization of the cubic approximation based on the one proposed by \cite{venturini2017a}. A set of different speeds $S$ is defined that can be used by ships to travel between ports. The set of speeds correspond to reasonable and realistic speeds in a range around the design speed. Then, for each of the selected speeds $\delta \in S$ and ship $i$, a fuel consumption value ($\gamma_{i,\delta}$) {measured in tons/nautical mile} can be calculated based on the cubic approximation using the following equation (\ref{eq:fuelCdisc}).
{
\begin{equation}\label{eq:fuelCdisc}
    \gamma_{i,\delta} = \frac{F(i,\delta)}{\delta} = \frac{\Big(\frac{\delta}{\delta_i}\Big)^3 \Gamma_i}{\delta}
\end{equation}
}

\subsection{Cost Structure}\label{sec:costs}

{The MPBAP aims at optimizing the operational costs for both carriers and terminal operators.
This Section defines the main costs involved in the problem context and describes to which stakeholder the costs are related.
}
\begin{figure}[th]
    \centering
    \includegraphics[width=0.9\textwidth]{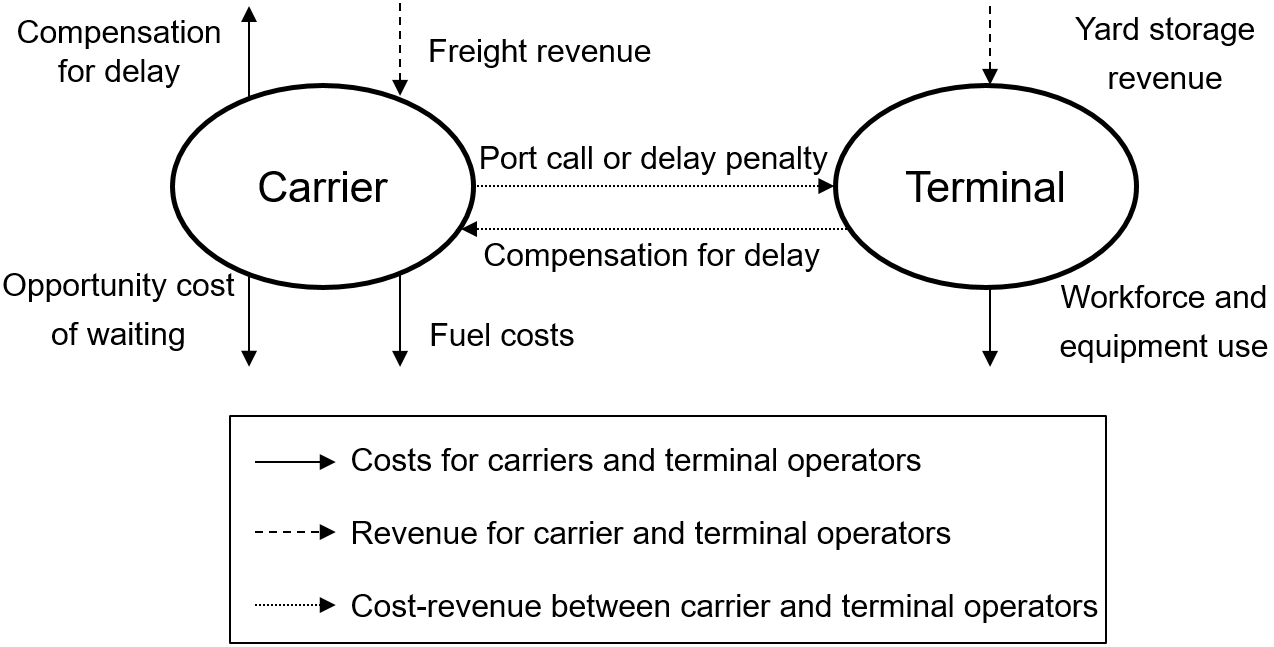}
    \caption{General overview of main costs and revenue for the shipping carriers and terminal operators.}
    \label{fig:CostStruct}
\end{figure}
{An overview of the main sources of cost and revenue for both carriers and terminal operators is shown in Figure \ref{fig:CostStruct}.}

{As mentioned in Section \ref{sec:fuelCon}, the main cost driver for a shipping line company is the fuel consumption which usually accounts for more than 50\% of the carrier's total costs \citep{fagerholt2015b}. Another carrier related cost is the waiting time at anchorage (i.e., waiting to berth at port). As described by \cite{chang2012a}, 
the waiting cost is not only the direct cost of being for longer time at a port, but also the resulting loss of potential income (i.e., opportunity cost).
Regarding the service time at port, this is usually pre-established by a contract or when booking the port call. The cost may differ based on multiple factors such as the number for containers to be loaded and unloaded (i.e., quay crane moves) or the size of the ship and number of cranes required. In this case, the cost can be considered constant for the carrier regardless of the resulting quality of the terminal's planning.
Finally, there are also delay costs associated to the carrier. 
Ending the service after the expected finish time at a port may result in additional payments to the shippers for the delay on the delivery of their cargo.}
In addition there may be other costs arising from delays. On one hand, if the delay at the terminal is caused by the ship arriving late, the carrier may be subject to a fine or delay penalty to the terminal. On the other hand, if the planned service time for a ship gets extended due to, for example, a breakdown of a quay crane or a poor berth allocation plan, the carrier affected by the delay may be entitled to a compensation from the terminal operator. It can be noticed, that these delay costs are paid from the carrier to the terminal operator or vice versa. This means that a cost for one party becomes a {revenue} for the other one. 
{The main premise of this problem is that carriers and terminal operators jointly plan their operations and, therefore, the internal delay costs do not exist and can be excluded from the objective of the problem.} 

The main costs impacting the planning of the terminal in this problem are both the handling and delay costs. 
We identify the use of resources to be directly proportional to the number of workforce and quay crane usage hours. The fewer shifts needed to serve the vessel, the greater the profit is for the terminal. Therefore, both an increasing handling time by the vessels or an increased delay will require additional workforce.
{As mentioned before in this section,}
one of the premises of the MPBAP is that the planning decisions are agreed between the carrier and the terminal operator based on the overall best solution for all. 
{This collaborative optimization removes the concept of delay between the participating players. However, we do consider a delay cost for the terminal operators in the objective of the problem.}
We study the problem from a tactical point of view but assume that, for instance, {workforce} planning at the terminal is performed beforehand. In this scenario, the suggested optimal berth allocation plan for a given terminal may require more {workforce} than initially planned. This will directly translate in the use of additional resources to cover for the additional handling operations that can be computed as delay costs.

{All in all, the MPBAP covers the costs depicted with a continuous line in Figure \ref{fig:CostStruct}. Thus, the objective of the problem focuses on minimizing the fuel consumption and the costs related to waiting, handling and delay time.}

\subsection{Mixed-integer problem formulation. The \cite{venturini2017a} model} \label{sec:mipForm}
The solution method presented in this paper is based on a mixed integer problem (MIP) formulation from \cite{venturini2017a}, which we now briefly present. We first list the notation used in the model:

\begin{table}[th!]
    \centering
\begin{tabular}{ p{3.0cm}  p{12.0cm} }
  \hline
  \multicolumn{2}{c}{Sets and parameters} \\
  \hline
  $N$ & Set of ships \\
  $P$ & Set of ports \\
  $P_i$ & Set of ports to be visited by ship $i \in N$ sorted in visiting order \\
  $B_p$ & Set of berths at port $p \in P$ \\
  $V^{p,b}$ & Set of vertices, $V^{p,b} = N \cup \{ o(p,b), d(p,b)\}$, with $o(p,b) =$ origin node for arcs and $d(p,b) =$ destination node for arcs, both defined for every port $p \in P$ and berth $b \in B_p$ \\ 
  $A^{p,b}$ & Set of arcs $(i,j)$ with $i,j \in {V^{p,b}}, i \neq j$ \\
  $S$ & Set of speeds \\
  $Start^p_i$ & Minimum starting time of activities for ship $i \in N$ at port $p \in P_i$ \\
  $EFT^p_i$ & Expected finishing time of activities for ship $i \in N$ at port $p \in P_i$ \\
  $s^{p,b}$ & Starting time of activities for berth $b \in B_p$ at port $p \in P$ \\
  $e^{p,b}$ & Ending time of activities for berth $b \in B_p$ at port $p \in P$ \\
  $h_i^{p,b}$ & Handling time of ship $i \in N$ at berth $b \in B_p$ at port $p \in P$ \\
  $d^{p,p'}$ & Distance between {pair of ports $p,p' \in P$} \\
  $P_{iL}$ & The last port to be visited by ship $i \in N$ in the route \\
  $\gamma_{i,\delta}$ & Fuel consumption per unit of distance for ship $i \in N$ at speed $\delta \in S$ \\
  $\Delta_{\delta}$ & Travelling time per unit of distance when travelling at speed $\delta \in S$ \\
  $M1^{p,b}$ & Big-M value, $M1^{p,b} = e^{p,b}$ \\
  $M2^{p,b}_i$ & Big-M value, $M2^{p,b}_i = e^{p,b} - h^{p,b}_i$ \\
  $F_c$ & Fuel consumption cost in \$ per ton \\
  $H_c$ & Handling activities cost in \$ per hour\\
  {$I_c$} & Idleness cost in \$ per hour\\
  {$D_c$} & Delay cost in \$ per hour\\
  \hline
\end{tabular}
\end{table}
\begin{table}[th!]
    \centering
\begin{tabular}{ p{3.0cm}  p{12.0cm} }
  \hline
  \multicolumn{2}{c}{Decision variables} \\
  \hline
  $y^{p,b}_{i,j} \in \mathbb{B}$ & 1 if ship $j$ immediately succeeds ship $i$ at berth $b \in B_p$ at port $p \in P$ where $(i,j) \in V^{p,b}$; 0 otherwise \\
  {$v_{i,\delta}^{p} \in \mathbb{B}$} & 1 if ship $i \in N$ sails from port $p$ to some other port $p' (p,p' \in P_i:= p \prec p')$ at speed $\delta \in S$; 0 otherwise \\
  $T^{p,b}_{i} \in \mathbb{Z}^+$ & Time at which ship $i \in N$ berths at berth $b \in B_p$ at port $p \in P_i$ (berthing time) \\
  $T^{p,b}_{o(p,b)} \in \mathbb{Z}^+$ & Time at which berth $b \in B_p$ at port $p \in P_i$ starts berthing ships (i.e., arrival time of the first ship to the berth) \\
  $T^{p,b}_{d(p,b)} \in \mathbb{Z}^+$ & Time at which berth $b \in B_p$ at port $p \in P_i$ finishes berthing ships (i.e., departure time of the last ship from the berth) \\
  $T^{p}_{i} \in \mathbb{Z}^+$ & Time at which port $p \in P_i$ opens activities for ship $i \in N$ \\
  $\Delta EFT^{p}_{i} \in \mathbb{Z}^+$ & Difference between effective finishing time and $EFT^p_i$ for ship $i \in N$ at port $p \in P_i$ \\
  \hline
\end{tabular}
\end{table}

The mathematical model is presented below:
\begin{equation}\label{mip:obj}
\begin{split}
    \min \sum_{i \in N} \sum_{p, p^{\prime} \in P_{i} \backslash\left\{P_{i L}\right\}:\left\{p \prec p'\right\}}  I_{c} \left( \sum_{b \in B_{p'}}{T^{p', b}_i} - \sum_{b \in B_{p}}T^{p, b}_i + \sum_{b \in B_{p}}h_{i}^{p, b}(\sum_{j \in N \cup\{d(p, b)\}}y_{i, j}^{p, b}) - \sum_{\delta \in S} {\Delta_{\delta}} d^{p, p'} v^{p}_{i, \delta} \right)\\ 
    +\sum_{i \in N} \sum_{p \in P_{i}} \sum_{b \in B_{p}} H_{c}(h_{i}^{p, b} \sum_{j \in N \cup\{d(p, b)\}} y_{i, j}^{p, b}) +\sum_{i \in N} \sum_{p \in P_{i}} D_{c} \Delta E F T_{i}^{p}+\sum_{i \in N} \sum_{p, p^{\prime} \in P_{i} \backslash\left\{P_{i L}\right\}:\left\{p \prec p'\right\}} \sum_{\delta \in S} F_{c}(\gamma_{i, \delta} d^{p, p'} v^{p}_{i, \delta})
\end{split}
\end{equation}
subject to:

\begin{align}
    \sum_{b \in B_{p} } \sum_{j \in N\cup\{ d(p, b)\}} y_{i, j}^{p, b} &=1 \quad \forall i \in N, \forall p \in P_{i}\label{mip:con1berth} \\
     \sum_{j \in N \cup\{d(p, b)\}} y_{o(p, b), j}^{p, b} &=1 \quad \forall p \in P, \forall b \in B_{p}\label{mip:con1orig} \\ 
    \sum_{j \in N \cup\{o(p, b)\}} y_{j, d(p, b)}^{p, b} &=1 \quad \forall p \in P, \forall b \in B_{p}\label{mip:con1dest} \\
    \sum_{j \in N \cup\{d(p, b)\}} y_{i, j}^{p, b}-\sum_{j \in N \cup\{o(p, b)\}} y_{j, i}^{p, b} &=0 \quad \forall i \in N, \forall p \in P_{i}, \forall b \in B_{p}\label{mip:conFlow} \\
    T_{i}^{p, b}+h_{i}^{p, b} - \left(1-y_{i, j}^{p, b}\right) M 1^{p, b} &\leqslant T_{j}^{p, b} \quad \forall(i, j) \in A^{p, b}, \forall p \in\left\{P_{i} \cap P_{j}\right\}, \forall b \in B_{p}\label{mip:conH} \\
    \sum_{b \in B_{p}} T_{i}^{p, b}+\sum_{b \in B_{p}} h_{i}^{p, b}(\sum_{j \in N \cup \{d(p,b) \}} y_{i, j}^{p, b})+\sum_{\delta \in S} \Delta_{\delta} d^{p, p'} v^{p}_{i, \delta} &\leqslant T_{i}^{p'} \quad \forall i \in N, \forall p, p^{\prime} \in P_{i} \backslash\left\{P_{i L}\right\}:\left\{p \prec p^{\prime}\right\}\label{mip:conTT} \\
    T_{i}^{p} &\geqslant {Start_{i}^{p}} \quad \forall i \in N, \forall p \in P_{i}\label{mip:conStart} \\
     \sum_{b \in B_{p}}T_i^{p, b}+\sum_{b \in B_{p}} h_{i}^{p, b}(\sum_{j \in N\cup \{d(p, b)\}} y_{i, j}^{p, b})-E F T_{i}^{p} &\leqslant \Delta E F T_{i}^{p}  \quad \forall i \in N, \forall p \in P_{i}\label{mip:conDelay} \\
    \sum_{b \in B_{p}} T_{i}^{p, b} &\geqslant T_{i}^{p} \quad \forall i \in N, \forall p \in P_{i}\label{mip:conBerth} \\
   (\sum_{j \in N\cup \{d(p, b)\}} y_{i, j}^{p, b}+\sum_{j \in N \cup(o(p, b))} y_{j, i}^{p, b}) M 2_i^{p, b} &\geqslant T_{i}^{p, b}  \quad \forall i \in N, \forall p \in P_{i}, \forall b \in B_{p}\label{mip:conRestBerth} \\
    T_{o(p, b)}^{p, b} &\geqslant s^{p, b} \quad \forall p \in P, \forall b \in B_{p}\label{mip:conBo} \\
    T_{d(p, b)}^{p, b} &\leqslant e^{p, b} \quad \forall p \in P, \forall b \in B_{p}\label{mip:conBd} \\
    \sum_{\delta \in S} v^{p}_{i, \delta} &=1 \quad \forall i \in N, \forall p \in P_{i} \backslash\left\{P_{i L}\right\}\label{mip:conSpeed} \\
    y_{i, j}^{p, b} &\in\{0,1\} \quad \forall(i, j) \in A^{p, b}, \forall p \in P, \forall b \in B_{p}\label{mip:conX} \\
    v^{p}_{i, \delta} &\in\{0,1\} \quad \forall i \in N, \forall p \in P_{i}, \forall \delta \in S\label{mip:conV} \\
   \Delta E F T_{i}^{p}, T_{i}^{p} &{\in \mathbb{Z}^+} \quad \forall i \in N, \forall p \in P_{i}\label{mip:conD} \\
    T_{o(p, b)}^{p, b}, T_{d(p, b)}^{p, b} &{\in \mathbb{Z}^+} \quad \forall p \in P, \forall b \in B_{p}\label{mip:conT1} \\
    T_{i}^{p, b} &{\in \mathbb{Z}^+} \quad \forall i \in N, \forall p \in P_{i}, \forall b \in B_{p}\label{mip:conT2}
\end{align}

The objective function (\ref{mip:obj}) minimizes the cost, both for the terminal operators and the liner shipping company. It consists of the four cost elements {described in Section \ref{sec:costs}, namely, the cost of waiting at the port, the vessels' handling cost, the cost of delays, and the total cost of the fuel consumed when sailing between ports. The waiting time is computed as the positive difference between the berthing time and the arrival time whereas the delay is computed as the positive difference between the actual and expected berthing finish time.}
Constraints (\ref{mip:con1berth}) ensure that each ship berths at only one berth at each port in its route. 
Constraints (\ref{mip:con1orig}) and (\ref{mip:con1dest}) denote that at each berth and each port, only one arc leaves the origin and one arrives at the destination respectively.
The flow conservation for all arcs at each berth and each port is ensured by constraint (\ref{mip:conFlow}).
Constraints (\ref{mip:conH}) guarantee that if ship $j$ is berthing right after ship $i$, it waits until the handling is completed. The big-M values for these constraints can be tightened to the time when the berth closes ($e^{p,b}$).
Constraints (\ref{mip:conTT}) ensure for each ship that the activities at the next port in the route do not commence before the ship arrives to the port. The left-hand side of the constraint computes the arrival time to the next port travelling at a chosen speed.
The start of activities for each ship at each port must also start after the minimum allowed time ($Start^p_i$) as indicated in constraints (\ref{mip:conStart}). This also ensures that a ship cannot start berthing if it arrives too early. Both constraints (\ref{mip:conTT}) and (\ref{mip:conStart}) set a lower bound (LB) for the variable $T^p_i$.
Constraints (\ref{mip:conDelay}) compute and set the delay ($\Delta E F T_{i}^{p}$) for each ship at each port.
Constraints (\ref{mip:conBerth}) ensure that the berthing time of each ship at each port starts after the activities for that ship are open at the port.
The values of the berthing time variables for the not chosen berths are set to 0 by constraints (\ref{mip:conRestBerth}).
Constraints (\ref{mip:conBo}) and (\ref{mip:conBd}) ensure that all berthing periods occur within the time window of each berth.
Constraints (\ref{mip:conSpeed}) ensure that exactly one speed is selected to travel between each
pair of consecutive ports (leg) in the route.
The domains for all the decision variables are defined in (\ref{mip:conX})-(\ref{mip:conT2}).
{We notice that a formulation where the time-based variables are defined as non-negative real numbers (i.e., $\mathbb{R}^+$) is also valid. However, we maintain the integer property of the variables for a fair comparison with the presented methods and the formulation presented in \cite{venturini2017a}.}

This formulation contains a few modifications to the original model presented in \cite{venturini2017a} (referred to as \textit{original} model). 
In the original model a set of additional variables for the arrival of a ship to a port is stated. These variables have been omitted in this formulation since the arrival time of a ship to the next port in the route is directly dependent on the departure time from the previous port and the sailing speed between ports. This calculation is given by the left-hand side of constraints ($\ref{mip:conTT}$), which then can be used to replace arrival time variables (e.g., in the objective function).
The delay calculation constraints ($\ref{mip:conDelay}$) use the berthing time ($T^{p,b}_i$) instead of the port opening time for the ship ($T^p_i$).
The big-M value of constraints (\ref{mip:conH}) is set to the closing time of the berth ($e^{p,b}$) instead of $e^{p,b} - min_{c \in (i,j)}\{Start^p_c\}$.

\cite{venturini2017a} enhance the original formulation by adding multiple sets of valid inequalities. These enhancements have also been implemented for the computational comparison. The reader is referred to the original publication for additional details.

\subsection{Network flow formulation}\label{sec:netForm}

\begin{figure}[th]
    \centering
    \includegraphics[width=\textwidth]{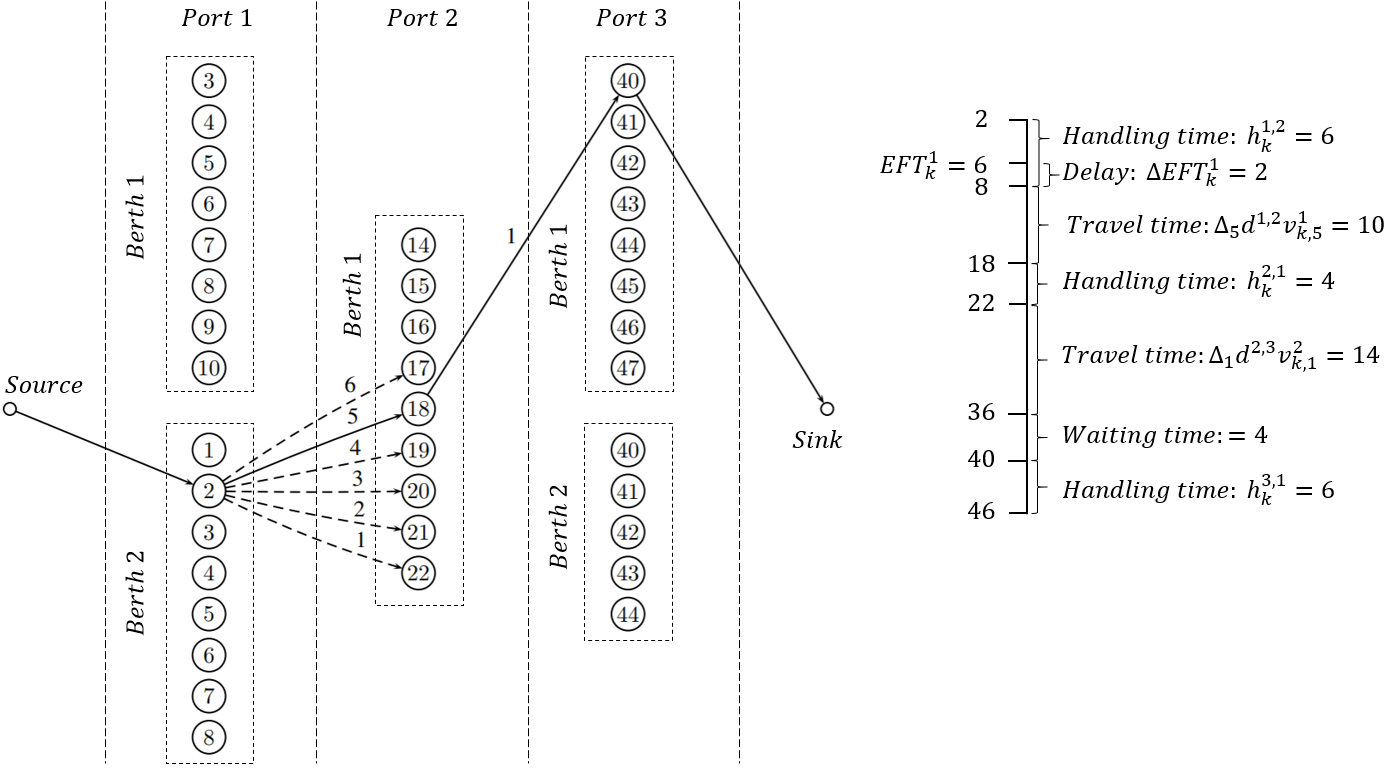}
    \caption{{Example voyage for ship $k$ and corresponding timeline. The number in the nodes indicate the berthing time and the number on the arcs denotes the speed level chosen. Alternative sailing options are denoted with dashed arcs. The rest of arcs in the graph are not displayed for simplicity.}}
    \label{fig:graphPPex}
\end{figure}

The MPBAP can also be modeled as a network flow problem
using a graph representation where each node represents a feasible berthing time at each port and berth and arcs enable the different combinations of berthing times along the route. This setup allows us to obtain a feasible voyage for a given ship by choosing a path along the ports in the graph. Figure \ref{fig:graphPPex} shows an illustrative example of such a path. It consists of three ports with either one or two berthing positions in each of them.    

\begin{figure}[]
    \centering
    \includegraphics[width=0.35\textwidth]{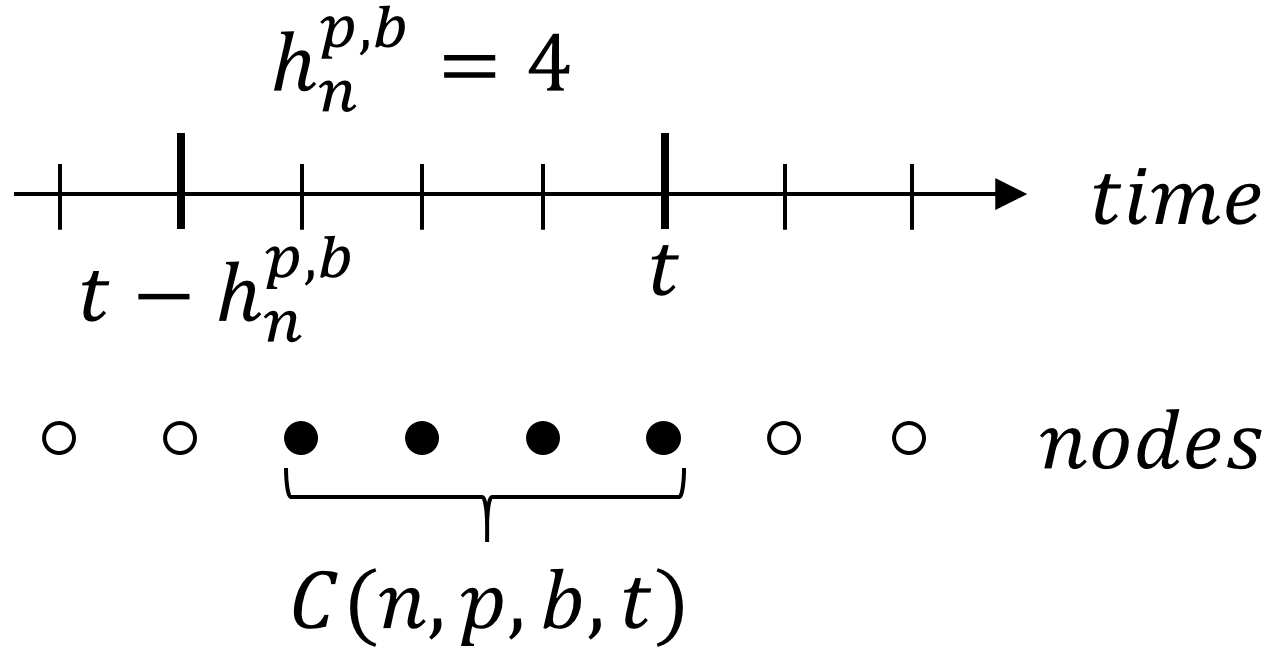}
    \caption{An example of the set $C(n,p,b,t)$ where the nodes depicted belong to $V(p,b)$ and refer to the time instant directly above. $h_n^{p,b}$ denotes the handling time for ship $n$.}
    \label{fig:confT}
\end{figure}

Let $G = (O,A)$ be a directed and acyclic graph formed by the sets of nodes $O$ and arcs $A$. Additionaly, we define the subset of arcs $A^k \subseteq A$ which denote the arcs available for a given ship $k\in N$.
Within the node set, we denote $o,d \in O$ as artificial source and sink nodes respectively. Let $\delta^+_k(u)$ be the set of nodes
that can be reached by following a single outgoing arc $a\in A^k$ from node $u \in O$ for ship $k \in N$.
Likewise, let $\delta^-_k(u)$ be the set of nodes that can be reached by following a single incoming arc $a\in A^k$ from node $u \in O$ for ship $k \in N$. Additionally, $\theta(u)$ denote the berthing time related to node $u\in O\backslash\{o,d\}$ and let $V(p,b) \subseteq O$ be the set of nodes corresponding to port $p\in P$ and berth $b \in B_p$.
We use the notation $[x;y]$ to define an interval between $x$ and $y$ where $y$ is included and $[x;y)$ where $y$ is not.
For each ship $n\in N$ port $p\in P$ berth $b \in B_p$ and operating time instant $t \in [s^{p,b};e^{p,b})$, we define the set $C(n,p,b,t) \subseteq V(p,b)$ that denote the graph nodes for ship $n$ whose berthing periods cover time $t$ (i.e., nodes that are \textit{in conflict} with any ship berthing at time $t$).
This basically corresponds to the nodes of the previous $h_n^{p,b}-1$ time instants and including the node related to time $t$.
An example is depicted in Figure \ref{fig:confT} and 
the expression can be stated 
as follows:
\begin{equation*}
    C(n,p,b,t):= \bigg\{v \in V(p,b) \Big| \theta(v)\in \Big[ \max \Big( t-h_n^{p,b}+1,s^{p,b} \Big);\min\Big(t,e^{p,b}\Big)\Big] \bigg\}
\end{equation*}    

Finally, let $x^k_{i,j}$ be a binary variable deciding if arc $(i,j) \in A^k$ is selected for ship $k\in N$ and let $c_{i,j}$ be the weight associated to the same arc. 

\begin{align}
    \min & \sum_{k \in N} \sum_{(i,j) \in A^k} c_{i,j} x^k_{i,j}\label{eq:Nobj} \\ 
    \sum_{j \in \delta_k^+(o)} x^k_{o,j} &=1 \quad \forall k \in N\label{eq:No} \\
     \sum_{i \in \delta_k^-(d)} x^k_{i,d} &=1 \quad \forall k \in N\label{eq:Nd} \\
      \sum_{i \in \delta_k^-(j)} x^k_{i,j} - \sum_{i \in \delta_k^+(j)} x^k_{j,i} &=0  \quad \forall j \in O \backslash\{o,d\},  k \in N\label{eq:Nflow} \\
     \sum_{k \in N} \sum_{i \in C(k,p,b,t)} \sum_{j \in \delta^+_k(i)} x^k_{i,j} &\leq 1 \quad \forall p\in P, b\in B_p, t \in [s^{p,b};e^{p,b})\label{eq:Nconf} \\
    x^k_{i,j} &\in \{0,1\} \quad \forall (i,j) \in A, k \in N \label{eq:Nbin}
\end{align}
The objective remains the same, and in this case the objective function (\ref{eq:Nobj}) minimizes the cost of the selected arcs.
Constraints (\ref{eq:No}) and (\ref{eq:Nd}) ensure that, for each ship, only one arc leaves from the source node and arrives to the sink node respectively.
Constraints (\ref{eq:Nflow}) enforce flow conservation ensuring that for each node, except the source and sink ones, there are as many incoming as outgoing arcs.
Constraints (\ref{eq:Nconf}) avoid overlapping of berthing periods in the same position by at most allowing one ship to be berthing at each time instant.
Finally, constraints (\ref{eq:Nbin}) define the binary property of the variable.

\subsection{Generalized Set Partitioning Problem formulation}\label{sec:gsppForm}

It is noted that all constraints of the network flow formulation (\ref{eq:Nobj})-(\ref{eq:Nbin}) except constraint (\ref{eq:Nconf}) are independent between ships. Exploiting the structure of the formulation, we can apply Dantzig-Wolfe decomposition (\cite{dantzig1960a}) and transform it into a generalized set partitioning problem (GSPP) formulation where constraint (\ref{eq:Nconf}) is handled in the master problem and each variable (i.e., column) refers to a whole feasible schedule of a ship along its route.  According to \cite{jans2010a}, the pure binary nature of the variables of the network flow formulation allows us to impose binary conditions on the variables of the new master problem.

{The set of all columns is comprised in $\Omega$ and} the decision variable $\lambda_j$ is set to 1 if column $j \in \Omega$ is chosen as part of the solution and 0 otherwise.
We denote $c_j$ as the cost related to column $j \in \Omega$. In order to replicate the objective of the MIP formulation, this cost consists of the idleness, handling cost, delay and bunker consumption cost of the ship denoted by the column.
Let $A^i_j$ be a parameter that is equal to 1 if column $j \in \Omega$ corresponds to ship $i \in N$ and 0 otherwise.
Likewise, let $Q^{p,b,t}_j$ be a parameter that is equal to 1 if the ship of column $j \in \Omega$ is occupying berth $b \in B_p$ at time instant $t \in [s^{p,b};e^{p,b})$ at port $p \in P$ and 0 otherwise. 
\begin{align}
    \min &\sum_{j \in \Omega} c_{j} \lambda_{j}\label{GSPP:obj} \\
    \sum_{j \in \Omega} A^i_{j} \lambda_{j} &=1 \quad \forall i \in N\label{GSPP:conv} \\
    \sum_{j \in \Omega} Q^{p,b,t}_{j} \lambda_j &\leqslant 1 \quad \forall p \in P, b \in B_p,t \in [s^{p,b};e^{p,b})\label{GSPP:pbt} \\
    \lambda_j &\in \{0,1\} \quad \forall j \in \Omega\label{GSPP:bin}
\end{align}
The objective function (\ref{GSPP:obj}) minimizes the cost $c_j$ of the columns.
Constraints (\ref{GSPP:conv}) ensure that one column is selected for each ship.
Constraints (\ref{GSPP:pbt}) guarantee that, at each time instant, there is at most one ship berthing at each berth of a port.
Finally, constraints (\ref{GSPP:bin}) set the binary property of the decision variables.

\section{Solution method} \label{sec:SolM}

To solve (\ref{GSPP:obj})-(\ref{GSPP:bin}), we propose a solution method based on a column generation procedure that, combined with branching, additional valid inequalities and symmetry breaking methods, results in a \textit{branch-and-cut-and-price} algorithm.

\subsection{Delayed Column generation} \label{sec:CG}
A common way of solving the GSPP formulation is by adding all the columns in advance. A successful example of this approach for the BAP can be found in \cite{buhrkal2011a}. For the BAP instances presented, the amount of columns is manageable and can be easily pre-processed. However, in the MPBAP, the amount of columns increase exponentially with the multiple sailing speeds and ports for a ship. This makes the pre-processing intractable even for a few ports. Therefore, more dynamic strategies for handling the columns need to be explored. One efficient procedure is the so-called \textit{delayed column generation}. This procedure relies on the premise that most of the variables will not be part of the optimal solution and, therefore, have a value of zero. Then, the focus is only on generating columns that have the potential to improve the objective value. This is done by relaxing and splitting the main problem into two, the master and subproblem. The restricted master problem (RMP) is the linear relaxation of the original formulation containing only a subset of the variables. The subproblem (or pricing problem) is used to identify the new variables. In our case, the relaxed version of the GSPP becomes the RMP and we define $N$ independent subproblems, one per ship. The subproblem is defined as a shortest path problem in the network defined in Section \ref{sec:netForm} which can be solved in polynomial time. Since the graph is directed and acyclic (DAG), it can be solved by a DAG shortest path algorithm (see \cite{cormen1996a} or \cite{magnanti1993a}).
The pricing problem aims at minimizing the reduced cost of a given path.
At each iteration, after solving the RMP, the dual values of the RMP constraints are used to solve the pricing problems. We denote $\alpha_k$ to the dual variable for ship $k \in N$ associated to constraint (\ref{GSPP:conv}). Likewise, we denote $\mu_{p,b,t}$ to the dual variable for port $p\in P$, berth $b\in B_p$ and time $t\in [s^{p,b};e^{p,b})$ associated to constraint (\ref{GSPP:pbt}). Let $\Bar{\alpha}_k, \Bar{\mu}_{p,b,t}$ be the dual solution values for the RMP and let $\Lambda_j$ be a sequence of \textit{(port,berth,time)} elements. Each of these elements refers to the port, berth and time of a graph node visited by column $j\in \Omega$.
The reduced cost $\hat{c}_j$ for a specific path $j$ for ship $k\in N$ is computed as follows:
\begin{equation*}
    \hat{c}_j = c_j - (\sum_{(p,b,t) \in \Lambda_j} \sum_{t' \in [t;t+h^{p,b}_k)} \Bar{\mu}_{p,b,t'}) - \Bar{\alpha}_k
\end{equation*}
Finally, for each pricing problem, we add the path with the lowest reduced cost to the RMP only if $\hat{c}_j$ is negative (i.e, $\hat{c}_j <0$).

In fact, when the pricing problem is a pure shortest path problem, the LP bound arising from solving the GSPP with column generation and solving the LP relaxation of the network flow problem is the same. This indicates that the Dantzig-Wolfe decomposition does not provide any gain bound-wise. On the other hand, in cases of very dense networks with significantly more arcs than nodes as in our case, solving the GSPP with column generation is expected to be faster (e.g., see \cite{brouer2011a}). 

\subsection{Branching}
Since the decision variables of the RMP are linear, the solution at the root node is often fractional and branching methods are required in order to achieve integrality.
A major aspect of the branching procedure is selecting a branching candidate, whose branch children improve the lower bound the most.
The most common branching methods consider branching on a specific node or arc from the graph. These strategies can be effective in some cases but do not necessarily apply to our problem. For instance, when branching on a graph node, one child will enforce the graph node to be used in the subsequent branch-and-bound (B\&B) tree while the other child will forbid it. Considering the large amount of nodes for most instances in this problem, we can clearly see that the effect can be significant for the first child but rather minimal for the second. This often results in a highly unbalanced B\&B tree to explore.
In this study, we present a different branching strategy for the problem at hand that aims to be more effective than branching on a single graph node.

The proposed branching strategy states that, given a fractional solution, we compute, for each ship $n$ and port $p$, the average berthing 
time $t$ and the variance of these times 
among all solution columns. 
{As an example, consider a fractional solution containing two columns for ship 1. At port 1, these columns correspond to ship 1 berthing at 
time 4 and 6 
respectively. Then, for ship 1 and port 1, the average berthing 
time is $5$ 
whereas the sample standard deviation is $\sqrt{2}$.}
We define this average berthing
time
and variance as a \textit{candidate} which results in a total of $|N|\cdot|P|$ candidates. We then select the candidate whose variance of berthing 
times
is higher. The procedure is described in Algorithm \ref{alg:branchNode}. Each of the child branches will enforce ship $n$ to berth 
before or after time $t$ respectively
 at port $p$. It should be noticed that a fractional solution where ships berth at the same 
 time
 but at different berthing 
 positions
 can exist. In this {case}, we can obtain candidates with no variance resulting in an impractical branching. If that happens, the criterion is changed to branching on berthing
 positions
 instead of on berthing 
 times,
following the same procedure. In practice, this scenario is highly unlikely to happen and we have not experienced it in any experiments hitherto.
As a result, the proposed strategy opts for branching on a set of graph nodes instead of on a single one.

Finally, the B\&B tree is explored following a \textit{best first} policy.
This policy prioritizes the queue of unexplored nodes according to their bound. Thus, the next node to be explored is always the one with the \textit{best} (i.e., lowest) lower bound.

\begin{algorithm}
\DontPrintSemicolon
\KwData{$sol$: current solution.}
\KwResult{$Cand^*$: the candidate to branch on.}
\Begin{
$[\lambda] \gets sol$ \tcp*{classify solution columns ($\lambda$) by ship}
$Cand^* = \emptyset$ \tcp*{initialize best candidate}
$\sigma^* = 0$ \tcp*{initialize standard deviation of candidate's berthing times}
\For{\textit{ships} and \textit{ports}}{
    $[times] \gets \lambda(ship,port)$ \tcp*{set of solution berthing times at \textit{port} for \textit{ship}}
    $time \gets avg([times])$ \tcp*{get average of berthing times}
    \If(\tcp*[f]{compare the standard deviation with the current best}){$\sigma([times]) > \sigma^*$}{
        $\sigma^* \gets \sigma([times])$\;
        $Cand^* \gets time,port,ship$ \tcp*{update best candidate so far}
    }
}
}
\caption{Branching candidate selection} \label{alg:branchNode}
\end{algorithm}

\subsection{Valid inequalities}

In order to improve the lower bound, we propose a set of valid inequalities that can be added to the problem by separation. 

\begin{figure}[]
    \centering
    \includegraphics[width=0.8\textwidth]{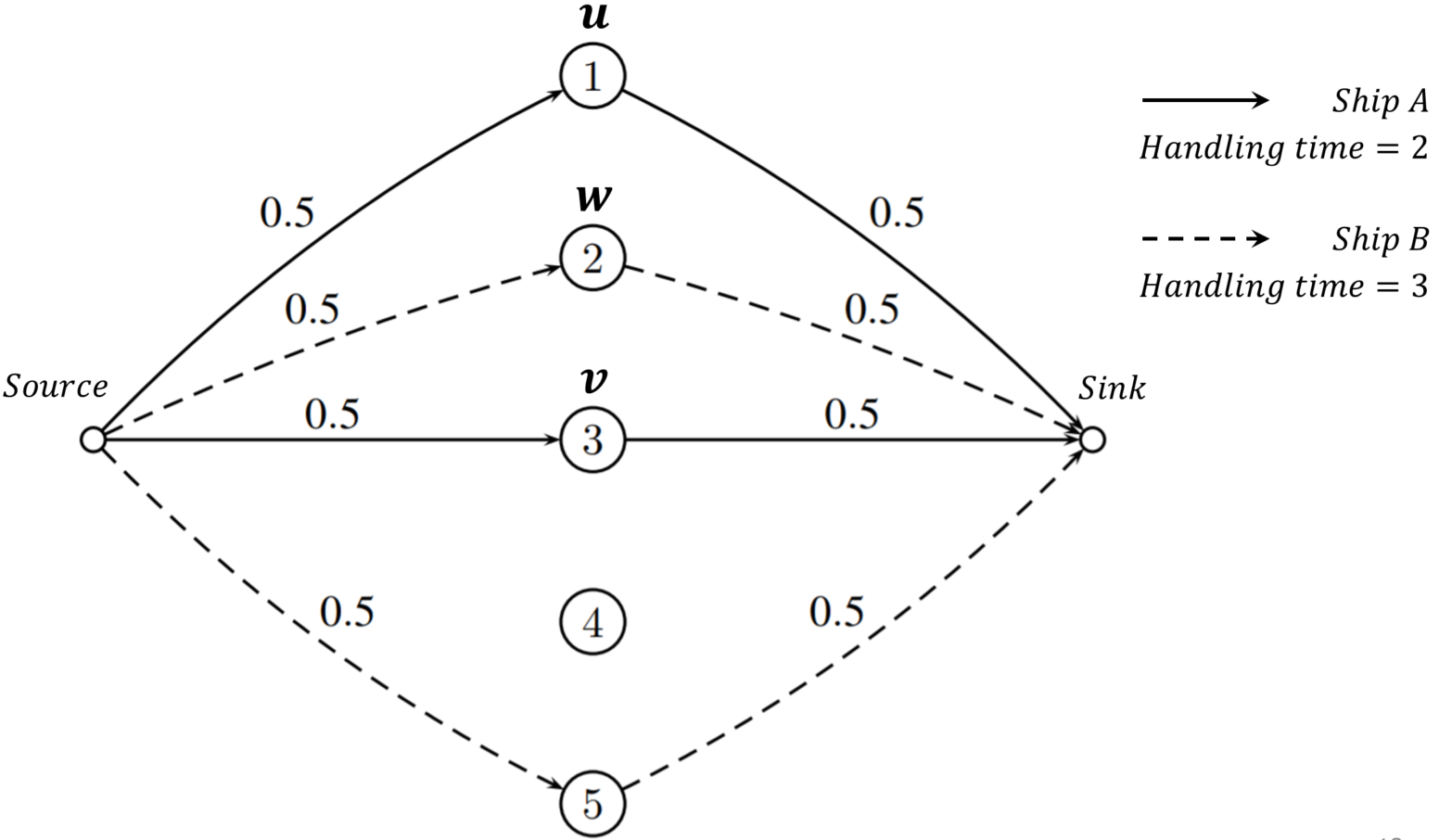}
    \caption{{Example LP solution of a problem with two ships, one port and one berth. The nodes represent berthing times and the numbers on the arcs denote the solution value of the arc variable $x^k_{i,j}$}. 
    }
    \label{fig:3nCut}
\end{figure}
Figure \ref{fig:3nCut} shows a small LP solution to a trivial problem 
{with two ships (i.e., continuous and dashed lines), one port and one berth}
where an example of a violated valid inequality can be found. We define $u,v$ as the two nodes corresponding to ship A (berthing at times 1 and 3) and let $w$ be the node of ship B berthing at time 2.
We observe that the arc from node $w$ is in conflict with the arcs from both nodes $u$ and $v$ due to overlapping berthing periods. In other words, the berthing period of ship B at node $w$ covers, at least partially, both berthing periods of ship A at nodes $u$ and $v$. The arcs from $u,v$ are also in conflict with each other as they belong to the same ship.
As a result, we notice that, at most, one outgoing arc can be chosen out of the ones from these three nodes. Since the solution values of the outgoing arcs sum to 1.5, this valid inequality would cut the example LP solution.
We aim at generalizing the definition of such a valid inequality and introduce the following proposition:

\begin{proposition}\label{prop:1}
Given two time instants $t_1,t_2 \in [s^{p,b};e^{p,b})$ where $t_1 < t_2$ and a port $p\in P$, berth $b\in B_p$ and ship $n\in N$, the following is a valid inequality:
\[
\sum_{ u \in \bigcup_{t\in [t_1;t_2]} C(n,p,b,t)}
\sum_{w \in \delta^+_n(u)} x^n_{u,w} + 
\sum_{m\in N \backslash\{n\}}\sum_{u \in C(m,p,b,t_1)\cap C(m,p,b,t_2)}\sum_{w \in \delta^+_m(u)} x^m_{u,w} \leqslant 1 
\]
\end{proposition}

\proof{Proof.}
The set $C(m,p,b,t)$ used in constraint (\ref{eq:Nconf}) defines the set of nodes for ship $m$ that are in conflict with time $t$ (see Section \ref{sec:netForm}). Based on this definition, the intersection set $C(m,p,b,t_1) \cap C(m,p,b,t_2)$ directly defines the set of nodes for ship $m$ that are in conflict with both time instants $t_1$ and $t_2$. Constraint (\ref{eq:Nconf}) indicates that at most one arc can be chosen out of the nodes from the sets $C(m,p,b,t)$ of all ships $m \in N$ and, therefore, the same applies to the intersection set $C(m,p,b,t_1) \cap C(m,p,b,t_2)$. 
By considering the intersection set $C(m,p,b,t_1) \cap C(m,p,b,t_2)$ for all ships except one $m\in N\backslash\{n\}$, the berthing period for ship $n$ is only required to be in conflict with either $t_1$ or $t_2$ and can be defined as the union of $C(n,p,b,t_1)\cup C(n,p,b,t_2)$. Considering these node sets, we can define the following valid inequality:
\begin{multline*}
    \sum_{u \in C(n,p,b,t_1)\cup C(n,p,b,t_2)}\sum_{w \in \delta^+_n(u)} x^n_{u,w} + 
\sum_{m\in N \backslash\{n\}}\sum_{u \in C(m,p,b,t_1)\cap C(m,p,b,t_2)}\sum_{w \in \delta^+_m(u)} x^m_{u,w} \leqslant 1  \\
\forall p\in P, b\in B_p, n\in N,t_1,t_2 \in [s^{p,b};e^{p,b}),t_1 < t_2 
\end{multline*}
Based on the assumption that a berthing period cannot be discontinued, the intersection set $C(m,p,b,t_1) \cap C(m,p,b,t_2)$ for any ship is not only in conflict with times $t_1$ and $t_2$ but with all the time instants in the period $[t_1;t_2]$. Therefore the interval for ship $n$ can be expanded to the union of $C(n,p,b,t)$ sets
for all time instants $t \in [t_1;t_2]$. An example of this set is shown in Figure \ref{fig:IntervShipN} and the resulting valid inequality can be formulated as follows:
\begin{figure}[]
    \centering
    \includegraphics[width=0.6\textwidth]{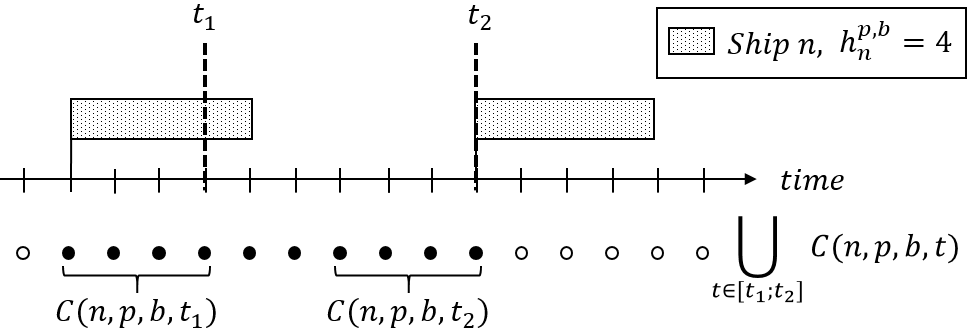}
    \caption{
    {In a valid inequality for ship $n$, port $p$, berth $b$ and time instants $t_1,t_2$, the filled nodes indicate the interval for ship $n$ with handling time $h_n^{p,b}$. The rectangles indicate the berthing period of ship $n$ at the earliest and latest possible berthing times in the interval.}
    }
    \label{fig:IntervShipN}
\end{figure}

\begin{multline}\label{eq:strongVI}
\sum_{u \in \bigcup_{t\in [t_1;t_2]} C(n,p,b,t)}\sum_{w \in \delta^+_n(u)} x^n_{u,w} + 
\sum_{m\in N \backslash\{n\}}\sum_{u \in C(m,p,b,t_1)\cap C(m,p,b,t_2)}\sum_{w \in \delta^+_m(u)} x^m_{u,w} \leqslant 1  \\
\forall p\in P, b\in B_p, n\in N,t_1,t_2 \in [s^{p,b};e^{p,b}),t_1 < t_2 
\end{multline}
\Halmos
\endproof

Returning to the example in Figure \ref{fig:3nCut}, the mentioned cut would be included in the proposed valid inequality (\ref{eq:strongVI}) for $n=A$, $t_1=2$ and $t_2=4$ where node $w$ would correspond to a node from the intersection sets $C(B,p,b,2) \cap C(B,p,b,4)$ and nodes $u,v$ for ship $A$ would correspond to berthing times covering $t_1$ and $t_2$ respectively and therefore belonging to the set $\bigcup_{t \in [2;4]} C(A,p,b,t)$.

We note that the inequality only is interesting when $C(m,p,b,t_1)\cap C(m,p,b,t_2)\neq\emptyset$. The size of the intersection set is dependent on the time instants $t_1,t_2$ used and we observe that this size increases when the $t_1,t_2$ are closer together in time.

These valid inequalities (\ref{eq:strongVI}) are added by separation after the column generation procedure concludes. Exploring the entire set of valid inequalities can be computationally intensive. Therefore, only valid inequalities based on berthing times from the LP solution are checked since the arcs from the related nodes are guaranteed to contain non-zero values and the resulting inequalities have a higher probability of being violated by the LP solution.
Given an LP solution, let $t^*_1$ and $t^*_2$ be two berthing times for ship $n$ at berth $b$ of port $p$ where $t^*_1 \leq t^*_2$. Let $t^*_3$ be a berthing time for another ship $m$ at the same berth $b$ of port $p$ whose berthing period is both in conflict with $t^*_1$ and $t^*_2$ for ship $n$. The conditions that  $t_3^*$ needs to satisfy to be in conflict with $t_1^*$ and $t_2^*$ are given by the following inequalities:
\begin{equation*}
  t^*_1 + h_n^{p,b} > t^*_3
\end{equation*}
\begin{equation*}
    t^*_2 < t^*_3 +  h^{p,b}_m
\end{equation*}
Based on these times, we can calculate time instants $t_1,t_2$ for a valid inequality that includes $t_1^*,t_2^*$ for ship $n$ and $t_3^*$ for ship $m$ as follows:
\begin{equation*}
    t_1 = t^*_1 + h^{p,b}_n - 1, \quad  t_2 = t^*_2
\end{equation*}
\begin{figure}[]
    \centering
    \includegraphics[width=0.8\textwidth]{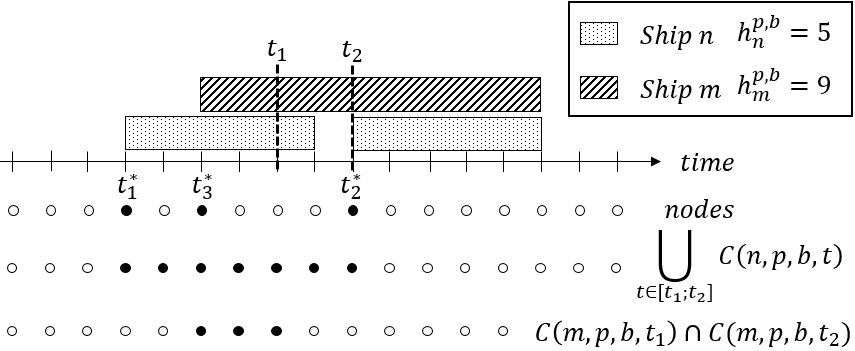}
    \caption{Example of times $t_1,t_2$ definition based on solution times $t^*_1,t^*_2$ for ship $n$ and $t^*_3$ for ship $m$. The bottom last two rows of filled nodes define the node interval for ship $n$ and $m$ respectively.}
    \label{fig:IntervalSolVals}
\end{figure}
An example of this calculation is shown in Figure \ref{fig:IntervalSolVals}.
It can be noticed that the interval for ship $n$ starts at time $t_1^*$ and ends at time $t_2^*$. If we add such a violated cut to the RMP, we risk finding a very similar solution in the next iteration where columns are shifted, for example, one time instant before $t_1^*$ or after $t^*_2$.
In order to avoid that, we aim at defining time instants $t_1,t_2$, so that the resulting intervals do not only cover solution nodes but also a number of neighboring nodes related to time instants immediately before and after the solution time.
We aim at expanding the interval between $t^*_1$ and $t^*_2$ as well as the one around $t^*_3$.
Based on the inequalities aforementioned to ensure that $t_1^*,t_2^*$ and $t_3^*$ relate to conflicting periods, we introduce the slack variables $\Delta^X$ and $\Delta^Y$ that would indicate how much we can modify the node intervals.
\begin{equation*}
  t^*_1 + h_n^{p,b} > t^*_3 + \Delta^X
\end{equation*}
\begin{equation*}
    t^*_2 + \Delta^Y < t^*_3 +  h^{p,b}_m
\end{equation*}
Both slack values are distributed equally between both intervals, which leads us to the following calculation of $t_1$ and $t_2$:
\begin{equation*}
    t_1 = t^*_1 - \frac{\Delta^X}{2} + h^{p,b}_n - 1, \quad  t_2 = t^*_2 + \frac{\Delta^Y}{2} 
\end{equation*}
Due to the discretization of the time horizon, if $\frac{\Delta^X}{2}$ or $\frac{\Delta^Y}{2}$ is fractional, they are rounded-up in the calculation of $t_1$ and $t_2$. Also, in the case that there is limited room for expansion in one of the intervals (e.g., operational time windows), the remaining slack is added to the other interval.

\begin{figure}[]
    \centering
    \includegraphics[width=0.8\textwidth]{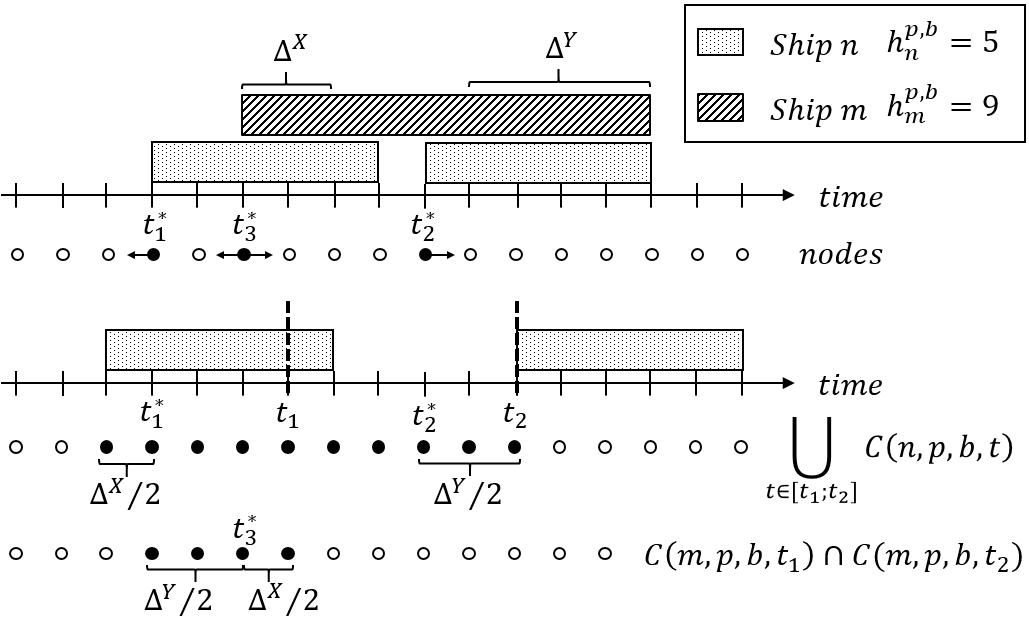}
    \caption{
    {Example of times $t_1,t_2$ selection based on solution times $t^*_1,t^*_2$ for ship $n$ and $t^*_3$ for ship $m$. The upper illustration depicts the LP solution for the three berthing times selected and the available slack and direction of expansion for the desired node intervals. The lower illustration depicts the resulting times $t_1,t_2$ for the valid inequality and the respective node intervals for ships $n$ and $m$.}
    }
    \label{fig:IntervalCalc}
\end{figure}
Figure \ref{fig:IntervalCalc} shows an example of the calculation of times $t_1$ and $t_2$ based on solution times $t^*_1,t^*_2$ and $t^*_3$ and slack variables $\Delta^X, \Delta^Y$.

In order to ensure $t_1 < t_2$, by substituting the above expressions, the criterion that $t_1^*,t_2^*,t_3^*$ need to fulfill in order to result in a valid inequality can be defined as follows:
\begin{equation*}
    t^*_1 - \frac{\Delta^X}{2} + h^{p,b}_n - 1 < t^*_2 + \frac{\Delta^Y}{2}
\end{equation*}
Not satisfying this inequality leads to a cut that, at best, is equal to constraint (\ref{GSPP:pbt}) which is already present in the RMP.

The entire cut separation process is summarized in Algorithm \ref{alg:cutSep}. The procedure requires the RMP model and an LP solution as input. From the solution, both the $\lambda^*$ solution values and the berthing times of the solution columns are extracted and classified by ship, port and berthing position. 
The cuts are checked by enumerating combinations of solution times $t^*_1,t^*_2$ and $t^*_3$. 
Only solution times whose berthing periods are in conflict are considered. This is the case if the berthing period of ship $m$ at time $t_3^*$ overlaps both berthing periods of ship $n$ at times $t_1^*$ and $t_2^*$ ($inConflict(t_1^*,t_2^*,t_3^*)$ in Algorithm \ref{alg:cutSep}). 
Then, the solution times are used to compute time instants $t_1,t_2$ distributing the slack available as aforementioned in this Section ($ t_1,t_2 \gets calcInterval(t_1^*,t_2^*,t_3^*)$ in Algorithm \ref{alg:cutSep}).
To check and add the violated cuts to the RMP, equation (\ref{eq:strongVI}) needs to be translated to the problem variables. 
The $x^k_{i,j}$ variables can be defined using $\lambda_p$ variables as follows:
\begin{equation*}
    x_{i,j}^k = \sum _{p\in \Omega} q^k_{i,j,p} \lambda_p
\end{equation*}
where parameter $q^k_{i,j,p}$ is 1 if graph arc $(i,j)\in A^k$ for ship $k$ is used by column $p \in \Omega$ and 0 otherwise. Applying this equality to equation (\ref{eq:strongVI}), we obtain the following version of the equation:
\begin{equation}\label{eq:VIrmp}
    \sum_{u \in \bigcup_{t\in [t_1;t_2]} C(n,p,b,t)}\sum_{w \in \delta^+_n(u)}\sum _{j\in \Omega} q^n_{u,w,j} \lambda_j
    +\sum_{m \in N\backslash\{n\}}\sum_{u \in C(m,p,b,t_1)\cap C(m,p,b,t_2)}\sum_{w\in\delta_{m}^{+}(u)}\sum _{j\in \Omega} q^m_{u,w,j} \lambda_j
    \leq 1
\end{equation}
For each cut inspected, the left-hand side of constraint (\ref{eq:VIrmp}) is computed and the valid inequality is added to the RMP if it is violated. 

\begin{algorithm}
\DontPrintSemicolon
\KwData{$sol,RMP$: current solution and model.}
\KwResult{$RMP$: updated model with separated cuts.}
\Begin{
$times[p,b,n] \gets sol$ \tcp*{classify solution times by port $p$, berth $b$ and ship $n$}
$[\lambda^*] \gets sol$ \tcp*{obtain solution values for columns}
\For(\tcp*[f]{cuts are based on a specific berth, port and ship}){$p\in P, b \in B_p, n \in N$}{
    \For(\tcp*[f]{loop over pairs of solution times for ship $n$}){$t^*_1,t^*_2 \in times[p,b,n]$}{
        \For{$m \in N,m\neq n$}{
            \For(\tcp*[f]{select a third time from a different ship}){$t^*_3 \in times[p,b,m]$}{
                \If(\tcp*[f]{check if berthing periods are in conflict}){inConflict($t^*_1,t^*_2,t^*_3$)}{
                    $t_1,t_2 \gets$ calcInterval($t^*_1,t^*_2,t^*_3$) \tcp*{compute $t_1,t_2$ for the valid inequality}
                    $violatedCut \gets$ checkCut($t_1,t_2,n,p,b,[\lambda^*]$) \tcp*{add violated cut to the RMP}
                }
            }
        }
    }
}
}
\caption{Cut separation} \label{alg:cutSep}
\end{algorithm}

These valid inequalities are relatively easy to handle in the reduced cost computation. For each valid inequality, its corresponding dual value needs to be subtracted in each of the nodes considered for each ship in the constraint. As an example, given a valid inequality for times $t_1,t_2$ where $t_1 < t_2$, port $p$, berth $b$ and ship $n$, its dual value needs to be subtracted in nodes $ \bigcup_{t \in [t_1;t_2]} C(n,p,b,t)$ for ship $n$ and in nodes $C(m,p,b,t_1) \cap C(m,p,b,t_2)$ for ship $m$ where $m\neq n$. A more mathematical definition of the updated reduced cost computation is given in Appendix \ref{app:a}. 

\subsection{Symmetry breaking}
In some instances, at each port, some of the berthing positions are identical in terms of their availability time window and the handling times for all ships. Identical berths may lead to many equivalent solutions, which may increase the solving time of the model.
Therefore, we propose adapting the model so it deals with berth types instead of individual berths in a similar procedure as the one stated in \cite{buhrkal2011a}.
Let $K_p$ be the set of berth types for port $p\in P$ and $\beta^k$ be the number of berthing positions of type $k\in K$ in the problem. {For each berth type $k \in K_p$ at port $p \in P$, $s^{p,k}$ and $e^{p,k}$ denote its opening and closing time respectively and the parameter $Q^{p,k,t}_{j}$ is 1 if the ship from column $j \in \Omega$ occupies berth type $k$ at time instant $t\in [s^{p,k};e^{p,k})$ at port $p$ and 0 otherwise.} We can therefore update the set of constraints (\ref{GSPP:pbt}) as follows:
\begin{equation}\label{GSPP:pkt}
\sum_{j \in \Omega} Q^{p,k,t}_{j} \lambda_j \leqslant \beta^k \quad \forall p \in P, k \in K_p,t \in [s^{p,k};e^{p,k})
\end{equation}
This adaptation has an equivalent impact in constraints (\ref{eq:Nconf}) from the network formulation where the right-hand side is also replaced by $\beta^k$. The valid inequality (\ref{eq:VIrmp}) from Proposition \ref{prop:1} can be updated similarly and it is described in Proposition \ref{prop:2} that can be found in Appendix \ref{app:b}. The resulting valid inequality is formulated as follows:
\begin{multline}\label{GSPP:VIbeta}
\sum_{u \in \bigcup_{t\in [t_1;t_2]} C(n,p,k,t)}\sum_{w \in \delta^+_n(u)}\sum _{j\in \Omega} q^n_{u,w,j} \lambda_j
    +\sum_{m \in N\backslash\{n\}}\sum_{u \in C(m,p,k,t_1)\cap C(m,p,k,t_2)}\sum_{w\in\delta_{m}^{+}(u)}\sum _{j\in \Omega} q^m_{u,w,j} \lambda_j
    \leq \beta^k \\
    \forall p\in P, k\in K_p, n\in N,t_1,t_2 \in [s^{p,k};e^{p,k}),t_1 < t_2 
\end{multline}
The reduced cost computation is also slightly modified where the dual variable $\mu_{p,k,t}$ of the modified constraint now is based on berth type $k \in K_p$ instead of berth $b\in B_p$.

We expect to see an improvement in the computational time as soon as there are two identical berths at a port. Likewise, we expect to see larger symmetry for the instances containing more berthing positions per port.

\section{Cooperative game theory}\label{sec:cgt}

The MPBAP is based on a strong collaboration between carriers and port operators and some of them, especially carriers, may be reticent to take part in such a collaboration scheme. In order to convince them that this form of collaboration is beneficial for all of them, we define a cooperative game.
The aim is to show that all stakeholders (i.e., carriers and terminals) can potentially benefit from a collaboration by distributing the overall {costs} efficiently. 
Our cooperative game is formed by a set of players {$\mathcal{P} = \{1,...,p\}$}, which in this case corresponds
{to both the carriers owning the ships and the terminal operators of the ports visited by the ships.}
The \textit{characteristic function} $\vartheta(S)$ measures the impact of a coalition of players $\mathcal{S}\subseteq \mathcal{P}$, which in this case is measured by the operational costs. The coalition formed by all players is known as the \textit{grand coalition}. 
It is normally assumed that the characteristic function satisfies:
\begin{align}
    \vartheta(\emptyset) &=0 \label{eq:void}\\
    \vartheta(\mathcal{S} \cup \mathcal{T}) &\leq \vartheta(\mathcal{S})+\vartheta(\mathcal{T}) \quad \forall \mathcal{S}, \mathcal{T} \subseteq \mathcal{P}, \quad \mathcal{S} \cap \mathcal{T}=\emptyset \label{eq:supera}
\end{align}
Equation (\ref{eq:void}) states that an empty coalition has a cost of zero, while equation (\ref{eq:supera}), known as {\textit{subadditivity}}, indicates that the costs of two separate coalitions $\mathcal{S},\mathcal{T} \subseteq \mathcal{P}$ cannot be lower than when acting together.
A solution to a cooperative game (i.e., \textit{imputation}) can be defined as $f = \{f_1,...,f_p\}$ where $f_i$ corresponds to the cost allocation of player $i$ in coalition $\mathcal{P}$. An imputation should satisfy the following conditions:
\begin{align}
f_{i} & \leq \vartheta(\{i\}) \quad \forall i \in \mathcal{P} \label{eq:proRat} \\
\sum_{i \in \mathcal{P}} f_{i} &=\vartheta(\mathcal{P})\label{eq:proEff}
\end{align}
The first condition is based on individual rationality and defines that the cost allocation for a player when being part of the grand coalition cannot be worse than the player's standalone cost. 
The second condition is based on group rationality and states that all the savings arising from a grand coalition are shared. This is the equivalent of saying that the sum of cost allocations needs to be equal to the total cost of the grand coalition and a solution fulfilling this condition is said to be \textit{efficient}.
Furthermore, we consider a solution to be \textit{stable}, if, for every coalition $\mathcal{S} \subseteq \mathcal{P}$, the sum of allocated cost of the players of the coalition is not higher than the cost of the coalition $\sum_{k \in \mathcal{S}} f_{k} \leq \vartheta(\mathcal{S})$.
We define the \textit{core} as the set of solutions that are both \textit{efficient} and \textit{stable}.
We see the core solutions as the 
most attractive and fair for all players.
Note, however, that the core may be empty in some cases.
{This means that a cost allocation that satisfies both the efficiency and stability properties does not exist. In other words, it means that a subset of the players in the grand coalition could do better by themselves (i.e., by forming a sub-coalition). If the core is empty, the grand coalition is unstable and there is a risk that it breaks apart. In practice, the grand coalition may stay together despite a non-core solution. For instance, it may be that a subset of players are not aware of the higher benefits of a specific sub-coalition or that players choose to stay in the coalition to reap more long-term benefits given future expectations.}
Next, we describe the two allocation methods we have used in this study.   

\subsection{Shapley value}
The Shapley value \citep{shapley1953value}, refers to the weighted average of each player's marginal contribution to each of the potential coalitions. Let $\Theta^i(\mathcal{S})$ be the marginal contribution of player $i$ to coalition $\mathcal{S}$, which is seen as the difference between the cost of the coalition including player $i$ and the coalition without the player:
\begin{equation}
    \Theta^i(\mathcal{S}) = \vartheta(\mathcal{S} \cup \{i\}) - \vartheta(\mathcal{S})
\end{equation}
Then, the cost allocated to participant $i$ is computed by the following expression:
\begin{equation}f_{i}=\sum_{\mathcal{S} \subseteq \mathcal{P} \backslash\{i\}} \frac{|\mathcal{S}| !|\mathcal{P} \backslash(\mathcal{S} \cup\{i\})| !}{|\mathcal{P}| !} \Theta^{i}(\mathcal{S})\end{equation}
where $| \cdot |$ refers to the number of players in the given coalition. Once the characteristic function $\vartheta(\mathcal{S})$ is calculated for all possible coalitions $\mathcal{S}$, it is a simple method to compute as it only requires applying a formula.
The Shapley value does not only provide \textit{efficient} solutions, it also contains other valuable properties. The solutions are \textit{symmetric} meaning that if two players contribute equally to the coalitions, they achieve the same savings. \textit{Anonimity} is also ensured, which states that the order or labelling of players does not have an impact on the assignment of savings. This property ensures a unique solution which avoids players to regret their choices and prevents additional negotiation processes. 
On the other hand, the Shapley value does not ensure the \textit{stability} property, meaning that the solution is not guaranteed to be part of the core.

\subsection{Equal profit method (EPM)}
The goal of the equal profit method \citep{frisk2010a} is to find the solution in the core that minimizes the maximal difference in relative savings between pairs of players. The relative saving of player $i$ is computed as $\frac{\vartheta(\{i\}) - f_i}{\vartheta(\{i\})}$. The method is formulated as the following linear programming model:
\begin{align}
     \min & \text{ }z \\
     z & \geq \frac{f_{i}}{\vartheta(\{i\})}-\frac{f_{j}}{\vartheta(\{j\})} \quad \forall i, j \in \mathcal{P} \label{eq:relSav}\\
     \sum_{i \in \mathcal{P}} f_{i} &=\vartheta(\mathcal{P}) \label{eq:eff}\\
     \sum_{i \in \mathcal{S}} f_{i} &\leq \vartheta(\mathcal{S}) \quad \forall \mathcal{S}\subseteq \mathcal{P} \label{eq:stable} \\
     f_{i} &\geq 0 \quad \forall i \in \mathcal{P}
\end{align}
Constraints (\ref{eq:relSav}) calculate the difference in relative savings between each pair of players and {restricts} $z$ to the largest of those differences.
Note that constraints (\ref{eq:eff}) and (\ref{eq:stable}) are the ones denoting the stability and efficiency properties which means that the EPM method only allows solutions lying in the core. 

\section{Computational results}\label{sec:results}
This section is divided in two. First, the performance of the proposed method is compared to a commercial solver on the set of instances from \cite{venturini2017a} and an additional generated set of harder instances. The second part covers the results of the {cost} allocation methods for the cooperative game.

\subsection{Instance results}\label{sec:instsRes}
Different versions of the algorithm have been tested varying the size of the B\&B tree where valid inequalities can be added. We consider (i) a pure \textit{branch-and-price} where cut separation is not performed at all, (ii) a partial \textit{branch-and-cut-and-price} where we only allow valid inequalities to be added in the root node, and (iii) a pure \textit{branch-and-cut-and-price} where cuts can be added in all the explored nodes. The RMP model solved is comprised by equations (\ref{GSPP:obj}),(\ref{GSPP:conv}),(\ref{GSPP:pkt}), the linear relaxation of (\ref{GSPP:bin}) and valid inequalities (\ref{GSPP:VIbeta}) that are added by separation.
The algorithm includes a running time-limit and, if it is reached and a gap between the lower and upper bounds still exists, the GSPP formulation problem is solved with all the generated columns in the B\&B tree. This helps tightening the upper bound but requires the integer problem to be solvable in reasonable time. The running time for solving the GSPP is set to 10\% of the algorithm running time. 
Two algorithm time limits of 5 minutes and 3 hours have been tested with an additional (if required) 30 seconds and 18 minutes respectively for solving the GSPP. The model has been entirely written in \textit{Julia} language \citep{bezanson2017a}, modelled using \textit{JuMP} \citep{DunningHuchetteLubin2017} and using \textit{CPLEX v. 12.9} as the solver, allowing 4 threads. It has been tested in an 2.20 GHz Intel Xeon Processor 2650v4 using 4 cores with 32 GB of memory per core.
The MIP formulation from \cite{venturini2017a} has been run in the same machine and solved with the same solver for a fair comparison.
The results are summarized in Tables \ref{tab:ResultsAlg5min}, \ref{tab:ResultsAlg3h}, \ref{tab:ResultsAlg5minHard} and \ref{tab:ResultsAlg3hHard}, that contain the performance comparison on the benchmark instances from \cite{venturini2017a} and the generated set of harder instances with both algorithm time limits.
An instance is represented indicating the number of ships $N$, the number of berthing positions per port $B$, the number of ports $P$ and if the time windows $TW$ are tight $T$ or loose $L$. As indicated in \cite{venturini2017a} a loose time window is approximately 3 times longer than a tight one.
In each instance, all ports have the same amount of berthing positions and all the ships follow the same route and have the same speed profiles but both the MIP and GSPP formulations can account for different amount of berthing positions per port, different ship routes and different ship types.
The set $S$ is discretized in 11 speed levels, covering the range 14-19 knots.
{A very low sulphur fuel oil (VLSFO) is used by the ships which is in accordance with the increasing need of ships to reduce their sulphur emissions. Its price ($F_c$) is computed as the average global price during the first quarter of 2021 corresponding to 500 \$/ton \citep{FuelPrice}.}
{Regarding the cost of the different operational aspects at port, the current literature does not provide a consensus on the costs of waiting, handling and delay time. Moreover, this may fluctuate significantly between ports and in many cases they are not made available to the public due to contractual agreements. \cite{meisel2009a} proposes a delay cost of 1000-3000 \$/hour depending on the ship size and a service cost per quay crane hour of 100 \$. They also consider a speeding-up cost to berth at an earlier time of 1000-3000 \$/hour which can resemble the waiting time cost considered in this study. \cite{venturini2017a} set the terminal handling cost weight to 200 \$/hour and charge an additional 300 \$/hour when there is a delay. They set the cost of waiting one hour at anchorage to 200 \$/hour. For the sake of a fair comparison, we use the same costs as \cite{venturini2017a} which correspond to $H_c = 200$, $D_c = 300$ and $I_c = 200$.}
$LB$ denotes the best lower bound found whereas $Z$ indicates the best integer solution (i.e., upper bound). The optimality gap is stated under the column $Gap$ and it is calculated using the optimal solution, or in the case that this is unknown, the best known solution. The computational time in seconds is given under column $T$.

\begin{landscape} 
\begin{table}[]
\centering
\caption{Computational results on instances from \cite{venturini2017a} with a total time limit of 5 minutes and 30 seconds. The MIP formulation is compared to the variants of the presented branch-and-cut-and-price method.
"-" means that no integer solution has been found within the time limit.
"*" means the time limit has been reached. 
The best running time is highlighted in bold for instances solved to optimality and the best optimality gap for the rest of instances.
}
\label{tab:ResultsAlg5min}
\scalebox{0.75}{
\begin{tabular}{crrrr|rrrr|rrrr|rrrr}
\multicolumn{1}{c|}{\textbf{Instance}} & \multicolumn{4}{c|}{\textbf{MIP formulation}} & \multicolumn{4}{c|}{\textbf{Branch   \& Price}} & \multicolumn{4}{c|}{\textbf{Branch   \& Cut (root node) \& Price}} & \multicolumn{4}{c}{\textbf{Branch   \& Cut \& Price}} \\ \hline
\multicolumn{1}{c|}{\textbf{N-B-P-TW}} & \multicolumn{1}{c}{\textbf{LB}} & \multicolumn{1}{c}{\textbf{Z}} & \multicolumn{1}{c}{\textbf{Gap (\%)}} & \multicolumn{1}{c|}{\textbf{T (s)}} & \multicolumn{1}{c}{\textbf{LB}} & \multicolumn{1}{c}{\textbf{Z}} & \multicolumn{1}{c}{\textbf{Gap (\%)}} & \multicolumn{1}{c|}{\textbf{T (s)}} & \multicolumn{1}{c}{\textbf{LB}} & \multicolumn{1}{c}{\textbf{Z}} & \multicolumn{1}{c}{\textbf{Gap (\%)}} & \multicolumn{1}{c|}{\textbf{T (s)}} & \multicolumn{1}{c}{\textbf{LB}} & \multicolumn{1}{c}{\textbf{Z}} & \multicolumn{1}{c}{\textbf{Gap (\%)}} & \multicolumn{1}{c}{\textbf{T (s)}} \\ \hline
\multicolumn{1}{c|}{4-3-3-L} & 296,600 & 296,600 & 0.00 & \textbf{0.1} & 296,600 & 296,600 & 0.00 & 0.5 & 296,600 & 296,600 & 0.00 & 0.2 & 296,600 & 296,600 & 0.00 & 0.2 \\
\multicolumn{1}{c|}{5-3-3-L} & 394,300 & 394,300 & 0.00 & \textbf{0.4} & 394,300 & 394,300 & 0.00 & 3.2 & 394,300 & 394,300 & 0.00 & 7.2 & 394,300 & 394,300 & 0.00 & 7.2 \\
\multicolumn{1}{c|}{6-3-3-L} & 421,720 & 421,720 & 0.00 & 2.3 & 421,720 & 421,720 & 0.00 & 0.8 & 421,720 & 421,720 & 0.00 & \textbf{0.5} & 421,720 & 421,720 & 0.00 & \textbf{0.5} \\
\multicolumn{1}{c|}{6-3-4-L} & 647,480 & 647,480 & 0.00 & 89.2 & 647,480 & 647,480 & 0.00 & \textbf{39.0} & 647,480 & 647,480 & 0.00 & 111.7 & 647,480 & 647,480 & 0.00 & 103.7 \\
\multicolumn{1}{c|}{10-4-4-L} & 1,014,437 & 1,060,900 & 3.80 & * & 1,053,030 & 1,054,700 & 0.14 & * & 1,053,092 & 1,055,300 & 0.13 & * & 1,053,295 & 1,055,000 & \textbf{0.11} & * \\
\multicolumn{1}{c|}{10-4-3-L} & 689,858 & 700,000 & 1.18 & * & 698,100 & 698,100 & 0.00 & 193.6 & 698,100 & 698,100 & 0.00 & \textbf{72.7} & 698,100 & 698,100 & 0.00 & 115.6 \\
\multicolumn{1}{c|}{4-4-4-L} & 405,120 & 405,120 & 0.00 & \textbf{0.3} & 405,120 & 405,120 & 0.00 & 0.6 & 405,120 & 405,120 & 0.00 & 0.6 & 405,120 & 405,120 & 0.00 & 0.6 \\
\multicolumn{1}{c|}{5-4-4-L} & 500,600 & 500,600 & 0.00 & \textbf{0.4} & 500,600 & 500,600 & 0.00 & 0.9 & 500,600 & 500,600 & 0.00 & 0.8 & 500,600 & 500,600 & 0.00 & 0.8 \\
\multicolumn{1}{c|}{6-4-4-L} & 599,980 & 599,980 & 0.00 & \textbf{1.2} & 599,980 & 599,980 & 0.00 & 7.6 & 599,980 & 599,980 & 0.00 & 5.1 & 599,980 & 599,980 & 0.00 & 5.3 \\
\multicolumn{1}{c|}{12-5-3-L} & 811,139 & 840,640 & 2.32 & * & 830,440 & 830,440 & 0.00 & 136.0 & 830,440 & 830,440 & 0.00 & \textbf{112.1} & 830,440 & 830,440 & 0.00 & 120.4 \\
\multicolumn{1}{c|}{10-6-3-L} & 680,600 & 680,600 & 0.00 & 219.4 & 680,600 & 680,600 & 0.00 & 9.7 & 680,600 & 680,600 & 0.00 & \textbf{5.1} & 680,600 & 680,600 & 0.00 & 5.8 \\
\multicolumn{1}{c|}{11-6-3-L} & 740,430 & 749,620 & 0.78 & * & 746,220 & 746,220 & 0.00 & 22.5 & 746,220 & 746,220 & 0.00 & \textbf{12.7} & 746,220 & 746,220 & 0.00 & 14.3 \\
\multicolumn{1}{c|}{12-6-3-L} & 805,930 & 810,740 & 0.48 & * & 809,840 & 809,840 & 0.00 & 112.6 & 809,840 & 809,840 & 0.00 & 79.0 & 809,840 & 809,840 & 0.00 & \textbf{72.5} \\
\multicolumn{1}{c|}{10-5-4-L} & 1,006,635 & 1,031,100 & 2.11 & * & 1,027,592 & 1,028,320 & 0.07 & * & 1,028,194 & 1,028,320 & \textbf{0.01} & * & 1,027,233 & 1,028,320 & 0.11 & * \\
\multicolumn{1}{c|}{15-10-3-L} & 1,006,000 & 1,006,200 & 0.02 & * & 1,006,200 & 1,006,200 & 0.00 & 46.4 & 1,006,200 & 1,006,200 & 0.00 & 27.8 & 1,006,200 & 1,006,200 & 0.00 & \textbf{25.0} \\
\multicolumn{1}{c|}{15-12-3-L} & 1,001,200 & 1,002,800 & 0.16 & * & 1,002,800 & 1,002,800 & 0.00 & \textbf{1.5} & 1,002,800 & 1,002,800 & 0.00 & \textbf{1.5} & 1,002,800 & 1,002,800 & 0.00 & \textbf{1.5} \\
\multicolumn{1}{c|}{15-10-4-L} & 1,459,400 & 1,459,600 & 0.01 & * & 1,459,600 & 1,459,600 & 0.00 & 23.8 & 1,459,600 & 1,459,600 & 0.00 & 40.2 & 1,459,600 & 1,459,600 & 0.00 & \textbf{13.5} \\
\multicolumn{1}{c|}{20-10-3-L} & 1,341,640 & - & 0.23 & * & 1,344,446 & 1,344,800 & 0.03 & * & 1,344,450 & 1,344,800 & 0.03 & * & 1,344,467 & 1,344,800 & \textbf{0.02} & * \\
\multicolumn{1}{c|}{20-12-3-L} & 1,331,640 & 1,343,000 & 0.36 & * & 1,336,400 & 1,336,400 & 0.00 & \textbf{2.1} & 1,336,400 & 1,336,400 & 0.00 & 2.2 & 1,336,400 & 1,336,400 & 0.00 & 2.2 \\
\multicolumn{1}{c|}{4-3-3-T} & 318,440 & 318,440 & 0.00 & 0.3 & 318,440 & 318,440 & 0.00 & 0.7 & 318,440 & 318,440 & 0.00 & \textbf{0.2} & 318,440 & 318,440 & 0.00 & \textbf{0.2} \\
\multicolumn{1}{c|}{5-3-3-T} & 405,240 & 405,240 & 0.00 & \textbf{0.6} & 405,240 & 405,240 & 0.00 & 1.5 & 405,240 & 405,240 & 0.00 & 1.3 & 405,240 & 405,240 & 0.00 & 1.1 \\
\multicolumn{1}{c|}{6-3-3-T} & 510,920 & 510,920 & 0.00 & 4.0 & 510,920 & 510,920 & 0.00 & 1.9 & 510,920 & 510,920 & 0.00 & 0.9 & 510,920 & 510,920 & 0.00 & \textbf{0.8} \\
\multicolumn{1}{c|}{6-3-4-T} & 993,460 & 993,460 & 0.00 & 3.5 & 993,460 & 993,460 & 0.00 & \textbf{1.2} & 993,460 & 993,460 & 0.00 & 1.3 & 993,460 & 993,460 & 0.00 & \textbf{1.2} \\
\multicolumn{1}{c|}{10-4-4-T} & 1,574,771 & 1,676,990 & 5.17 & * & 1,660,640 & 1,660,640 & 0.00 & 101.3 & 1,660,640 & 1,660,640 & 0.00 & 63.9 & 1,660,640 & 1,660,640 & 0.00 & \textbf{61.0} \\
\multicolumn{1}{c|}{10-4-3-T} & 973,445 & 1,023,890 & 4.77 & * & 1,022,200 & 1,022,200 & 0.00 & 12.3 & 1,022,200 & 1,022,200 & 0.00 & \textbf{6.6} & 1,022,200 & 1,022,200 & 0.00 & 7.9 \\
\multicolumn{1}{c|}{4-4-4-T} & 442,600 & 442,600 & 0.00 & 0.9 & 442,600 & 442,600 & 0.00 & 1.1 & 442,600 & 442,600 & 0.00 & 0.7 & 442,600 & 442,600 & 0.00 & \textbf{0.6} \\
\multicolumn{1}{c|}{5-4-4-T} & 576,010 & 576,010 & 0.00 & \textbf{4.1} & 576,010 & 576,010 & 0.00 & 10.3 & 576,010 & 576,010 & 0.00 & 6.0 & 576,010 & 576,010 & 0.00 & 6.6 \\
\multicolumn{1}{c|}{6-4-4-T} & 653,560 & 653,560 & 0.00 & 11.5 & 653,560 & 653,560 & 0.00 & 23.5 & 653,560 & 653,560 & 0.00 & \textbf{8.8} & 653,560 & 653,560 & 0.00 & 10.4 \\
\multicolumn{1}{c|}{12-5-3-T} & 811,240 & 835,740 & 2.31 & * & 830,440 & 830,440 & 0.00 & 128.0 & 830,440 & 830,440 & 0.00 & \textbf{68.3} & 830,440 & 830,440 & 0.00 & 96.1 \\
\multicolumn{1}{c|}{12-6-3-T} & 805,180 & 823,240 & 1.67 & * & 818,840 & 818,840 & 0.00 & 173.6 & 818,840 & 818,840 & 0.00 & 163.7 & 818,840 & 818,840 & 0.00 & \textbf{158.6} \\
\multicolumn{1}{c|}{10-5-4-T} & 1,117,723 & 1,147,530 & 2.31 & * & 1,144,160 & 1,144,160 & 0.00 & 141.4 & 1,144,160 & 1,144,160 & 0.00 & \textbf{68.6} & 1,144,160 & 1,144,160 & 0.00 & 72.6 \\
\multicolumn{1}{c|}{15-10-4-T} & 1,575,640 & 1,605,460 & 1.34 & * & 1,597,100 & 1,597,100 & 0.00 & \textbf{9.3} & 1,597,100 & 1,597,100 & 0.00 & 11.5 & 1,597,100 & 1,597,100 & 0.00 & 13.7 \\
\multicolumn{1}{c|}{20-10-3-T} & 1,551,597 & - & 4.78 & * & 1,629,000 & 1,629,500 & \textbf{0.03} & * & 1,629,000 & 1,629,500 & \textbf{0.03} & * & 1,629,000 & 1,629,500 & \textbf{0.03} & * \\
\multicolumn{1}{c|}{20-12-3-T} & 1,541,949 & 1,628,900 & 4.02 & * & 1,606,500 & 1,606,500 & 0.00 & 46.3 & 1,606,500 & 1,606,500 & 0.00 & 45.4 & 1,606,500 & 1,606,500 & 0.00 & \textbf{42.5} \\ \hline
\multicolumn{1}{c|}{\textbf{Average}} & \textbf{} & \textbf{} & \textbf{1.113} & \textbf{} & \textbf{} & \textbf{} & \textbf{0.0079} & \textbf{} & \textbf{} & \textbf{} & \textbf{0.0060} & \textbf{} & \textbf{} & \textbf{} & \textbf{0.0081} & \textbf{} \\
\multicolumn{2}{c}{\textbf{Optimal   solutions}} & \textbf{} & \textbf{15/34} & \textbf{} & \textbf{} & \textbf{} & \textbf{30/34} & \textbf{} & \textbf{} & \textbf{} & \textbf{30/34} & \textbf{} & \textbf{} & \textbf{} & \textbf{30/34} & \textbf{}
\end{tabular}
}
\end{table}
\end{landscape}

\begin{landscape}
\begin{table}[]
\centering
\caption{Computational results on instances from \cite{venturini2017a} with a total time limit of 3 hours and 18 minutes. The MIP formulation is compared to the variants of the presented branch-and-cut-and-price method.
"*" means the time limit has been reached. 
The best running time is highlighted in bold for instances solved to optimality and the best optimality gap for the rest of instances.
}
\label{tab:ResultsAlg3h}
\scalebox{0.74}{
\begin{tabular}{crrrr|rrrr|rrrr|rrrr}
\multicolumn{1}{c|}{\textbf{Instance}} & \multicolumn{4}{c|}{\textbf{MIP formulation}} & \multicolumn{4}{c|}{\textbf{Branch   \& Price}} & \multicolumn{4}{c|}{\textbf{Branch   \& Cut (root node) \& Price}} & \multicolumn{4}{c}{\textbf{Branch   \& Cut \& Price}} \\ \hline
\multicolumn{1}{c|}{\textbf{N-B-P-TW}} & \multicolumn{1}{c}{\textbf{LB}} & \multicolumn{1}{c}{\textbf{Z}} & \multicolumn{1}{c}{\textbf{Gap (\%)}} & \multicolumn{1}{c|}{\textbf{T (s)}} & \multicolumn{1}{c}{\textbf{LB}} & \multicolumn{1}{c}{\textbf{Z}} & \multicolumn{1}{c}{\textbf{Gap (\%)}} & \multicolumn{1}{c|}{\textbf{T (s)}} & \multicolumn{1}{c}{\textbf{LB}} & \multicolumn{1}{c}{\textbf{Z}} & \multicolumn{1}{c}{\textbf{Gap (\%)}} & \multicolumn{1}{c|}{\textbf{T (s)}} & \multicolumn{1}{c}{\textbf{LB}} & \multicolumn{1}{c}{\textbf{Z}} & \multicolumn{1}{c}{\textbf{Gap (\%)}} & \multicolumn{1}{c}{\textbf{T (s)}} \\ \hline
\multicolumn{1}{c|}{4-3-3-L} & 296,600 & 296,600 & 0.00 & \textbf{0.1} & 296,600 & 296,600 & 0.00 & 0.5 & 296,600 & 296,600 & 0.00 & 0.2 & 296,600 & 296,600 & 0.00 & 0.2 \\
\multicolumn{1}{c|}{5-3-3-L} & 394,300 & 394,300 & 0.00 & \textbf{0.4} & 394,300 & 394,300 & 0.00 & 3.2 & 394,300 & 394,300 & 0.00 & 7.2 & 394,300 & 394,300 & 0.00 & 7.2 \\
\multicolumn{1}{c|}{6-3-3-L} & 421,679 & 421,720 & 0.01 & 2.3 & 421,720 & 421,720 & 0.00 & 0.8 & 421,720 & 421,720 & 0.00 & \textbf{0.5} & 421,720 & 421,720 & 0.00 & \textbf{0.5} \\
\multicolumn{1}{c|}{6-3-4-L} & 647,423 & 647,480 & 0.01 & 89.2 & 647,480 & 647,480 & 0.00 & \textbf{39.0} & 647,480 & 647,480 & 0.00 & 111.7 & 647,480 & 647,480 & 0.00 & 103.7 \\
\multicolumn{1}{c|}{10-4-4-L} & 1,020,581 & 1,055,800 & 3.22 & * & 1,054,500 & 1,054,500 & 0.00 & \textbf{5563.8} & 1,054,500 & 1,054,500 & 0.00 & 6068.2 & 1,054,300 & 1,054,500 & \textbf{0.02} & * \\
\multicolumn{1}{c|}{10-4-3-L} & 694,451 & 699,000 & 0.52 & * & 698,100 & 698,100 & 0.00 & 193.6 & 698,100 & 698,100 & 0.00 & \textbf{72.7} & 698,100 & 698,100 & 0.00 & 115.6 \\
\multicolumn{1}{c|}{4-4-4-L} & 405,120 & 405,120 & 0.00 & \textbf{0.3} & 405,120 & 405,120 & 0.00 & 0.6 & 405,120 & 405,120 & 0.00 & 0.6 & 405,120 & 405,120 & 0.00 & 0.6 \\
\multicolumn{1}{c|}{5-4-4-L} & 500,600 & 500,600 & 0.00 & \textbf{0.4} & 500,600 & 500,600 & 0.00 & 0.9 & 500,600 & 500,600 & 0.00 & 0.8 & 500,600 & 500,600 & 0.00 & 0.8 \\
\multicolumn{1}{c|}{6-4-4-L} & 599,980 & 599,980 & 0.00 & \textbf{1.2} & 599,980 & 599,980 & 0.00 & 7.6 & 599,980 & 599,980 & 0.00 & 5.1 & 599,980 & 599,980 & 0.00 & 5.3 \\
\multicolumn{1}{c|}{12-5-3-L} & 813,713 & 834,740 & 2.01 & * & 830,440 & 830,440 & 0.00 & 136.0 & 830,440 & 830,440 & 0.00 & \textbf{112.1} & 830,440 & 830,440 & 0.00 & 120.4 \\
\multicolumn{1}{c|}{10-6-3-L} & 680,600 & 680,600 & 0.00 & 219.4 & 680,600 & 680,600 & 0.00 & 9.7 & 680,600 & 680,600 & 0.00 & \textbf{5.1} & 680,600 & 680,600 & 0.00 & 5.8 \\
\multicolumn{1}{c|}{11-6-3-L} & 746,220 & 746,220 & 0.00 & 8705.8 & 746,220 & 746,220 & 0.00 & 22.5 & 746,220 & 746,220 & 0.00 & \textbf{12.7} & 746,220 & 746,220 & 0.00 & 14.3 \\
\multicolumn{1}{c|}{12-6-3-L} & 809,840 & 809,840 & 0.00 & 5032.0 & 809,840 & 809,840 & 0.00 & 112.6 & 809,840 & 809,840 & 0.00 & 79.0 & 809,840 & 809,840 & 0.00 & \textbf{72.5} \\
\multicolumn{1}{c|}{10-5-4-L} & 1,013,114 & 1,029,300 & 1.48 & * & 1,028,320 & 1,028,320 & 0.00 & 589.0 & 1,028,320 & 1,028,320 & 0.00 & \textbf{366.1} & 1,028,320 & 1,028,320 & 0.00 & 1135.1 \\
\multicolumn{1}{c|}{15-10-3-L} & 1,006,200 & 1,006,200 & 0.00 & 3259.5 & 1,006,200 & 1,006,200 & 0.00 & 46.4 & 1,006,200 & 1,006,200 & 0.00 & 27.8 & 1,006,200 & 1,006,200 & 0.00 & \textbf{25.0} \\
\multicolumn{1}{c|}{15-12-3-L} & 1,002,240 & 1,002,800 & 0.06 & * & 1,002,800 & 1,002,800 & 0.00 & \textbf{1.5} & 1,002,800 & 1,002,800 & 0.00 & \textbf{1.5} & 1,002,800 & 1,002,800 & 0.00 & \textbf{1.5} \\
\multicolumn{1}{c|}{15-10-4-L} & 1,459,600 & 1,459,600 & 0.00 & 1703.1 & 1,459,600 & 1,459,600 & 0.00 & 23.8 & 1,459,600 & 1,459,600 & 0.00 & 40.2 & 1,459,600 & 1,459,600 & 0.00 & \textbf{13.5} \\
\multicolumn{1}{c|}{20-10-3-L} & 1,341,640 & 1,346,000 & 0.23 & * & 1,344,520 & 1,344,800 & 0.02 & * & 1,344,525 & 1,344,800 & 0.02 & * & 1,344,600 & 1,344,800 & \textbf{0.01} & * \\
\multicolumn{1}{c|}{20-12-3-L} & 1,331,680 & 1,337,400 & 0.35 & * & 1,336,400 & 1,336,400 & 0.00 & \textbf{2.1} & 1,336,400 & 1,336,400 & 0.00 & 2.2 & 1,336,400 & 1,336,400 & 0.00 & 2.2 \\
\multicolumn{1}{c|}{4-3-3-T} & 318,440 & 318,440 & 0.00 & 0.3 & 318,440 & 318,440 & 0.00 & 0.7 & 318,440 & 318,440 & 0.00 & \textbf{0.2} & 318,440 & 318,440 & 0.00 & \textbf{0.2} \\
\multicolumn{1}{c|}{5-3-3-T} & 405,240 & 405,240 & 0.00 & \textbf{0.6} & 405,240 & 405,240 & 0.00 & 1.5 & 405,240 & 405,240 & 0.00 & 1.3 & 405,240 & 405,240 & 0.00 & 1.1 \\
\multicolumn{1}{c|}{6-3-3-T} & 510,920 & 510,920 & 0.00 & 4.0 & 510,920 & 510,920 & 0.00 & 1.9 & 510,920 & 510,920 & 0.00 & 0.9 & 510,920 & 510,920 & 0.00 & \textbf{0.8} \\
\multicolumn{1}{c|}{6-3-4-T} & 993,460 & 993,460 & 0.00 & 3.5 & 993,460 & 993,460 & 0.00 & \textbf{1.2} & 993,460 & 993,460 & 0.00 & 1.3 & 993,460 & 993,460 & 0.00 & \textbf{1.2} \\
\multicolumn{1}{c|}{10-4-4-T} & 1,660,640 & 1,660,640 & 0.00 & 1660.0 & 1,660,640 & 1,660,640 & 0.00 & 101.3 & 1,660,640 & 1,660,640 & 0.00 & 63.9 & 1,660,640 & 1,660,640 & 0.00 & \textbf{61.0} \\
\multicolumn{1}{c|}{10-4-3-T} & 1,022,200 & 1,022,200 & 0.00 & 562.5 & 1,022,200 & 1,022,200 & 0.00 & 12.3 & 1,022,200 & 1,022,200 & 0.00 & \textbf{6.6} & 1,022,200 & 1,022,200 & 0.00 & 7.9 \\
\multicolumn{1}{c|}{4-4-4-T} & 442,600 & 442,600 & 0.00 & 0.9 & 442,600 & 442,600 & 0.00 & 1.1 & 442,600 & 442,600 & 0.00 & 0.7 & 442,600 & 442,600 & 0.00 & \textbf{0.6} \\
\multicolumn{1}{c|}{5-4-4-T} & 576,010 & 576,010 & 0.00 & \textbf{4.1} & 576,010 & 576,010 & 0.00 & 10.3 & 576,010 & 576,010 & 0.00 & 6.0 & 576,010 & 576,010 & 0.00 & 6.6 \\
\multicolumn{1}{c|}{6-4-4-T} & 653,560 & 653,560 & 0.00 & 11.5 & 653,560 & 653,560 & 0.00 & 23.5 & 653,560 & 653,560 & 0.00 & \textbf{8.8} & 653,560 & 653,560 & 0.00 & 10.4 \\
\multicolumn{1}{c|}{12-5-3-T} & 817,533 & 830,440 & 1.55 & * & 830,440 & 830,440 & 0.00 & 128.0 & 830,440 & 830,440 & 0.00 & \textbf{68.3} & 830,440 & 830,440 & 0.00 & 96.1 \\
\multicolumn{1}{c|}{12-6-3-T} & 810,476 & 821,540 & 1.02 & * & 818,840 & 818,840 & 0.00 & 173.6 & 818,840 & 818,840 & 0.00 & 163.7 & 818,840 & 818,840 & 0.00 & \textbf{158.6} \\
\multicolumn{1}{c|}{10-5-4-T} & 1,144,160 & 1,144,160 & 0.00 & 3649.8 & 1,144,160 & 1,144,160 & 0.00 & 141.4 & 1,144,160 & 1,144,160 & 0.00 & \textbf{68.6} & 1,144,160 & 1,144,160 & 0.00 & 72.6 \\
\multicolumn{1}{c|}{15-10-4-T} & 1,584,071 & 1,597,620 & 0.82 & * & 1,597,100 & 1,597,100 & 0.00 & \textbf{9.3} & 1,597,100 & 1,597,100 & 0.00 & 11.5 & 1,597,100 & 1,597,100 & 0.00 & 13.7 \\
\multicolumn{1}{c|}{20-10-3-T} & 1,552,283 & 1,634,900 & 4.74 & * & 1,629,380 & 1,629,500 & \textbf{0.01} & * & 1,629,300 & 1,629,500 & 0.01 & * & 1,629,380 & 1,629,500 & \textbf{0.01} & * \\
\multicolumn{1}{c|}{20-12-3-T} & 1,545,006 & 1,609,100 & 3.83 & * & 1,606,500 & 1,606,500 & 0.00 & 46.3 & 1,606,500 & 1,606,500 & 0.00 & 45.4 & 1,606,500 & 1,606,500 & 0.00 & \textbf{42.5} \\ \hline
\multicolumn{1}{c|}{\textbf{Average}} & \textbf{} & \textbf{} & \textbf{0.584} & \textbf{} & \textbf{} & \textbf{} & \textbf{0.0008} & \textbf{} & \textbf{} & \textbf{} & \textbf{0.0010} & \textbf{} & \textbf{} & \textbf{} & \textbf{0.0012} & \textbf{} \\
\multicolumn{2}{c}{\textbf{Optimal   solutions}} & \textbf{} & \textbf{24/34} & \textbf{} & \textbf{} & \textbf{} & \textbf{32/34} & \textbf{} & \textbf{} & \textbf{} & \textbf{32/34} & \textbf{} & \textbf{} & \textbf{} & \textbf{31/34} & \textbf{}
\end{tabular}
}
\end{table}
\end{landscape}

\begin{landscape}
\begin{table}[]
\centering
\caption{Computational results on the set of harder instances with a total time limit of 5 minutes and 30 seconds. The MIP formulation is compared to the variants of the presented branch-and-cut-and-price method.
"-" means that no integer solution has been found within the time limit. 
"*" means the time limit has been reached. 
The best running time is highlighted in bold for instances solved to optimality and the best optimality gap for the rest of instances.
}
\label{tab:ResultsAlg5minHard}
\scalebox{0.75}{
\begin{tabular}{crrrr|rrrr|rrrr|rrrr}
\multicolumn{1}{c|}{\textbf{Instance}} & \multicolumn{4}{c|}{\textbf{MIP formulation}} & \multicolumn{4}{c|}{\textbf{Branch   \& Price}} & \multicolumn{4}{c|}{\textbf{Branch   \& Cut (root node) \& Price}} & \multicolumn{4}{c}{\textbf{Branch   \& Cut \& Price}} \\ \hline
\multicolumn{1}{c|}{\textbf{N-B-P-TW}} & \multicolumn{1}{c}{\textbf{LB}} & \multicolumn{1}{c}{\textbf{Z}} & \multicolumn{1}{c}{\textbf{Gap (\%)}} & \multicolumn{1}{c|}{\textbf{T (s)}} & \multicolumn{1}{c}{\textbf{LB}} & \multicolumn{1}{c}{\textbf{Z}} & \multicolumn{1}{c}{\textbf{Gap (\%)}} & \multicolumn{1}{c|}{\textbf{T (s)}} & \multicolumn{1}{c}{\textbf{LB}} & \multicolumn{1}{c}{\textbf{Z}} & \multicolumn{1}{c}{\textbf{Gap (\%)}} & \multicolumn{1}{c|}{\textbf{T (s)}} & \multicolumn{1}{c}{\textbf{LB}} & \multicolumn{1}{c}{\textbf{Z}} & \multicolumn{1}{c}{\textbf{Gap (\%)}} & \multicolumn{1}{c}{\textbf{T (s)}} \\ \hline
\multicolumn{1}{c|}{25-12-3-L} & 1,668,080 & - & 0.56 & * & 1,677,200 & 1,677,400 & \textbf{0.01} & * & 1,677,200 & 1,677,400 & \textbf{0.01} & * & 1,677,200 & 1,677,400 & \textbf{0.01} & * \\
\multicolumn{1}{c|}{25-12-3-T} & 1,923,560 & - & 6.03 & * & 2,046,930 & 2,047,100 & 0.00 & * & 2,046,933 & 2,047,200 & \textbf{0.00} & * & 2,046,933 & 2,047,200 & \textbf{0.00} & * \\
\multicolumn{1}{c|}{12-5-4-L} & 1,201,546 & 1,253,160 & 3.50 & * & 1,239,277 & 1,245,360 & 0.47 & * & 1,240,826 & 1,247,760 & 0.35 & * & 1,240,862 & 1,247,060 & \textbf{0.35} & * \\
\multicolumn{1}{c|}{12-5-4-T} & 1,324,428 & - & 5.69 & * & 1,398,848 & 1,410,270 & 0.39 & * & 1,399,631 & 1,408,070 & \textbf{0.33} & * & 1,398,580 & 1,408,220 & 0.41 & * \\
\multicolumn{1}{c|}{30-12-3-L} & 1,997,400 & - & 0.95 & * & 2,016,178 & 2,016,600 & \textbf{0.02} & * & 2,016,178 & 2,016,600 & \textbf{0.02} & * & 2,016,178 & 2,016,600 & \textbf{0.02} & * \\
\multicolumn{1}{c|}{30-12-3-T} & 2,304,432 & - & 7.55 & * & 2,491,406 & 2,496,900 & \textbf{0.04} & * & 2,491,406 & 2,495,300 & \textbf{0.04} & * & 2,491,406 & 2,495,600 & \textbf{0.04} & * \\
\multicolumn{1}{c|}{20-12-4-L} & 1,934,640 & - & 0.47 & * & 1,943,486 & 1,943,800 & 0.02 & * & 1,943,500 & 1,943,800 & 0.02 & * & 1,943,800 & 1,943,800 & 0.00 & \textbf{243.7} \\
\multicolumn{1}{c|}{20-12-4-T} & 3,050,995 & - & 2.00 & * & 3,113,170 & 3,113,170 & 0.00 & 42.5 & 3,113,170 & 3,113,170 & 0.00 & 15.8 & 3,113,170 & 3,113,170 & 0.00 & \textbf{14.9} \\
\multicolumn{1}{c|}{15-8-4-L} & 1,471,903 & 1,508,500 & 1.68 & * & 1,495,225 & 1,497,100 & 0.12 & * & 1,496,131 & 1,497,300 & \textbf{0.06} & \textbf{*} & 1,496,118 & 1,497,100 & 0.06 & * \\
\multicolumn{1}{c|}{15-8-4-T} & 1,599,496 & - & 3.41 & * & 1,654,848 & 1,656,260 & 0.07 & * & 1,655,376 & 1,656,260 & 0.04 & * & 1,655,416 & 1,656,260 & \textbf{0.04} & * \\
\multicolumn{1}{c|}{25-12-4-L} & 2,419,800 & - & 0.85 & * & 2,439,377 & 2,440,700 & 0.05 & * & 2,439,531 & 2,440,600 & \textbf{0.04} & * & 2,439,380 & 2,441,500 & 0.05 & * \\
\multicolumn{1}{c|}{25-12-4-T} & 3,560,100 & - & 3.97 & * & 3,706,995 & 3,707,390 & 0.01 & * & 3,707,182 & 3,707,390 & 0.01 & * & 3,707,390 & 3,707,390 & \textbf{0.00} & * \\
\multicolumn{1}{c|}{30-15-4-L} & 2,905,000 & - & 0.47 & * & 2,918,400 & 2,918,800 & 0.01 & * & 2,918,420 & 2,918,800 & 0.01 & * & 2,918,436 & 2,918,800 & \textbf{0.01} & * \\
\multicolumn{1}{c|}{30-15-4-T} & 3,109,816 & - & 5.14 & * & 3,274,118 & 3,278,880 & 0.13 & * & 3,274,390 & 3,283,550 & 0.12 & * & 3,274,441 & 3,280,400 & \textbf{0.12} & * \\
\multicolumn{1}{c|}{40-15-3-L} & 2,658,400 & - & 1.15 & * & 2,689,060 & 2,689,200 & \textbf{0.01} & * & 2,689,060 & 2,689,200 & \textbf{0.01} & * & 2,689,060 & 2,689,200 & \textbf{0.01} & * \\
\multicolumn{1}{c|}{40-15-3-T} & 3,067,200 & - & 7.89 & * & 3,329,567 & 3,330,900 & \textbf{0.02} & * & 3,329,567 & 3,330,900 & \textbf{0.02} & * & 3,329,567 & 3,331,300 & \textbf{0.02} & * \\ \hline
\multicolumn{1}{c|}{\textbf{Average}} & \textbf{} & \textbf{} & \textbf{3.206} & \textbf{} & \textbf{} & \textbf{} & \textbf{0.0849} & \textbf{} & \textbf{} & \textbf{} & \textbf{0.0665} & \textbf{} & \textbf{} & \textbf{} & \textbf{0.0698} & \textbf{} \\
\multicolumn{2}{c}{\textbf{Optimal   solutions}} & \textbf{} & \textbf{0/16} & \textbf{} & \textbf{} & \textbf{} & \textbf{1/16} & \textbf{} & \textbf{} & \textbf{} & \textbf{1/16} & \textbf{} & \textbf{} & \textbf{} & \textbf{3/16} & \textbf{}
\end{tabular}
}
\end{table}
\end{landscape}

\begin{landscape}
\begin{table}[]
\centering
\caption{Computational results on the set of harder instances with a total time limit of 3 hours and 18 minutes. The MIP formulation is compared to the variants of the presented branch-and-cut-and-price method.
"-" means that no integer solution has been found within the time limit. 
"*" means the time limit has been reached. 
The best running time is highlighted in bold for instances solved to optimality and the best optimality gap for the rest of instances.
}
\label{tab:ResultsAlg3hHard}
\scalebox{0.74}{
\begin{tabular}{crrrr|rrrr|rrrr|rrrr}
\multicolumn{1}{c|}{\textbf{Instance}} & \multicolumn{4}{c|}{\textbf{MIP formulation}} & \multicolumn{4}{c|}{\textbf{Branch   \& Price}} & \multicolumn{4}{c|}{\textbf{Branch   \& Cut (root node) \& Price}} & \multicolumn{4}{c}{\textbf{Branch   \& Cut \& Price}} \\ \hline
\multicolumn{1}{c|}{\textbf{N-B-P-TW}} & \multicolumn{1}{c}{\textbf{LB}} & \multicolumn{1}{c}{\textbf{Z}} & \multicolumn{1}{c}{\textbf{Gap (\%)}} & \multicolumn{1}{c|}{\textbf{T (s)}} & \multicolumn{1}{c}{\textbf{LB}} & \multicolumn{1}{c}{\textbf{Z}} & \multicolumn{1}{c}{\textbf{Gap (\%)}} & \multicolumn{1}{c|}{\textbf{T (s)}} & \multicolumn{1}{c}{\textbf{LB}} & \multicolumn{1}{c}{\textbf{Z}} & \multicolumn{1}{c}{\textbf{Gap (\%)}} & \multicolumn{1}{c|}{\textbf{T (s)}} & \multicolumn{1}{c}{\textbf{LB}} & \multicolumn{1}{c}{\textbf{Z}} & \multicolumn{1}{c}{\textbf{Gap (\%)}} & \multicolumn{1}{c}{\textbf{T (s)}} \\ \hline
\multicolumn{1}{c|}{25-12-3-L} & 1,668,080 & 1,679,400 & 0.56 & * & 1,677,200 & 1,677,400 & \textbf{0.01} & * & 1,677,200 & 1,677,400 & \textbf{0.01} & * & 1,677,200 & 1,677,400 & \textbf{0.01} & * \\
\multicolumn{1}{c|}{25-12-3-T} & 1,932,150 & - & 5.61 & * & 2,046,930 & 2,047,100 & 0.00 & * & 2,046,933 & 2,047,100 & \textbf{0.00} & * & 2,046,933 & 2,047,000 & \textbf{0.00} & * \\
\multicolumn{1}{c|}{12-5-4-L} & 1,205,096 & 1,253,160 & 3.22 & * & 1,244,427 & 1,245,660 & 0.06 & * & 1,245,110 & 1,245,160 & \textbf{0.00} & * & 1,244,080 & 1,245,160 & 0.09 & * \\
\multicolumn{1}{c|}{12-5-4-T} & 1,348,176 & 1,407,980 & 4.00 & * & 1,404,280 & 1,404,280 & 0.00 & \textbf{5126.1} & 1,404,280 & 1,404,280 & 0.00 & 9530.6 & 1,403,243 & 1,405,760 & 0.07 & * \\
\multicolumn{1}{c|}{30-12-3-L} & 1,998,483 & - & 0.90 & * & 2,016,178 & 2,016,600 & \textbf{0.02} & * & 2,016,178 & 2,016,600 & \textbf{0.02} & * & 2,016,178 & 2,016,600 & \textbf{0.02} & * \\
\multicolumn{1}{c|}{30-12-3-T} & 2,311,705 & - & 7.25 & * & 2,491,406 & 2,492,500 & \textbf{0.04} & * & 2,491,406 & 2,492,500 & \textbf{0.04} & * & 2,491,406 & 2,492,500 & \textbf{0.04} & * \\
\multicolumn{1}{c|}{20-12-4-L} & 1,934,640 & - & 0.47 & * & 1,943,800 & 1,943,800 & 0.00 & 2544.5 & 1,943,800 & 1,943,800 & 0.00 & 1202.7 & 1,943,800 & 1,943,800 & 0.00 & \textbf{243.7} \\
\multicolumn{1}{c|}{20-12-4-T} & 3,055,040 & - & 1.87 & * & 3,113,170 & 3,113,170 & 0.00 & 42.5 & 3,113,170 & 3,113,170 & 0.00 & 15.8 & 3,113,170 & 3,113,170 & 0.00 & \textbf{14.9} \\
\multicolumn{1}{c|}{15-8-4-L} & 1,475,646 & 1,497,000 & 1.43 & * & 1,497,000 & 1,497,000 & 0.00 & 2870.3 & 1,497,000 & 1,497,000 & 0.00 & \textbf{968.3} & 1,497,000 & 1,497,000 & 0.00 & 1831.1 \\
\multicolumn{1}{c|}{15-8-4-T} & 1,608,977 & 1,671,220 & 2.84 & * & 1,656,040 & 1,656,040 & 0.00 & \textbf{746.7.3} & 1,656,040 & 1,656,040 & 0.00 & 1365.a & 1,656,040 & 1,656,040 & 0.00 & 1534.0 \\
\multicolumn{1}{c|}{25-12-4-L} & 2,420,034 & - & 0.84 & * & 2,439,594 & 2,440,500 & 0.04 & * & 2,439,796 & 2,440,500 & 0.03 & * & 2,439,881 & 2,440,500 & \textbf{0.03} & * \\
\multicolumn{1}{c|}{25-12-4-T} & 3,572,971 & - & 3.63 & * & 3,707,390 & 3,707,390 & 0.00 & 604.3 & 3,707,390 & 3,707,390 & 0.00 & 501.1 & 3,707,390 & 3,707,390 & 0.00 & \textbf{440.8} \\
\multicolumn{1}{c|}{30-15-4-L} & 2,905,000 & - & 0.47 & * & 2,918,400 & 2,918,800 & 0.01 & * & 2,918,425 & 2,918,800 & 0.01 & * & 2,918,441 & 2,918,600 & \textbf{0.01} & * \\
\multicolumn{1}{c|}{30-15-4-T} & 3,123,445 & - & 4.72 & * & 3,274,905 & 3,278,280 & 0.10 & * & 3,275,095 & 3,278,280 & 0.10 & * & 3,275,112 & 3,278,940 & \textbf{0.10} & * \\
\multicolumn{1}{c|}{40-15-3-L} & 2,658,440 & - & 1.14 & * & 2,689,060 & 2,689,200 & \textbf{0.01} & * & 2,689,060 & 2,689,200 & \textbf{0.01} & * & 2,689,060 & 2,689,200 & \textbf{0.01} & * \\
\multicolumn{1}{c|}{40-15-3-T} & 3,071,346 & - & 7.77 & * & 3,329,567 & 3,330,300 & \textbf{0.02} & * & 3,329,567 & 3,330,200 & \textbf{0.02} & * & 3,329,567 & 3,330,100 & \textbf{0.02} & * \\ \hline
\multicolumn{1}{c|}{\textbf{Average}} & \textbf{} & \textbf{} & \textbf{2.919} & \textbf{} & \textbf{} & \textbf{} & \textbf{0.0192} & \textbf{} & \textbf{} & \textbf{} & \textbf{0.0148} & \textbf{} & \textbf{} & \textbf{} & \textbf{0.0243} & \textbf{} \\
\multicolumn{2}{c}{\textbf{Optimal   solutions}} & \textbf{} & \textbf{0/16} & \textbf{} & \textbf{} & \textbf{} & \textbf{6/16} & \textbf{} & \textbf{} & \textbf{} & \textbf{6/16} & \textbf{} & \textbf{} & \textbf{} & \textbf{5/16} & \textbf{}
\end{tabular}
}
\end{table}
\end{landscape}

The results show a better performance of the solution methods based on the GSPP formulation.
All variants of the \textit{branch-and-cut-and-price} method are able to find optimal or near-optimal solutions in less than 6 minutes.
Among the new proposed methods, the one where cutting is performed {at the root node} shows a better performance.
All the solution method variants outperform CPLEX in all instances that require more than
{4 seconds}
to solve and show similar running times for the faster ones. The difference in performance is more notable on the set of harder instances where CPLEX is not able to find a feasible integer solution in {11} out of the 16 instances within 3 hours and 18 minutes and the average optimality gap is above {2.9} \%. Within 5 minutes and 30 seconds, the proposed new methods not only find feasible solutions to all instances but also achieve an optimality gap of {0.03} \%. This gap is further reduced to less than {0.01} \% with a time limit of 3 hours and 18 minutes. The good quality of the solutions in such a short computational time is attractive from an operational point of view where suboptimal solutions usually are not a problem and possible disruptions require rapid computations for new plans.
{Regarding the impact of solving the GSPP at the end, it is higher when the time limits are low. The GSPP is only solved if the method still has not proven optimality and, from those cases, it is found that the GSPP improves the upper bound in 84-100\% of the cases with the 5 minutes and 30 seconds time limit and in 50-57\% of the cases with the 3 hours and 18 minutes, depending on the method. In the vast majority of these cases, an integer solution is not found in the B\&B tree.}

\begin{table}[]
\centering
\caption{Performance of presented methods on a subset of five instances.}
\label{tab:extraInfo}
\resizebox{\textwidth}{!}{%
\begin{tabular}{c|rrrrrrrrrrr}
\textbf{Instance} & \multicolumn{11}{c}{\textbf{Branch \& Price}} \\ \hline
\textbf{N-B-P-TW} & \multicolumn{1}{c}{\textbf{\begin{tabular}[c]{@{}c@{}}Gap\\ (\%)\end{tabular}}} & \multicolumn{1}{c}{\textbf{\begin{tabular}[c]{@{}c@{}}T\\ (s)\end{tabular}}} & \multicolumn{1}{c}{\textbf{\begin{tabular}[c]{@{}c@{}}T RMP\\ (\%)\end{tabular}}} & \multicolumn{1}{c}{\textbf{\begin{tabular}[c]{@{}c@{}}T PP\\ (\%)\end{tabular}}} & \multicolumn{1}{c}{\textbf{\begin{tabular}[c]{@{}c@{}}T Sep\\ (\%)\end{tabular}}} & \multicolumn{1}{c}{\textbf{\begin{tabular}[c]{@{}c@{}}T Branch\\ (\%)\end{tabular}}} & \multicolumn{1}{c}{\textbf{\begin{tabular}[c]{@{}c@{}}T GSPP\\ (\%)\end{tabular}}} & \multicolumn{1}{c}{\textbf{Nodes}} & \multicolumn{1}{c}{\textbf{CG Its}} & \multicolumn{1}{c}{\textbf{Cols}} & \multicolumn{1}{c}{\textbf{Cuts}} \\ \hline
12-6-3-T & 0.00 & 174 & 75.8 & 19.6 & 0.0 & 1.5 & 0.0 & 219 & 3,820 & 16,928 & 0 \\
10-4-4-L & 0.00 & 5,564 & 89.2 & 10.3 & 0.0 & 0.1 & 0.0 & 233 & 14,745 & 91,489 & 0 \\
20-10-3-L & 0.02 & 11,031 & 28.7 & 28.8 & 0.0 & 18.7 & 0.0 & 13,907 & 45,145 & 112,557 & 0 \\
15-8-4-L & 0.00 & 2,870 & 48.3 & 49.8 & 0.0 & 0.6 & 0.0 & 327 & 7,553 & 48,756 & 0 \\
40-15-3-T & 0.02 & 11,097 & 74.3 & 10.5 & 0.0 & 2.6 & 2.1 & 899 & 6,678 & 99,529 & 0 \\ \hline
\textbf{Instance} & \multicolumn{11}{c}{\textbf{Branch \& Cut (root node) \& Price}} \\ \hline
\textbf{N-B-P-TW} & \multicolumn{1}{c}{\textbf{\begin{tabular}[c]{@{}c@{}}Gap\\ (\%)\end{tabular}}} & \multicolumn{1}{c}{\textbf{\begin{tabular}[c]{@{}c@{}}T\\ (s)\end{tabular}}} & \multicolumn{1}{c}{\textbf{\begin{tabular}[c]{@{}c@{}}T RMP\\ (\%)\end{tabular}}} & \multicolumn{1}{c}{\textbf{\begin{tabular}[c]{@{}c@{}}T PP\\ (\%)\end{tabular}}} & \multicolumn{1}{c}{\textbf{\begin{tabular}[c]{@{}c@{}}T Sep\\ (\%)\end{tabular}}} & \multicolumn{1}{c}{\textbf{\begin{tabular}[c]{@{}c@{}}T Branch\\ (\%)\end{tabular}}} & \multicolumn{1}{c}{\textbf{\begin{tabular}[c]{@{}c@{}}T GSPP\\ (\%)\end{tabular}}} & \multicolumn{1}{c}{\textbf{Nodes}} & \multicolumn{1}{c}{\textbf{CG Its}} & \multicolumn{1}{c}{\textbf{Cols}} & \multicolumn{1}{c}{\textbf{Cuts}} \\ \hline
12-6-3-T & 0.00 & 164 & 88.3 & 9.1 & 0.4 & 0.9 & 0.0 & 45 & 1,356 & 7,215 & 567 \\
10-4-4-L & 0.00 & 6,068 & 95.3 & 4.4 & 0.1 & 0.1 & 0.0 & 81 & 5,667 & 35,916 & 382 \\
20-10-3-L & 0.02 & 11,089 & 28.5 & 26.8 & 4.1 & 17.8 & 0.0 & 15,421 & 49,256 & 131,616 & 24 \\
15-8-4-L & 0.00 & 968 & 69.8 & 27.8 & 0.4 & 1.1 & 0.0 & 49 & 1,417 & 10,629 & 861 \\
40-15-3-T & 0.02 & 11,022 & 68.9 & 12.2 & 1.4 & 3.3 & 1.3 & 903 & 6,390 & 94,066 & 20 \\ \hline
\textbf{Instance} & \multicolumn{11}{c}{\textbf{Branch \& Cut \& Price}} \\ \hline
\textbf{N-B-P-TW} & \multicolumn{1}{c}{\textbf{\begin{tabular}[c]{@{}c@{}}Gap\\ (\%)\end{tabular}}} & \multicolumn{1}{c}{\textbf{\begin{tabular}[c]{@{}c@{}}T\\ (s)\end{tabular}}} & \multicolumn{1}{c}{\textbf{\begin{tabular}[c]{@{}c@{}}T RMP\\ (\%)\end{tabular}}} & \multicolumn{1}{c}{\textbf{\begin{tabular}[c]{@{}c@{}}T PP\\ (\%)\end{tabular}}} & \multicolumn{1}{c}{\textbf{\begin{tabular}[c]{@{}c@{}}T Sep\\ (\%)\end{tabular}}} & \multicolumn{1}{c}{\textbf{\begin{tabular}[c]{@{}c@{}}T Branch\\ (\%)\end{tabular}}} & \multicolumn{1}{c}{\textbf{\begin{tabular}[c]{@{}c@{}}T GSPP\\ (\%)\end{tabular}}} & \multicolumn{1}{c}{\textbf{Nodes}} & \multicolumn{1}{c}{\textbf{CG Its}} & \multicolumn{1}{c}{\textbf{Cols}} & \multicolumn{1}{c}{\textbf{Cuts}} \\ \hline
12-6-3-T & 0.00 & 159 & 85.2 & 10.9 & 1.6 & 1.0 & 0.0 & 43 & 1,426 & 7,071 & 1,382 \\
10-4-4-L & 0.02 & 11,120 & 93.6 & 3.3 & 0.2 & 0.0 & 2.7 & 61 & 6,755 & 38,497 & 3,716 \\
20-10-3-L & 0.01 & 11,062 & 28.0 & 28.7 & 11.6 & 12.9 & 0.0 & 12,681 & 49,949 & 118,332 & 32,198 \\
15-8-4-L & 0.00 & 1,831 & 81.1 & 17.5 & 0.6 & 0.5 & 0.0 & 31 & 1,660 & 11,727 & 2,236 \\
40-15-3-T & 0.02 & 11,180 & 73.1 & 9.6 & 6.5 & 1.2 & 3.1 & 515 & 5,547 & 79,677 & 7,124
\end{tabular}
}
\end{table}
Table \ref{tab:extraInfo} provides a more detailed comparison of the proposed method variants for 5 instances that aim to be representative of the entire set of 50 instances. The first column indicates the instance, the second and third column recap the optimality gap and computational time spent given the time-limit of 3 hours and 18 minutes. The fourth to eighth columns indicate the percentage amount of time spent by the algorithm in the RMP, pricing problems (PPs), cut separation process (Sep), branching procedure (Branch), and the final GSPP model respectively.
{The pricing problems are solved in parallel on the four cores used.}
It should also be noticed that the branching time not only includes the selection of the branching candidate but also, the child nodes creation, which in the case of our algorithm, requires intensive data structure manipulation.
The number of nodes explored in the B\&B tree is displayed in the ninth column. The last three columns indicate the number of column generation iterations, generated columns and added cuts respectively.

The RMP takes most of the time for most of the instances, and the cut separation has an insignificant impact except when it is applied in every B\&B node. 
The time spent in branching procedures grows in accordance to both the size of the RMP and B\&B tree.
The short RMP solving times and large amount of columns generated for instance $20{-}10{-}3{-}L$ suggest that the RMP is easy to solve and the existence of many equivalent or similar solutions.
This increases the impact of other internal operations in the algorithm.
The number of B\&B nodes explored grows inversely to the amount of nodes where cutting is allowed.
It can be observed that the full \textit{branch-and-cut-price} performs more column generation iterations than the one with only cutting in the root node but it also requires longer computational times. 
As it can be observed, the time percentages do not sum exactly to 100\%. The remaining time accounts to diverse internal operations in the implementation which are not strictly linked to any of the main parts of the algorithm.
This also suggests that there is room for improvement in the implementation of the algorithm.

\begin{table}[]
\centering
\caption{Optimality gap to optimal or best known solution at the root node.}
\label{tab:rootNode}
\resizebox{\textwidth}{!}{%
\begin{tabular}{l|rr|rr|rr|rr}
\multicolumn{1}{c|}{\textbf{Instance}} & \multicolumn{2}{c|}{\textbf{MIP formulation}} & \multicolumn{2}{c|}{\textbf{Without cuts}} & \multicolumn{2}{c|}{\textbf{With sol-based cuts}} & \multicolumn{2}{c}{\textbf{With all cuts}} \\ \hline
\multicolumn{1}{c|}{\textbf{N-B-P-TW}} & \multicolumn{1}{c}{\textbf{Gap (\%)}} & \multicolumn{1}{c|}{\textbf{T (s)}} & \multicolumn{1}{c}{\textbf{Gap (\%)}} & \multicolumn{1}{c|}{\textbf{T (s)}} & \multicolumn{1}{c}{\textbf{Gap (\%)}} & \multicolumn{1}{c|}{\textbf{T (s)}} & \multicolumn{1}{c}{\textbf{Gap (\%)}} & \multicolumn{1}{c}{\textbf{T (s)}} \\ \hline
12-6-3-T & 2.03 & 0.3 & 0.36 & 2 & 0.22 & 5 & 0.22 & 87 \\
10-4-4-L & 4.69 & 0.1 & 0.22 & 15 & 0.15 & 102 & 0.15 & 1,137 \\
20-10-3-L & 0.36 & 0.2 & 0.03 & 2 & 0.03 & 2 & 0.03 & 143 \\
15-8-4-L & 2.11 & 0.2 & 0.25 & 19 & 0.12 & 43 & 0.12 & 499 \\
40-15-3-T & 7.89 & 1.2 & 0.02 & 280 & 0.02 & 285 & 0.02 & 685 \\ \hline
\multicolumn{1}{c|}{\textbf{Average 50 instances}} & \textbf{3.97} & \textbf{0.3} & \textbf{0.23} & \textbf{18} & \textbf{0.0865} & \textbf{25} & \textbf{0.0863} & \textbf{315}
\end{tabular}
}
\end{table}
The effectiveness of the aforementioned cut separation process is displayed in Table \ref{tab:rootNode}. The optimality gap of the LP solution at the root node is shown for the subset of five instances studied in detail together with the average across all 50 instances. The second column denotes the LP solution at the root node for the MIP formulation.
The third column refers to the presented method without adding any cuts whereas the fourth column considers the proposed cut separation procedure (Algorithm \ref{alg:cutSep}) based on solution values (i.e., \textit{sol-based}). 
This procedure only checks a subset of the valid inequalities which we believe that contains most, if not all, of the violated ones. In any case, we can find all violated inequalities by simple enumeration. This case, where all violated valid inequalities (i.e., \textit{all cuts}) are added, has also been tested and the results are shown in the last column.
The improvement in the lower bound is significant for the proposed methods where the cut separation is able to further improve it achieving an average gap of {0.09} \%. Adding all possible cuts only leads to an average improvement of {0.0002} \% in the bound. However, the algorithm requires {12} times more computational time to solve the root node. It is therefore decided to discard this variant of the separation procedure given the slow performance and the insignificant gain.

\begin{table}[]
\centering
\caption{Solving time comparison between the network flow problem and the \textit{branch-and-price} method.}
\label{tab:nfp}
\begin{tabular}{l|rr|r|r}
\multicolumn{1}{c|}{\textbf{Instance}} &
  \multicolumn{2}{c|}{\textbf{Graph size}} &
  \multicolumn{1}{c|}{\textbf{Network flow problem}} &
  \multicolumn{1}{c}{\textbf{Branch \& Price}} \\ \hline
\multicolumn{1}{c|}{\textbf{N-B-P-TW}} &
  \multicolumn{1}{c}{\textbf{Nodes}} &
  \multicolumn{1}{c|}{\textbf{Arcs}} &
  \multicolumn{1}{c|}{\textbf{T(s)}} &
  \multicolumn{1}{c}{\textbf{T(s)}} \\ \hline
4-3-3-L & 1,621 & 44,989    & 5.4  & 0.5  \\
5-3-3-L & 2,711 & 1,669,872 & 383 & 3.2  \\
6-3-3-L & 2,711 & 2,353,886 & 517 & 0.8  \\
6-3-4-L & 5,414 & 8,659,488 & 2,504 & 39.0
\end{tabular}
\end{table}
As mentioned in section \ref{sec:CG}, when having a pure shortest path as a pricing problem, solving the LP relaxation of the network flow problem gives the same bound as column generation on the GSPP but the network flow problem is expected to require more time and memory resources on instances with dense graphs. In order to verify that, the network flow problem has been solved for the first four instances which are considered among the \textit{easiest} ones from the entire set. The solving times of the network flow problem and the \textit{branch-and-price} method are compared in Table \ref{tab:nfp}. The complexity of the graph is shown by the high solving times for the network flow problem where the proposed model is between 10 and more than 500 times faster. The rest of instances have not been further analyzed as most of them were reaching the memory limit. The number of nodes and arcs for all the instances are documented in Table \ref{tab:InstGraphSize} in Appendix \ref{app:c}.

\begin{table}[]
\centering
\caption{Average performance of different branching strategies across all 50 instances using the \textit{branch-and-cut-and-price} method with only cutting allowed in the root node.}
\label{tab:varAlgSumm}
\begin{tabular}{c|cccc}
\textbf{Branching strategy} & \textbf{\begin{tabular}[c]{@{}c@{}}Best first on\\ single node\end{tabular}} & \textbf{\begin{tabular}[c]{@{}c@{}}Strong branching\\ on berthing time\end{tabular}} & \textbf{\begin{tabular}[c]{@{}c@{}}Best first on\\ berthing position\end{tabular}} & \textbf{\begin{tabular}[c]{@{}c@{}}Best first on\\ berthing time\end{tabular}} \\ \hline
\textbf{Average gap (\%)} & 0.058 & 0.011 & 0.009 & \textbf{0.005} \\
\textbf{Optimal solutions} & 16/50 & 37/50 & \textbf{39/50} & 38/50
\end{tabular}
\end{table}
Apart from the presented methods, slight variations have been tested which helped {to select} the best algorithm procedure.
For instance, we have tried to generate all columns \textit{a priori} without success. The complexity of the problem and the exponentially large numbers of columns make it intractable.
Regarding branching procedures, an alternative method of exploring the B\&B tree known as \textit{strong branching} has been tested. This strategy requires to select a number of candidates (between 5 and 10 in our case) and compute, or at least estimate, the lower bounds at the child nodes. For each candidate, a weighted sum of the child bounds is computed and the candidate with the best weighted sum is selected. In this case, a weight of 0.75 is set for the child with the lowest bound and a weight of 0.25 for the other child. This method has proven to create better branches 
{and, for example, is able to find an optimal solution to instance \textit{20-10-3-T} in less than 20 minutes. Nevertheless, the overall worse optimality gap and additional}
time consumed exploring more nodes has lead us to discard it. A different branching strategy has been tested where the branching is done on berthing positions instead of on berthing times at a port. 
{This strategy is able to solve all the instances from \cite{venturini2017a} to optimality and one more instance in overall than the presented method. However, the overall worse optimality gap indicates a poorer performance on the set of harder instances (see Table \ref{tab:varAlgs} in Appendix \ref{app:c}).}
In addition, a trivial branching on a single node has also been tested to compare the effectiveness of the proposed branching strategy. The solution values of the columns are added on the graph nodes of the respective paths, computing in this way the "usage" of each graph node. Then, the graph node whose value is closer to 0.5 (i.e., \textit{most fractional}) is the one selected to branch on. A summary of the performance of these alternative branching strategies is displayed in Table \ref{tab:varAlgSumm} and the results for all instances can be found
in Table \ref{tab:varAlgs} in Appendix \ref{app:c}.

\subsection{Cooperative game theory results}
The two methods for allocating the {costs} have been tested in the same set of instances.
\begin{table}[th]
\centering
\caption{Carrier ship share and priority for the instances.}
\label{tab:Carrier}
\begin{tabular}{c|ccc}
\textbf{Carrier} & \textbf{A} & \textbf{B} & \textbf{C} \\ \hline
\textbf{\% of ships} & 50 & 25 & 25 \\
\textbf{Priority} & 1 & 2 & 3
\end{tabular}
\end{table}
Three carriers $A, B$ and $C$ have been defined for all instances each of them with an assigned priority and a number of ships (see Table \ref{tab:Carrier}).
This priority is often given in accordance to the handling volume \citep{imai2003a}.
For instance, carrier A can be seen as a large carrier and often this translates in more power of decision and a higher priority at the port.
\begin{table}[th]
\centering
\caption{Terminal denomination.}
\label{tab:Terminals}
\begin{tabular}{c|cccc}
\textbf{Terminal} & \textbf{D} & \textbf{E} & \textbf{F} & \textbf{G} \\ \hline
\textbf{Visit position} & 1 & 2 & 3 & 4
\end{tabular}
\end{table}
{The terminal operator at each visited port is also a player in the game. In this study, all ships visit the ports in the same order as shown in Figure \ref{tab:Terminals}, but the game can equally be applied to instances with different visit orders. Thus, depending on the instance's number of ports, the game is formed by either 6 or 7 players. The number of possible coalitions is given by $2^{|\mathcal{P}|}$ where $|\mathcal{P}|$ denotes the number of players, and in this case, corresponds to 64 or 128 coalitions respectively.}

{The overall cooperative game is based on what we denote as the \textit{standalone solution}. This solution reflects the scenario 
where a single carrier negotiates with a single terminal at a time in order to decide the schedule for the carriers' vessels.
}
{We apply a greedy heuristic to compute this solution 
where we optimize and fix the schedule of a carrier's ships one port at a time. The sequence of carriers and ports used by the heuristic is given by the carrier's priority at port (see Table \ref{tab:Carrier}) and the position of the port visited (see Table \ref{tab:Terminals}). 
For example, 
assume ship 1 is carrier A's only ship and 
visits first port D and then port E where it has the highest priority. We then optimize the schedule of ship 1 for port D, fix the decisions and then optimize the schedule of ship 1 for port E.
Once the port visits of a carrier's ships are scheduled, the ships of the next carrier with the highest priority are scheduled within the remaining available berthing positions and time windows.}

This sequential planning approach resembles the actual procedure in some ports when the carriers book the port calls and they are assigned based on different priority schemes.
Due to the heuristic nature of the process, some of the carriers may not find a feasible schedule.
In order to avoid this and still ensure a fair comparison, the operational time windows of all berthing positions have been increased by 20\% in the tests performed in this section.

{The different coalitions $\mathcal{S}\subseteq \mathcal{P}$ can be classified into three groups, depending on the type of players forming it:
\begin{itemize}
    \item {Coalitions formed by carriers only. For the problem at hand this form of collaboration does not provide any planning advantage as the carriers require collaboration with the terminal operators to improve their planning. Therefore, the solution of this type of coalition corresponds to the standalone solution.}
    \item {Coalitions formed by terminal operators only. For the problem at hand this form of collaboration does not provide any planning advantage as the terminal operators require collaboration with the carriers to improve their planning. Therefore, the solution of this type of coalition corresponds to the standalone solution.}
    \item {Coalitions formed by both carriers and terminal operators. This type of coalition is the basis for the MPBAP. In order to compute a solution to a given coalition, we assume that the planning of all players that are not part of the coalition are kept fixed as in the standalone solution. Then, the MPBAP is solved for the coalition given the available berthing space and time at the terminals. Note that when more players are part of the coalition, fewer port calls of the standalone solution need to be fixed. In addition, it can be noticed that optimal solutions to coalitions formed by a single carrier and a single operator are equivalent to the standalone solution. We believe this minimal collaboration resembles the real-life port call booking process for carriers.}
\end{itemize}}

{The premise of fixing the port calls of non-collaborators, ensures that, in the worst case, the standalone solution is feasible for any coalition $\mathcal{S} \subseteq \mathcal{P}$.}
{As indicated in Section \ref{sec:costs}, carriers and terminal operators have different operational costs. This is also reflected in the characteristic function, where the costs of each player are measured differently. On one side, the carrier's cost comprise the fuel consumption costs, waiting time costs and half of the delay costs. On the other side, the terminal operator's cost comprise the handling costs and the remaining half of the delay costs. The process of quantifying the delay costs in this type of problems is complex and it has been decided to equally split the delay costs between carriers and terminal operators.}

{As shown in Section \ref{sec:instsRes}, the grand coalition scenario for some instances are not solved to optimality, but the proposed methods are able to find solutions within a very small gap in less than 6 minutes (see Tables \ref{tab:ResultsAlg5min} and \ref{tab:ResultsAlg5minHard}). It is assumed, that all the subcoalition scenarios are at most, as hard to solve as the grand coalition one and, therefore, a time limit of 5 minutes, and an additional 30 seconds to solve the GSPP, is set to solve each coalition scenario of the game.}

\begin{table}[]
\centering
\caption{Comparison of the two cost allocation methods across instances with three ports.}
\label{tab:compCGT3p}
\begin{tabular}{cc|r|rrr|rrr}
\multicolumn{1}{l}{} & \textbf{Player} & \multicolumn{1}{c|}{\textbf{Cost}} & \multicolumn{3}{c|}{\textbf{Shapley value}} & \multicolumn{3}{c}{\textbf{Equal profit method (EPM)}} \\ \cline{2-9} 
\multicolumn{1}{l}{} & $\mathcal{S}$ & \multicolumn{1}{c|}{$\vartheta(\mathcal{S})$} & \multicolumn{1}{c}{$f_i$} & \multicolumn{1}{c}{\textbf{\begin{tabular}[c]{@{}c@{}}\small Relative \\ \small savings (\%)\end{tabular}}} & \multicolumn{1}{c|}{\textbf{\begin{tabular}[c]{@{}c@{}} \small\% of total \\ \small costs\end{tabular}}} & \multicolumn{1}{c}{$f_i$} & \multicolumn{1}{c}{\textbf{\begin{tabular}[c]{@{}c@{}} \small Relative \\ \small savings (\%)\end{tabular}}} & \multicolumn{1}{c}{\textbf{\begin{tabular}[c]{@{}c@{}}\small \% of total \\ \small costs\end{tabular}}} \\ \hline
\multirow{3}{*}{\rotatebox[origin=c]{90}{\textbf{\small Carrier}}} & \textbf{A} & 485,261 & 478,464 & 1.5 & 39.1 & 471,601 & 3.3 & 38.4 \\
 & \textbf{B} & 247,053 & 242,800 & 2.2 & 20.0 & 240,685 & 3.1 & 19.8 \\
 & \textbf{C} & 251,485 & 246,914 & 2.4 & 20.9 & 245,092 & 3.1 & 20.7 \\ \hline
\multirow{4}{*}{\rotatebox[origin=c]{90}{\textbf{\tiny \ \ \ \ Terminal}}} & \textbf{D} & 88,635 & 79,042 & 10.7 & 7.0 & 86,058 & 3.4 & 7.5 \\
 & \textbf{E} & 75,221 & 70,789 & 6.9 & 6.4 & 73,059 & 3.3 & 6.7 \\
 & \textbf{F} & 75,296 & 71,546 & 5.8 & 6.6 & 73,060 & 3.2 & 6.8 \\ \hline
 & \textbf{\begin{tabular}[c]{@{}c@{}}Grand \\ coalition\end{tabular}} & \multicolumn{1}{l|}{1,189,555} &  & \multicolumn{1}{c}{} & \multicolumn{1}{c|}{} & \multicolumn{1}{c}{} & \multicolumn{1}{c}{} & \multicolumn{1}{c}{}
\end{tabular}
\end{table}

\begin{table}[]
\centering
\caption{Comparison of the two cost allocation methods across instances with four ports.}
\label{tab:compCGT4p}
\begin{tabular}{cc|r|rrr|rrr}
\multicolumn{1}{l}{} & \textbf{Player} & \multicolumn{1}{c|}{\textbf{Cost}} & \multicolumn{3}{c|}{\textbf{Shapley value}} & \multicolumn{3}{c}{\textbf{Equal profit method (EPM)}} \\ \cline{2-9} 
\multicolumn{1}{l}{} & $\mathcal{S}$ & \multicolumn{1}{c|}{$\vartheta(\mathcal{S})$} & \multicolumn{1}{c}{$f_i$} & \multicolumn{1}{c}{\textbf{\begin{tabular}[c]{@{}c@{}}\small Relative \\ \small savings (\%)\end{tabular}}} & \multicolumn{1}{c|}{\textbf{\begin{tabular}[c]{@{}c@{}} \small\% of total \\ \small costs\end{tabular}}} & \multicolumn{1}{c}{$f_i$} & \multicolumn{1}{c}{\textbf{\begin{tabular}[c]{@{}c@{}} \small Relative \\ \small savings (\%)\end{tabular}}} & \multicolumn{1}{c}{\textbf{\begin{tabular}[c]{@{}c@{}}\small \% of total \\ \small costs\end{tabular}}} \\ \hline
\multirow{3}{*}{\rotatebox[origin=c]{90}{\textbf{\small Carrier}}} & \textbf{A} & 352,695 & 347,232 & 1.7 & 37.7 & 339,450 & 4.0 & 36.8 \\
 & \textbf{B} & 187,320 & 181,789 & 3.1 & 20.5 & 179,622 & 4.3 & 20.2 \\
 & \textbf{C} & 192,634 & 187,148 & 3.1 & 20.9 & 184,727 & 4.4 & 20.6 \\ \hline
\multirow{4}{*}{\rotatebox[origin=c]{90}{\textbf{\small Terminal}}} & \textbf{D} & 46,264 & 40,908 & 11.7 & 4.6 & 44,464 & 4.0 & 5.0 \\
 & \textbf{E} & 44,079 & 38,207 & 14.1 & 4.1 & 42,215 & 4.4 & 4.6 \\
 & \textbf{F} & 66,496 & 60,718 & 8.5 & 6.3 & 63,629 & 3.9 & 6.6 \\
 & \textbf{G} & 60,346 & 55,765 & 7.6 & 5.9 & 57,659 & 3.9 & 6.1 \\ \cline{2-9} 
\multicolumn{1}{l}{} & \textbf{\begin{tabular}[c]{@{}c@{}}Grand \\ coalition\end{tabular}} & 911,766 &  & \multicolumn{1}{c}{} & \multicolumn{1}{c|}{} & \multicolumn{1}{c}{} & \multicolumn{1}{c}{} & \multicolumn{1}{c}{}
\end{tabular}
\end{table}

{All instances have a non-empty core, meaning that both efficient and stable solutions can be found in all scenarios.
}
{Tables \ref{tab:compCGT3p} and \ref{tab:compCGT4p} show the average cost allocations to each of the carriers and terminal operators across instances with 3 and 4 ports respectively (24 out of the 50 instances have 4 ports). For each allocation method we display three columns, (i) the first column indicates the cost allocation to the carrier when being part of the grand coalition, (ii) the second column computes the percentual savings compared with the player's standalone cost and (iii) the third column shows the percentage of the overall costs allocated to each player. Both allocation methods show that significant savings can be achieved by all of the players involved.}
In fact, player $A$, which in theory may be the least interested in engaging in such grand coalitions due to its high priority at all ports, achieves significant savings.
{The same applies to terminal D, which is the first one visited by the ships, and can benefit substantially by the overall better planning of the rest of terminals.}
{The differences in the allocation strategy used by the two methods are noticeable. The EPM tends to equalize the relative savings of all players whereas the Shapley value is prone to balance the absolute savings.}

{
As mentioned in Section \ref{sec1}, we conceive that the solution to the MPBAP and the cost allocation methods could be provided as a service by third party software companies. 
To establish the side payments in practice,
players would need to commit to the service for a pre-established period and agree with the potential savings estimated by the third party.
This is required in order to define
the number of participants in the collaboration.
Moreover, 
to measure the savings of the MPBAP solution, we need an estimate of the standalone costs of each player in a non-collaborative scenario. The third party could estimate this cost for each player using the current planning software and we assume that the player would agree to that estimate. Once the agreement is in place, the third party could be used as a proxy for the side payments, which could be performed on a regular basis. Based on the actual costs, each player would need to make or receive a payment to align with the projected savings that the player has agreed to. 
}

{
Similar collaboration is already taking place in real-life in tramp shipping where it is common that a number of ship owners place their ships in a shipping pool under the control of a pool administrator (the analog to the third party coordinator) that takes over most of the business decisions regarding the ships and is responsible for distributing earnings to ship owners \citep{packard1995shipping, wen2019a}.
}

\section{Conclusions and future work}\label{sec:conclusion}

A novel solution method based on a GSPP formulation has been presented for the MPBAP. The method exploits a graph formulation for defining the berthing plan of a ship along its route.
This, combined together with delayed column generation, additional valid inequalities and symmetry breaking constraints results in an efficient algorithm able to find optimal or near-optimal solutions to wide range of instances outperforming the capacity of commercial solvers.

In addition, the graph formulation adds flexibility as many additional constraints can be easily integrated with simple alterations in the graph.
For instance, a finer discretization of the berthing positions would allow to approximate the continuous version of the MPBAP better. 
Considering a continuous berth is a more realistic approach and allows to increase the usage of the quay.
Transhipments are also an important aspect of the operations at port and the fulfillment of them are crucial in some cases (e.g., when transporting perishable food). Our model could eventually account for that by limiting the time window of the ships involved in the transhipment and penalizing late arrivals of the incoming ship or too early departures of the outgoing ship.
Nevertheless, this could be better modelled if the relative arrival and departure times are considered. That case is harder to incorporate in the presented model and it would require additional constraints for each transhipment in the RMP.
The transit times between ports could be further improved by considering the time needed to enter and leave the port which is usually performed at a slower speed \citep{reinhardt2016a}.

The instances solved reflect the size of real-life scenarios to a large extent. However, some of the instance parameters could be further improved. This comprises improving the size of the vessel time windows, having different routes for the ships, different amount of berthing positions per port and different ship types.

{Alternative branching methods have also shown great potential, especially branching on berthing positions as opposed to branching on berthing times. A natural next step would be to explore a branching method that combines both of them. For instance, one could test both methods simultaneously when branching and select the one with better bounds at the children nodes.}

A natural extension of the problem could be to integrate the berth allocation with the quay crane assignment problem (QAP). Studies such as \cite{iris2015a} and \cite{iris2017a} have already shown the effectiveness of heuristic and exact methods based on a GSPP formulation for the integrated problem in one terminal.

Last but not least, the benefits for both ship carriers and terminal operators are verified defining a cooperative game and using {cost} allocation methods to distribute the {costs} of such collaboration fairly. 
{The results of} the cooperative game strengthen the viability of such a decision tool and can encourage carriers and port operators to participate in collaborative schemes.

\ACKNOWLEDGMENT{%
This work is partially supported by The Danish Maritime Fund.
The authors are grateful to Dr. Cagatay Iris for the useful discussion about the model and instances
{and to Prof. Harilaos Psaraftis and Dr. Thalis Zis for their valuable insights.}
{Also, the authors are thankful to the associate editor and to two anonymous reviewers for their constructive remarks that have helped improve the manuscript.}
}

%
\begin{APPENDICES}{}

\section{Reduced cost computation including valid inequalities (\ref{eq:VIrmp})}\label{app:a}

We denote $\beta^{n,p,b}_{t_1,t_2}$ to the dual variable of constraint (\ref{eq:VIrmp}) for ship $n\in N$, port $p\in P$, berth $b\in B_p$ and times $t_1,t_2 \in [s^{p,b};e^{p,b}],t_1 < t_2$ and let $\Bar{\beta}^{n,p,b}_{t_1,t_2}$ be its value for the RMP solution.
Let $w(p,b,t) \in O$ be the graph node related to berthing at port $p \in P$ at position $b \in B_p$ at time $t \in [s^{p,b};e^{p,b}]$ and let $\Phi = \{(n,p,b,t_1,t_2) \}$ be the set of constraints (\ref{eq:VIrmp}) added to the RMP denoted by $(n,p,b,t_1,t_2)$ elements where $n\in N$, $p\in P$, $b\in B_p$ and $t_1,t_2 \in [s^{p,b};e^{p,b}],t_1 < t_2$.
Additionally, let $\Phi(k,p,b,t) \subseteq \Phi$ be the set of valid inequalities that include arcs from the graph node $w(p,b,t)$ for a given ship $k$. By definition the range of nodes for each ship within a valid inequality differ if $k=n$ or $k\neq n$. Therefore we denote $\Phi_{k=n}(k,p,b,t), \Phi_{k\neq n}(k,p,b,t) \subseteq \Phi(k,p,b,t)$ to the subset of cuts $(n,p,b,t_1,t_2)$ where $k=n$ and $k\neq n$ respectively, that together form the entire set $\Phi(k,p,b,t) = \Phi_{k=n}(k,p,b,t) \cup  \Phi_{k\neq n}(k,p,b,t)$ and are defined mathematically as follows:
\begin{equation*}
    \Phi_{k=n}(k,p,b,t) = \Big\{ (n,p,b,t_1,t_2) | k=n, w(p,b,t) \in \bigcup_{t \in [t_1;t_2]} C(k,p,b,t) \Big\}
\end{equation*}
\begin{equation*}
    \Phi_{k\neq n}(k,p,b,t) = \Big\{ (n,p,b,t_1,t_2) | k \neq n, w(p,b,t) \in C(k,p,b,t_1) \cap C(k,p,b,t_2) \Big\}
\end{equation*}
The computation of the reduced cost $\hat{c}_j$ for column $j$ of ship $k \in N$ is updated as follows:
\begin{equation*}
    \hat{c}_j = c_j - \Big(\sum_{(p,b,t) \in \Lambda_j} (\sum_{t' \in [t;t+h^{p,b}_k)} \Bar{\mu}_{p,b,t'}) - (\sum_{(n,p,b,t_1,t_2) \in \Phi(k,p,b,t)} \Bar{\beta}^{n,p,b}_{t_1,t_2}) \Big) - \Bar{\alpha}_k
\end{equation*}

\section{Adaption of the proposed valid inequality considering berth types}\label{app:b}

\begin{proposition}\label{prop:2}
Given two time instants $t_1,t_2 \in [s^{p,k},e^{p,k})$ where $t_1 < t_2$ and a port $p\in P$, berth type $k\in K_p$ and ship $n\in N$, the following is a valid inequality:
\[
\sum_{ u \in \bigcup_{t\in [t_1;t_2]} C(n,p,k,t)}
\sum_{w \in \delta^+_n(u)} x^n_{u,w} + 
\sum_{m\in N \backslash\{n\}}\sum_{u \in C(m,p,k,t_1)\cap C(m,p,k,t_2)}\sum_{w \in \delta^+_m(u)} x^m_{u,w} \leqslant \beta^k 
\]
\end{proposition}
\proof{Proof.}
Constraint (\ref{GSPP:pkt}) has been adapted from constraint (\ref{GSPP:pbt}) which is a direct translation from constraint (\ref{eq:Nconf}) from the network-flow formulation. Therefore, constraint (\ref{GSPP:pkt}) can be formulated as a network-flow problem constraint as follows:
\begin{equation*}
    \sum_{m \in N} \sum_{i \in C(m,p,k,t)} \sum_{j \in \delta^+_m(i)} x^m_{i,j} \leq \beta^k \quad \forall p\in P, k\in K_p, t \in [s^{p,k};e^{p,k})
\end{equation*}
where the $C(m,p,k,t)$ defines the set of nodes for ship $m\in N$ that are in conflict with time $t\in [s^{p,k};e^{p,k})$ in berth type $k\in K_p$ of port $p\in P$. Based on this definition, the intersection set $C(m,p,k,t_1) \cap C(m,p,k,t_2)$ directly defines the set of nodes for ship $m$ that are in conflict with both time instants $t_1$ and $t_2$. Constraint (\ref{GSPP:pkt}) indicates that at most $\beta^k$ (i.e., number of berths of type $k \in B_k$) arcs can be chosen out of the nodes from the sets $C(m,p,k,t)$ of all ships $m \in N$ and, therefore, the same applies to the intersection sets $C(m,p,k,t_1) \cap C(m,p,k,t_2)$. 
Based on the premise that each ship can only berth in one position, we can relax the requirement of being in conflict with both $t_1$ and $t_2$ for a single ship $n$ and only require it to be in conflict with $t_1$ or $t_2$. In practice, this means, on one hand, that if ship $n$ berths at a period covering $t_1$ or $t_2$, then, at most $\beta^k-1$ ships $m \in N \backslash \{n\}$ can berth in a period covering both $t_1$ and $t_2$. On the other hand, if $\beta^k$ ships $m \in N \backslash \{n\}$ are berthing at times whose periods cover $t_1$ and $t_2$, then ship $n$ is not able to berth at a period covering $t_1$ or $t_2$.   
The relaxed node interval for ship $n$ can be defined as the union of $C(n,p,k,t_1)\cup C(n,p,k,t_2)$. Considering these node sets, we can define the following valid inequality:
\begin{multline*}
    \sum_{u \in C(n,p,k,t_1)\cup C(n,p,k,t_2)}\sum_{w \in \delta^+_n(u)} x^n_{u,w} + 
\sum_{m\in N \backslash\{n\}}\sum_{u \in C(m,p,k,t_1)\cap C(m,p,k,t_2)}\sum_{w \in \delta^+_m(u)} x^m_{u,w} \leqslant \beta^k  \\
\forall p\in P, k\in K_p, n\in N,t_1,t_2 \in [s^{p,k},e^{p,k}),t_1 < t_2 
\end{multline*}

Based on the assumption that a berthing period cannot be discontinued, the intersection set $C(m,p,k,t_1) \cap C(m,p,k,t_2)$ for any ship is not only in conflict with times $t_1$ and $t_2$ but with all the time instants in the period $[t_1;t_2]$. Therefore the interval for ship $n$ can be expanded to the union of $C(n,p,k,t)$ sets for all time instants $t \in [t_1;t_2]$ and the resulting valid inequality can be formulated as follows:

\begin{multline*}
\sum_{u \in \bigcup_{t\in [t_1;t_2]} C(n,p,k,t)}\sum_{w \in \delta^+_n(u)} x^n_{u,w} + 
\sum_{m\in N \backslash\{n\}}\sum_{u \in C(m,p,k,t_1)\cap C(m,p,k,t_2)}\sum_{w \in \delta^+_m(u)} x^m_{u,w} \leqslant \beta^k  \\
\forall p\in P, k\in K_p, n\in N,t_1,t_2 \in [s^{p,k},e^{p,k}),t_1 < t_2 
\end{multline*}
\Halmos
\endproof

\section{Additional computational results}\label{app:c}

\begin{table}[]
\centering
\caption{Number of nodes and arcs in graph $G$ for each of the instances. An horizontal line is used to indicate the separation between the set of benchmark instances by \cite{venturini2017a} and the newly generated set of harder instances.}
\label{tab:InstGraphSize}
\scalebox{1}{
\begin{tabular}{l|rr|l|rr}
\multicolumn{1}{c|}{\textbf{Instance}} &
  \multicolumn{2}{c|}{\textbf{Graph size}} &
  \multicolumn{1}{c|}{\textbf{Instance}} &
  \multicolumn{2}{c|}{\textbf{Graph size}} \\ \hline
\multicolumn{1}{c|}{\textbf{N-B-P-TW}} &
  \multicolumn{1}{c}{\textbf{Nodes}} &
  \multicolumn{1}{c|}{\textbf{Arcs}} &
  \multicolumn{1}{c|}{\textbf{N-B-P-TW}} &
  \multicolumn{1}{c}{\textbf{Nodes}} &
  \multicolumn{1}{c}{\textbf{Arcs}} \\ \hline
4-3-3-L   & 1,621  & 44,989     & 4-4-4-T   & 3,948  & 1,239,853   \\
5-3-3-L   & 2,711  & 1,669,872  & 5-4-4-T   & 3,948  & 1,498,010   \\
6-3-3-L   & 2,711  & 2,353,886  & 6-4-4-T   & 3,948  & 1,318,353   \\
6-3-4-L   & 5,414  & 8,659,488  & 12-5-3-T  & 3,217  & 2,572,191   \\
10-4-4-L  & 7,218  & 25,069,050 & 12-6-3-T  & 3,860  & 3,586,436   \\
10-4-3-L  & 3,614  & 6,583,065  & 10-5-4-T  & 5,722  & 8,069,850   \\
4-4-4-L   & 6,418  & 5,631,290  & 15-10-4-T & 8,870  & 25,003,591  \\
5-4-4-L   & 6,418  & 8,436,182  & 20-10-3-T & 5,120  & 11,456,682  \\
6-4-4-L   & 6,418  & 10,062,500 & 20-12-3-T & 6,826  & 21,299,213  \\ \cline{4-6} 
12-5-3-L  & 4,517  & 12,072,289 & 25-12-3-L & 7,226  & 65,709,632  \\
10-6-3-L  & 5,420  & 15,698,842 & 25-12-3-T & 6,826  & 26,586,128  \\
11-6-3-L    & 5,420  & 15,640,325 & 12-5-4-L  & 9,022  & 45,812,954  \\
12-6-3-L    & 5,420  & 17,117,528 & 12-5-4-T  & 5,722  & 9,550,769   \\
10-5-4-L    & 9,022  & 38,507,345 & 30-12-3-L & 7,226  & 78,917,800  \\
15-10-3-L & 5,420  & 21,596,125 & 30-12-3-T & 6,826  & 31,947,951  \\
15-12-3-L   & 7,226  & 37,072,022 & 20-12-4-L & 15,638 & 171,445,244 \\
15-10-4-L & 12,430 & 76,196,427 & 20-12-4-T & 11,158 & 56,721,679  \\
20-10-3-L & 5,420  & 28,790,640 & 15-8-4-L  & 14,434 & 144,240,102 \\
20-12-3-L & 7,226  & 52,608,061 & 15-8-4-T  & 9,154  & 29,624,264  \\
4-3-3-T   & 1,621  & 229,372    & 25-12-4-L & 15,638 & 214,550,093 \\
5-3-3-T   & 2,711  & 2,696,410  & 25-12-4-T & 11,158 & 70,841,345  \\
6-3-3-T   & 2,711  & 2,538,662  & 30-15-4-L & 16,541 & 300,541,288 \\
6-3-4-T   & 3,464  & 2,307,063  & 30-15-4-T & 11,801 & 98,992,726  \\
10-4-4-T  & 4,618  & 6,633,974  & 40-15-3-L & 8,129  & 133,728,684 \\
10-4-3-T  & 2,614  & 1,825,953  & 40-15-3-T & 7,679  & 54,138,993 
\end{tabular}
}
\end{table}

\begin{table}[]
\centering
\caption{Results of solution methods with alternative branching strategies. The underlying algorithm is a \textit{branch-and-cut-and-price} where cutting is only allowed at the root node. "*" means the time limit of 3 hours and 18 minutes has been reached.}
\label{tab:varAlgs}
\scalebox{0.7}{
\begin{tabular}{l|rrrr|rrrr|rrrr}
\textbf{Instance} & \multicolumn{4}{c|}{\textbf{Best first on single graph node}} & \multicolumn{4}{c|}{\textbf{Strong branching on berthing time}} & \multicolumn{4}{c}{\textbf{Best first on berthing position}} \\ \hline
\textbf{N-B-P-TW} & \textbf{LB} & \textbf{Z} & \textbf{Gap (\%)} & \textbf{T (s)} & \textbf{LB} & \textbf{Z} & \textbf{Gap (\%)} & \textbf{T (s)} & \textbf{LB} & \textbf{Z} & \textbf{Gap (\%)} & \textbf{T (s)} \\ \hline
4-3-3-L & 296,600 & 296,600 & 0.00 & 0.2 & 296,600 & 296,600 & 0.00 & 0.2 & 296,600 & 296,600 & 0.00 & 0.2 \\
5-3-3-L & 394,300 & 394,300 & 0.00 & 8.9 & 394,300 & 394,300 & 0.00 & 9.2 & 394,300 & 394,300 & 0.00 & 8.1 \\
6-3-3-L & 421,720 & 421,720 & 0.00 & 0.5 & 421,720 & 421,720 & 0.00 & 0.6 & 421,720 & 421,720 & 0.00 & 0.5 \\
6-3-4-L & 647,149 & 647,480 & 0.05 & * & 647,480 & 647,480 & 0.00 & 1,024.4 & 647,480 & 647,480 & 0.00 & 93.7 \\
10-4-4-L & 1,053,050 & 1,054,500 & 0.14 & * & 1,054,250 & 1,054,500 & 0.02 & * & 1,054,500 & 1,054,500 & 0.00 & 9,628.2 \\
10-4-3-L & 697,057 & 698,100 & 0.15 & * & 698,100 & 698,100 & 0.00 & 443.5 & 698,100 & 698,100 & 0.00 & 173.0 \\
4-4-4-L & 405,120 & 405,120 & 0.00 & 0.9 & 405,120 & 405,120 & 0.00 & 0.7 & 405,120 & 405,120 & 0.00 & 0.6 \\
5-4-4-L & 500,600 & 500,600 & 0.00 & 1.2 & 500,600 & 500,600 & 0.00 & 0.9 & 500,600 & 500,600 & 0.00 & 0.9 \\
6-4-4-L & 599,780 & 599,980 & 0.03 & * & 599,980 & 599,980 & 0.00 & 26.9 & 599,980 & 599,980 & 0.00 & 7.4 \\
12-5-3-L & 829,897 & 830,440 & 0.07 & * & 830,440 & 830,440 & 0.00 & 646.1 & 830,440 & 830,440 & 0.00 & 143.3 \\
10-6-3-L & 680,550 & 680,600 & 0.01 & * & 680,600 & 680,600 & 0.00 & 34.4 & 680,600 & 680,600 & 0.00 & 15.1 \\
11-6-3-L & 745,960 & 746,220 & 0.03 & * & 746,220 & 746,220 & 0.00 & 52.4 & 746,220 & 746,220 & 0.00 & 16.9 \\
12-6-3-L & 809,093 & 810,040 & 0.09 & * & 809,840 & 809,840 & 0.00 & 301.7 & 809,840 & 809,840 & 0.00 & 54.3 \\
10-5-4-L & 1,026,521 & 1,028,320 & 0.17 & * & 1,028,320 & 1,028,320 & 0.00 & 2,472.5 & 1,028,320 & 1,028,320 & 0.00 & 5,991.3 \\
15-10-3-L & 1,006,000 & 1,006,200 & 0.02 & * & 1,006,200 & 1,006,200 & 0.00 & 59.8 & 1,006,200 & 1,006,200 & 0.00 & 4.8 \\
15-12-3-L & 1,002,800 & 1,002,800 & 0.00 & 2.0 & 1,002,800 & 1,002,800 & 0.00 & 1.8 & 1,002,800 & 1,002,800 & 0.00 & 2.3 \\
15-10-4-L & 1,459,600 & 1,459,600 & 0.00 & 339.9 & 1,459,600 & 1,459,600 & 0.00 & 62.6 & 1,459,600 & 1,459,600 & 0.00 & 14.3 \\
20-10-3-L & 1,344,437 & 1,344,800 & 0.03 & * & 1,344,583 & 1,344,800 & 0.02 & * & 1,344,800 & 1,344,800 & 0.00 & 163.3 \\
20-12-3-L & 1,336,400 & 1,336,400 & 0.00 & 3.4 & 1,336,400 & 1,336,400 & 0.00 & 2.9 & 1,336,400 & 1,336,400 & 0.00 & 2.2 \\
4-3-3-T & 318,440 & 318,440 & 0.00 & 0.5 & 318,440 & 318,440 & 0.00 & 0.3 & 318,440 & 318,440 & 0.00 & 0.3 \\
5-3-3-T & 405,240 & 405,240 & 0.00 & 2.5 & 405,240 & 405,240 & 0.00 & 3.6 & 405,240 & 405,240 & 0.00 & 7.4 \\
6-3-3-T & 510,920 & 510,920 & 0.00 & 1.5 & 510,920 & 510,920 & 0.00 & 1.1 & 510,920 & 510,920 & 0.00 & 0.8 \\
6-3-4-T & 993,460 & 993,460 & 0.00 & 2.5 & 993,460 & 993,460 & 0.00 & 1.8 & 993,460 & 993,460 & 0.00 & 1.2 \\
10-4-4-T & 1,660,640 & 1,660,640 & 0.00 & 1,193.6 & 1,660,640 & 1,660,640 & 0.00 & 329.9 & 1,660,640 & 1,660,640 & 0.00 & 132.9 \\
10-4-3-T & 1,022,200 & 1,022,200 & 0.00 & 49.1 & 1,022,200 & 1,022,200 & 0.00 & 56.8 & 1,022,200 & 1,022,200 & 0.00 & 9.1 \\
4-4-4-T & 442,600 & 442,600 & 0.00 & 1.3 & 442,600 & 442,600 & 0.00 & 1.0 & 442,600 & 442,600 & 0.00 & 0.8 \\
5-4-4-T & 575,350 & 576,010 & 0.11 & * & 576,010 & 576,010 & 0.00 & 26.7 & 576,010 & 576,010 & 0.00 & 9.6 \\
6-4-4-T & 651,480 & 654,040 & 0.32 & * & 653,560 & 653,560 & 0.00 & 39.0 & 653,560 & 653,560 & 0.00 & 10.6 \\
12-5-3-T & 830,067 & 830,440 & 0.04 & * & 830,440 & 830,440 & 0.00 & 355.7 & 830,440 & 830,440 & 0.00 & 100.0 \\
12-6-3-T & 817,602 & 819,040 & 0.15 & * & 818,840 & 818,840 & 0.00 & 369.4 & 818,840 & 818,840 & 0.00 & 141.1 \\
10-5-4-T & 1,143,431 & 1,144,160 & 0.06 & * & 1,144,160 & 1,144,160 & 0.00 & 218.5 & 1,144,160 & 1,144,160 & 0.00 & 53.0 \\
15-10-4-T & 1,596,310 & 1,597,100 & 0.05 & * & 1,597,100 & 1,597,100 & 0.00 & 66.0 & 1,597,100 & 1,597,100 & 0.00 & 21.1 \\
20-10-3-T & 1,629,000 & 1,629,500 & 0.03 & * & 1,629,500 & 1,629,500 & 0.00 & 1,138.7 & 1,629,500 & 1,629,500 & 0.00 & 408.7 \\
20-12-3-T & 1,606,500 & 1,606,500 & 0.00 & 50.5 & 1,606,500 & 1,606,500 & 0.00 & 80.1 & 1,606,500 & 1,606,500 & 0.00 & 38.3 \\ \hline
\textbf{Average} &  &  & \textbf{0.0461} &  &  &  & \textbf{0.0012} &  &  &  & \textbf{0.0000} &  \\ \hline
25-12-3-L & 1,677,200 & 1,677,400 & 0.01 & * & 1,677,234 & 1,677,400 & 0.01 & * & 1,677,233 & 1,677,400 & 0.01 & * \\
25-12-3-T & 2,046,933 & 2,047,000 & 0.00 & * & 2,046,933 & 2,047,000 & 0.00 & * & 2,046,933 & 2,047,000 & 0.00 & * \\
12-5-4-L & 1,240,248 & 1,245,160 & 0.39 & * & 1,243,262 & 1,245,160 & 0.15 & * & 1,243,019 & 1,245,160 & 0.17 & * \\
12-5-4-T & 1,398,454 & 1,404,640 & 0.41 & * & 1,402,329 & 1,404,520 & 0.14 & * & 1,403,550 & 1,405,580 & 0.05 & * \\
30-12-3-L & 2,016,178 & 2,016,600 & 0.02 & * & 2,016,178 & 2,016,600 & 0.02 & * & 2,016,245 & 2,016,600 & 0.02 & * \\
30-12-3-T & 2,491,406 & 2,492,500 & 0.04 & * & 2,491,490 & 2,492,600 & 0.04 & * & 2,491,437 & 2,492,600 & 0.04 & * \\
20-12-4-L & 1,943,545 & 1,943,800 & 0.01 & * & 1,943,800 & 1,943,800 & 0.00 & 778.0 & 1,943,800 & 1,943,800 & 0.00 & 635.2 \\
20-12-4-T & 3,113,085 & 3,113,170 & 0.00 & * & 3,113,170 & 3,113,170 & 0.00 & 75.1 & 3,113,170 & 3,113,170 & 0.00 & 40.8 \\
15-8-4-L & 1,495,243 & 1,497,100 & 0.12 & * & 1,497,000 & 1,497,000 & 0.00 & 3,228.6 & 1,497,000 & 1,497,000 & 0.00 & 4,309.9 \\
15-8-4-T & 1,655,075 & 1,656,040 & 0.06 & * & 1,656,040 & 1,656,040 & 0.00 & 4,356.3 & 1,656,040 & 1,656,040 & 0.00 & 792.3 \\
25-12-4-L & 2,439,300 & 2,440,500 & 0.05 & * & 2,439,810 & 2,440,500 & 0.03 & * & 2,439,585 & 2,440,800 & 0.04 & * \\
25-12-4-T & 3,705,997 & 3,707,790 & 0.04 & * & 3,707,390 & 3,707,390 & 0.00 & 578.0 & 3,707,390 & 3,707,390 & 0.00 & 6,483.4 \\
30-15-4-L & 2,918,409 & 2,918,800 & 0.01 & * & 2,918,510 & 2,918,600 & 0.00 & * & 2,918,464 & 2,918,800 & 0.00 & * \\
30-15-4-T & 3,274,228 & 3,278,860 & 0.12 & * & 3,275,538 & 3,278,480 & 0.08 & * & 3,275,004 & 3,278,280 & 0.10 & * \\
40-15-3-L & 2,689,060 & 2,689,200 & 0.01 & * & 2,689,060 & 2,689,200 & 0.01 & * & 2,689,060 & 2,689,200 & 0.01 & * \\
40-15-3-T & 3,329,567 & 3,330,500 & 0.02 & * & 3,329,567 & 3,330,500 & 0.02 & * & 3,329,567 & 3,330,300 & 0.02 & * \\ \hline
\textbf{Average} &  &  & \textbf{0.0824} &  &  &  & \textbf{0.0314} &  &  &  & \textbf{0.0288} & \\
\multicolumn{2}{l}{\textbf{Optimal solutions}} & \textbf{} & \textbf{16/50} & \textbf{} & \textbf{} & \textbf{} & \textbf{37/50} & \textbf{} & \textbf{} & \textbf{} & \textbf{39/50} & \textbf{} 
\end{tabular}
}
\end{table}

\end{APPENDICES}
%
%


\bibliographystyle{informs2014trsc} 
\bibliography{sample.bib} 

\begin{thebibliography}{80}
\providecommand{\natexlab}[1]{#1}
\providecommand{\url}[1]{\texttt{#1}}
\providecommand{\urlprefix}{URL }

\bibitem[{{APM Terminals}(2021)}]{apmt}
{APM Terminals}, 2021 \emph{{APM Terminals}}.
  \url{https://www.apmterminals.com/}, accessed: 2021-09-15.

\bibitem[{Bektaş et~al.(2019)Bektaş, Ehmke, Psaraftis, \protect\BIBand{}
  Puchinger}]{bekta2019a}
Bektaş T, Ehmke JF, Psaraftis HN, Puchinger J, 2019 \emph{The role of
  operational research in green freight transportation}. \emph{European Journal
  of Operational Research} 274(3):807--823,
  \urlprefix\url{http://dx.doi.org/10.1016/j.ejor.2018.06.001}.

\bibitem[{Bezanson et~al.(2017)Bezanson, Edelman, Karpinski, \protect\BIBand{}
  Shah}]{bezanson2017a}
Bezanson J, Edelman A, Karpinski S, Shah VB, 2017 \emph{Julia: A fresh approach
  to numerical computing}. \emph{Siam Review} 59(1):65--98,
  \urlprefix\url{http://dx.doi.org/10.1137/141000671}.

\bibitem[{Bierwirth \protect\BIBand{} Meisel(2015)}]{bierwirth2015a}
Bierwirth C, Meisel F, 2015 \emph{A follow-up survey of berth allocation and
  quay crane scheduling problems in container terminals}. \emph{European
  Journal of Operational Research} 244(3):12689, 675--689,
  \urlprefix\url{http://dx.doi.org/10.1016/j.ejor.2014.12.030}.

\bibitem[{Brouer, Pisinger, \protect\BIBand{} Spoorendonk(2011)}]{brouer2011a}
Brouer BD, Pisinger D, Spoorendonk S, 2011 \emph{Liner shipping cargo
  allocation with repositioning of empty containers}. \emph{INFOR Journal}
  49(2):109--124, \urlprefix\url{http://dx.doi.org/10.3138/infor.49.2.109}.

\bibitem[{Buhrkal et~al.(2011)Buhrkal, Zuglian, Røpke, Larsen,
  \protect\BIBand{} Lusby}]{buhrkal2011a}
Buhrkal KF, Zuglian S, Røpke S, Larsen J, Lusby RM, 2011 \emph{Models for the
  discrete berth allocation problem: A computational comparison}.
  \emph{Transportation Research. Part E: Logistics and Transportation Review}
  47(4):461--473, \urlprefix\url{http://dx.doi.org/10.1016/j.tre.2010.11.016}.

\bibitem[{Carlo, Vis, \protect\BIBand{} Roodbergen(2014)}]{carlo2014a}
Carlo HJ, Vis IF, Roodbergen KJ, 2014 \emph{Transport operations in container
  terminals: Literature overview, trends, research directions and
  classification scheme}. \emph{European Journal of Operational Research}
  236(1):1--13, \urlprefix\url{http://dx.doi.org/10.1016/j.ejor.2013.11.023}.

\bibitem[{Chang et~al.(2012)Chang, Tongzon, Luo, \protect\BIBand{}
  Lee}]{chang2012a}
Chang YT, Tongzon J, Luo M, Lee PTW, 2012 \emph{Estimation of optimal handling
  capacity of a container port: An economic approach}. \emph{Transport Reviews}
  32(2):241--258,
  \urlprefix\url{http://dx.doi.org/10.1080/01441647.2011.644346}.

\bibitem[{Cheong et~al.(2010)Cheong, Tan, Liu, \protect\BIBand{}
  Lin}]{cheong2010a}
Cheong CY, Tan KC, Liu DK, Lin CJ, 2010 \emph{Multi-objective and prioritized
  berth allocation in container ports}. \emph{Annals of Operations Research}
  180(1):63--103, \urlprefix\url{http://dx.doi.org/10.1007/s10479-008-0493-0}.

\bibitem[{Cordeau et~al.(2005)Cordeau, Laporte, Legato, \protect\BIBand{}
  Moccia}]{cordeau2005a}
Cordeau JF, Laporte G, Legato P, Moccia L, 2005 \emph{Models and tabu search
  heuristics for the berth-allocation problem}. \emph{Transportation Science}
  39(4):526--538, \urlprefix\url{http://dx.doi.org/10.1287/trsc.1050.0120}.

\bibitem[{Cormen, Leiserson, \protect\BIBand{} Rivest(1996)}]{cormen1996a}
Cormen T, Leiserson C, Rivest R, 1996 \emph{Introduction to algorithms} (MIT
  Press,), ISBN 0070131430.

\bibitem[{Corry \protect\BIBand{} Bierwirth(2019)}]{corry2019a}
Corry P, Bierwirth C, 2019 \emph{The berth allocation problem with channel
  restrictions}. \emph{Transportation Science} 53(3):708--727,
  \urlprefix\url{http://dx.doi.org/10.1287/trsc.2018.0865}.

\bibitem[{Dantzig \protect\BIBand{} Wolfe(1960)}]{dantzig1960a}
Dantzig G, Wolfe P, 1960 \emph{Decomposition principle for linear-programs}.
  \emph{Operations Research} 8(1):101--111,
  \urlprefix\url{http://dx.doi.org/10.1287/opre.8.1.101}.

\bibitem[{Du et~al.(2015)Du, Chen, Lam, Xu, \protect\BIBand{} Cao}]{du2015a}
Du Y, Chen Q, Lam JSL, Xu Y, Cao JX, 2015 \emph{Modeling the impacts of tides
  and the virtual arrival policy in berth allocation}. \emph{Transportation
  Science} 49(4):939--956,
  \urlprefix\url{http://dx.doi.org/10.1287/trsc.2014.0568}.

\bibitem[{Du et~al.(2011)Du, Chen, Quan, Long, \protect\BIBand{}
  Fung}]{du2011a}
Du Y, Chen Q, Quan X, Long L, Fung RY, 2011 \emph{Berth allocation considering
  fuel consumption and vessel emissions}. \emph{Transportation Research Part E:
  Logistics and Transportation Review} 47(6):1021--1037,
  \urlprefix\url{http://dx.doi.org/10.1016/j.tre.2011.05.011}.

\bibitem[{Dulebenets(2018)}]{dulebenets2018b}
Dulebenets MA, 2018 \emph{A comprehensive multi-objective optimization model
  for the vessel scheduling problem in liner shipping}. \emph{International
  Journal of Production Economics} 196:293--318,
  \urlprefix\url{http://dx.doi.org/10.1016/j.ijpe.2017.10.027}.

\bibitem[{Dulebenets(2019)}]{dulebenets2019b}
Dulebenets MA, 2019 \emph{Minimizing the total liner shipping route service
  costs via application of an efficient collaborative agreement}. \emph{Ieee
  Transactions on Intelligent Transportation Systems} 20(1):8315131,
  \urlprefix\url{http://dx.doi.org/10.1109/TITS.2018.2801823}.

\bibitem[{Dulebenets, Golias, \protect\BIBand{} Mishra(2018)}]{dulebenets2018a}
Dulebenets MA, Golias MM, Mishra S, 2018 \emph{A collaborative agreement for
  berth allocation under excessive demand}. \emph{Engineering Applications of
  Artificial Intelligence} 69:76--92,
  \urlprefix\url{http://dx.doi.org/10.1016/j.engappai.2017.11.009}.

\bibitem[{Dulebenets et~al.(2019)Dulebenets, Pasha, Abioye, \protect\BIBand{}
  Kavoosi}]{dulebenets2019a}
Dulebenets MA, Pasha J, Abioye OF, Kavoosi M, 2019 \emph{Vessel scheduling in
  liner shipping: a critical literature review and future research needs}.
  \emph{Flexible Services and Manufacturing Journal} 33(1):43--106,
  \urlprefix\url{http://dx.doi.org/10.1007/s10696-019-09367-2}.

\bibitem[{Dunning, Huchette, \protect\BIBand{}
  Lubin(2017)}]{DunningHuchetteLubin2017}
Dunning I, Huchette J, Lubin M, 2017 \emph{Jump: A modeling language for
  mathematical optimization}. \emph{SIAM Review} 59(2):295--320,
  \urlprefix\url{http://dx.doi.org/10.1137/15M1020575}.

\bibitem[{Ergun, Kuyzu, \protect\BIBand{} Savelsbergh(2007)}]{ergun2007a}
Ergun O, Kuyzu G, Savelsbergh M, 2007 \emph{Reducing truckload transportation
  costs through collaboration}. \emph{Transportation Science} 41(2):206--221,
  \urlprefix\url{http://dx.doi.org/10.1287/trsc.1060.0169}.

\bibitem[{Fagerholt(2001)}]{fagerholt2001a}
Fagerholt K, 2001 \emph{Ship scheduling with soft time windows: An optimisation
  based approach}. \emph{European Journal of Operational Research}
  131(3):559--571,
  \urlprefix\url{http://dx.doi.org/10.1016/S0377-2217(00)00098-9}.

\bibitem[{Fagerholt et~al.(2015)Fagerholt, Gausel, Rakke, \protect\BIBand{}
  Psaraftis}]{fagerholt2015a}
Fagerholt K, Gausel NT, Rakke JG, Psaraftis HN, 2015 \emph{Maritime routing and
  speed optimization with emission control areas}. \emph{Transportation
  Research, Part C: Emerging Technologies} 52:57--73,
  \urlprefix\url{http://dx.doi.org/10.1016/j.trc.2014.12.010}.

\bibitem[{Fagerholt, Laporte, \protect\BIBand{} Norstad(2010)}]{fagerholt2010a}
Fagerholt K, Laporte G, Norstad I, 2010 \emph{Reducing fuel emissions by
  optimizing speed on shipping routes}. \emph{Journal of the Operational
  Research Society} 61(3):523--529,
  \urlprefix\url{http://dx.doi.org/10.1057/jors.2009.77}.

\bibitem[{Fagerholt \protect\BIBand{} Psaraftis(2015)}]{fagerholt2015b}
Fagerholt K, Psaraftis HN, 2015 \emph{On two speed optimization problems for
  ships that sail in and out of emission control areas}. \emph{Transportation
  Research. Part D: Transport and Environment} 39:56--64,
  \urlprefix\url{http://dx.doi.org/10.1016/j.trd.2015.06.005}.

\bibitem[{Frisk et~al.(2010)Frisk, Göthe-Lundgren, Jörnsten,
  \protect\BIBand{} Rönnqvist}]{frisk2010a}
Frisk M, Göthe-Lundgren M, Jörnsten K, Rönnqvist M, 2010 \emph{Cost
  allocation in collaborative forest transportation}. \emph{European Journal of
  Operational Research} 205(2):448--458,
  \urlprefix\url{http://dx.doi.org/10.1016/j.ejor.2010.01.015}.

\bibitem[{Guan \protect\BIBand{} Cheung(2005)}]{guan2005a}
Guan Y, Cheung RK, 2005 \emph{The berth allocation problem: Models and solution
  methods}. \emph{Container Terminals and Automated Transport Systems:
  Logistics Control Issues and Quantitative Decision Support} 141--158,
  \urlprefix\url{http://dx.doi.org/10.1007/3-540-26686-0_6}.

\bibitem[{Hansen \protect\BIBand{} Oguz(2003)}]{hansen2015a}
Hansen P, Oguz C, 2003 \emph{A note on formulation of the static and dynamic
  berth allocation problems}. \emph{Les Cahiers Du Gerad} 30:1--17.

\bibitem[{Imai, Nishimura, \protect\BIBand{} Papadimitriou(2003)}]{imai2003a}
Imai A, Nishimura E, Papadimitriou S, 2003 \emph{Berth allocation with service
  priority}. \emph{Transportation Research Part B: Methodological}
  37(5):437--457,
  \urlprefix\url{http://dx.doi.org/10.1016/S0191-2615(02)00023-1}.

\bibitem[{Imai et~al.(2005)Imai, Sun, Nishimura, \protect\BIBand{}
  Papadimitriou}]{imai2005a}
Imai A, Sun X, Nishimura E, Papadimitriou S, 2005 \emph{Berth allocation in a
  container port: using a continuous location space approach}.
  \emph{Transportation Research Part B-methodological} 39(3):199--221,
  \urlprefix\url{http://dx.doi.org/10.1016/j.trb.2004.04.004}.

\bibitem[{IMO(2018)}]{IMO2018}
IMO, 2018 \emph{Initial {IMO} strategy on reduction of {GHG} emissions from
  ships}. Technical Report MEPC.304(72), International Maritime Organization,
  \urlprefix\url{http://www.imo.org/en/OurWork/Environment/PollutionPrevention/AirPollution/Pages/GHG-Emissions.aspx},
  (Accessed on 04.05.2020).

\bibitem[{IMO(2020)}]{IMO2020}
IMO, 2020 \emph{Reduction of {GHG} emissions from ships. {F}ourth {IMO GHG}
  study 2020.} Technical Report MEPC.304(72), International Maritime
  Organization,
  \urlprefix\url{https://docs.imo.org/Documents/Detail.aspx?did=125134},
  (Accessed on 22.08.2020).

\bibitem[{Iris, Pacino, \protect\BIBand{} Røpke(2017)}]{iris2017a}
Iris C, Pacino D, Røpke S, 2017 \emph{Improved formulations and an adaptive
  large neighborhood search heuristic for the integrated berth allocation and
  quay crane assignment problem}. \emph{Transportation Research. Part E:
  Logistics and Transportation Review} 105:123--147,
  \urlprefix\url{http://dx.doi.org/10.1016/j.tre.2017.06.013}.

\bibitem[{Iris et~al.(2015)Iris, Pacino, Røpke, \protect\BIBand{}
  Larsen}]{iris2015a}
Iris C, Pacino D, Røpke S, Larsen A, 2015 \emph{Integrated berth allocation
  and quay crane assignment problem: Set partitioning models and computational
  results}. \emph{Transportation Research. Part E: Logistics and Transportation
  Review} 81:75--97,
  \urlprefix\url{http://dx.doi.org/10.1016/j.tre.2015.06.008}.

\bibitem[{Jans(2010)}]{jans2010a}
Jans R, 2010 \emph{Classification of dantzig-wolfe reformulations for binary
  mixed integer programming problems}. \emph{European Journal of Operational
  Research} 204(2):251--254,
  \urlprefix\url{http://dx.doi.org/10.1016/j.ejor.2009.11.014}.

\bibitem[{Kontovas \protect\BIBand{} Psaraftis(2011)}]{kontovas2011a}
Kontovas C, Psaraftis, 2011 \emph{Reduction of emissions along the maritime
  inter modal container chain: Operational models and policies}. \emph{Maritime
  Policy and Management} 38(4):451--469,
  \urlprefix\url{http://dx.doi.org/10.1080/03088839.2011.588262}.

\bibitem[{Kordić et~al.(2016)Kordić, Davidović, Kovač, \protect\BIBand{}
  Dragović}]{kordi2016a}
Kordić S, Davidović T, Kovač N, Dragović B, 2016 \emph{Combinatorial
  approach to exactly solving discrete and hybrid berth allocation problem}.
  \emph{Applied Mathematical Modelling} 40(21-22):8952--8973,
  \urlprefix\url{http://dx.doi.org/10.1016/j.apm.2016.05.004}.

\bibitem[{Krajewska et~al.(2008)Krajewska, Kopfer, Laporte, Ropke,
  \protect\BIBand{} Zaccour}]{krajewska2008a}
Krajewska MA, Kopfer H, Laporte G, Ropke S, Zaccour G, 2008 \emph{Horizontal
  cooperation among freight carriers: request allocation and profit sharing}.
  \emph{Journal of the Operational Research Society} 59(11):1483--1491,
  \urlprefix\url{http://dx.doi.org/10.1057/palgrave.jors.2602489}.

\bibitem[{Kramer et~al.(2019)Kramer, Lalla-Ruiz, Iori, \protect\BIBand{}
  Voß}]{kramer2019a}
Kramer A, Lalla-Ruiz E, Iori M, Voß S, 2019 \emph{Novel formulations and
  modeling enhancements for the dynamic berth allocation problem}.
  \emph{European Journal of Operational Research} 278(1):170--185,
  \urlprefix\url{http://dx.doi.org/10.1016/j.ejor.2019.03.036}.

\bibitem[{Lalla-Ruiz et~al.(2016)Lalla-Ruiz, Expósito-Izquierdo,
  Melián-Batista, \protect\BIBand{} Moreno-Vega}]{lalla-ruiz2016b}
Lalla-Ruiz E, Expósito-Izquierdo C, Melián-Batista B, Moreno-Vega JM, 2016
  \emph{A set-partitioning-based model for the berth allocation problem under
  time-dependent limitations}. \emph{European Journal of Operational Research}
  250(3):1001--1012,
  \urlprefix\url{http://dx.doi.org/10.1016/j.ejor.2015.10.021}.

\bibitem[{Lalla-Ruiz, Melián-Batista, \protect\BIBand{}
  Moreno-Vega(2016)}]{lalla-ruiz2016a}
Lalla-Ruiz E, Melián-Batista B, Moreno-Vega JM, 2016 \emph{A cooperative
  search for berth scheduling}. \emph{Knowledge Engineering Review}
  31(05):498--507, \urlprefix\url{http://dx.doi.org/10.1017/s0269888916000266}.

\bibitem[{Lim(1998)}]{lim1998a}
Lim A, 1998 \emph{The berth planning problem}. \emph{Operations Research
  Letters} 22(2-3):105--110,
  \urlprefix\url{http://dx.doi.org/10.1016/S0167-6377(98)00010-8}.

\bibitem[{{Maersk}(2021)}]{maersk}
{Maersk}, 2021 \emph{{Maersk}}. \url{https://www.maersk.com/}, accessed:
  2021-09-15.

\bibitem[{Magnanti, Orlin, \protect\BIBand{} Ahuja(1993)}]{magnanti1993a}
Magnanti TL, Orlin JB, Ahuja RK, 1993 \emph{Network flows : theory, algorithms,
  and applications} (Prentice-Hall), ISBN 013617549x, 9780136175490.

\bibitem[{Meisel \protect\BIBand{} Bierwirth(2009)}]{meisel2009a}
Meisel F, Bierwirth C, 2009 \emph{Heuristics for the integration of crane
  productivity in the berth allocation problem}. \emph{Transportation Research
  Part E: Logistics and Transportation Review} 45(1):196--209,
  \urlprefix\url{http://dx.doi.org/10.1016/j.tre.2008.03.001}.

\bibitem[{Meng \protect\BIBand{} Wang(2011)}]{meng2011a}
Meng Q, Wang S, 2011 \emph{Optimal operating strategy for a long-haul liner
  service route}. \emph{European Journal of Operational Research}
  215(1):105--114, \urlprefix\url{http://dx.doi.org/10.1016/j.ejor.2011.05.057,
  10.1016/j.ejor.2011.05.057}.

\bibitem[{{Navis}(2021)}]{navis}
{Navis}, 2021 \emph{{Navis}}. \url{https://www.navis.com/}, accessed:
  2021-09-07.

\bibitem[{Notteboom et~al.(2017)Notteboom, Parola, Satta, \protect\BIBand{}
  Pallis}]{notteboom2017a}
Notteboom TE, Parola F, Satta G, Pallis AA, 2017 \emph{The relationship between
  port choice and terminal involvement of alliance members in container
  shipping}. \emph{Journal of Transport Geography} 64:158--173,
  \urlprefix\url{http://dx.doi.org/10.1016/j.jtrangeo.2017.09.002}.

\bibitem[{Notteboom \protect\BIBand{} Vernimmen(2009)}]{notteboom2009a}
Notteboom TE, Vernimmen B, 2009 \emph{The effect of high fuel costs on liner
  service configuration in container shipping}. \emph{Journal of Transport
  Geography} 17(5):325--337,
  \urlprefix\url{http://dx.doi.org/10.1016/j.jtrangeo.2008.05.003}.

\bibitem[{\"Ozener, Ergun, \protect\BIBand{} Savelsbergh(2011)}]{oezener2011a}
\"Ozener OO, Ergun O, Savelsbergh M, 2011 \emph{Lane-exchange mechanisms for
  truckload carrier collaboration}. \emph{Transportation Science} 45(1):1--17,
  \urlprefix\url{http://dx.doi.org/10.1287/trsc.1100.0327}.

\bibitem[{Packard(1995)}]{packard1995shipping}
Packard W, 1995 \emph{Shipping Pools}. Business of Shipping S (Lloyd's of
  London Press), ISBN 9781850445128,
  \urlprefix\url{https://books.google.dk/books?id=g7ArAQAAIAAJ}.

\bibitem[{Pang \protect\BIBand{} Liu(2014)}]{pang2014a}
Pang KW, Liu J, 2014 \emph{An integrated model for ship routing with
  transshipment and berth allocation}. \emph{Iie Transactions (institute of
  Industrial Engineers)} 46(12):1357--1370,
  \urlprefix\url{http://dx.doi.org/10.1080/0740817X.2014.889334}.

\bibitem[{Portchain(2021)}]{portchain}
Portchain, 2021 \url{https://www.portchain.com}, accessed: 2021-09-06.

\bibitem[{Psaraftis \protect\BIBand{} Kontovas(2013)}]{psaraftis2013a}
Psaraftis HN, Kontovas CA, 2013 \emph{Speed models for energy-efficient
  maritime transportation: A taxonomy and survey}. \emph{Transportation
  Research Part C-emerging Technologies} 26:331--351,
  \urlprefix\url{http://dx.doi.org/10.1016/j.trc.2012.09.012}.

\bibitem[{Psaraftis \protect\BIBand{}
  Kontovas(2015{\natexlab{a}})}]{psaraftis2015b}
Psaraftis HN, Kontovas CA, 2015{\natexlab{a}} \emph{Green maritime
  transportation: Speed and route optimization}. \emph{International Series in
  Operations Research and Management Science} 226:299--349,
  \urlprefix\url{http://dx.doi.org/10.1007/978-3-319-17175-3_9}.

\bibitem[{Psaraftis \protect\BIBand{}
  Kontovas(2015{\natexlab{b}})}]{psaraftis2015a}
Psaraftis HN, Kontovas CA, 2015{\natexlab{b}} \emph{Slow steaming in maritime
  transportation: Fundamentals, trade-offs, and decision models}.
  \emph{Handbook of Ocean Container Transport Logistics} 315--358,
  \urlprefix\url{http://dx.doi.org/10.1007/978-3-319-11891-8_11}.

\bibitem[{Pujats, Golias, \protect\BIBand{} Konur(2020)}]{pujats2020a}
Pujats K, Golias M, Konur D, 2020 \emph{A review of game theory applications
  for seaport cooperation and competition}. \emph{Journal of Marine Science and
  Engineering} 8(2):100, \urlprefix\url{http://dx.doi.org/10.3390/jmse8020100}.

\bibitem[{{RBS}(2021)}]{rbs}
{RBS}, 2021 \emph{{Realtime Business Solutions}}.
  \url{https://www.rbs-tops.com/}, accessed: 2021-09-06.

\bibitem[{Reinhardt et~al.(2016)Reinhardt, Plum, Pisinger, Sigurd,
  \protect\BIBand{} Vial}]{reinhardt2016a}
Reinhardt LB, Plum CE, Pisinger D, Sigurd MM, Vial GT, 2016 \emph{The liner
  shipping berth scheduling problem with transit times}. \emph{Transportation
  Research. Part E: Logistics and Transportation Review} 86:116--128,
  \urlprefix\url{http://dx.doi.org/10.1016/j.tre.2015.12.006}.

\bibitem[{Saadaoui, Umang, \protect\BIBand{}
  Frejinger(2015)}]{saadaoui2015column}
Saadaoui Y, Umang N, Frejinger E, 2015 \emph{A column generation framework for
  berth scheduling at port terminals} (CIRRELT, Centre interuniversitaire de
  recherche sur les r{\'e}seaux d'entreprise~…).

\bibitem[{Saeed \protect\BIBand{} Larsen(2010)}]{saeed2010a}
Saeed N, Larsen OI, 2010 \emph{An application of cooperative game among
  container terminals of one port}. \emph{European Journal of Operational
  Research} 203(2):393--403,
  \urlprefix\url{http://dx.doi.org/10.1016/j.ejor.2009.07.019}.

\bibitem[{{Sealytix}(2021)}]{sealytix}
{Sealytix}, 2021 \emph{{Sealytix}}. \url{https://www.sealytix.com/}, accessed:
  2021-09-20.

\bibitem[{Shapley(1953)}]{shapley1953value}
Shapley LS, 1953 \emph{A value for n-person games}. \emph{Contributions to the
  Theory of Games} 2(28):307--317.

\bibitem[{{Ship \& Bunker}(2021)}]{FuelPrice}
{Ship \& Bunker}, 2021 \emph{{Global Average Bunker Price - VLSFO}}.
  \url{https://shipandbunker.com/prices/av/global/av-glb-global-average-bunker-price#VLSFO},
  accessed: 2021-03-21.

\bibitem[{Song(2003)}]{song2003a}
Song DW, 2003 \emph{Port co-opetition in concept and practice}. \emph{Maritime
  Policy and Management} 30(1):29--44,
  \urlprefix\url{http://dx.doi.org/10.1080/0308883032000051612}.

\bibitem[{Song \protect\BIBand{} Panayides(2002)}]{song2002a}
Song DW, Panayides PM, 2002 \emph{A conceptual application of cooperative game
  theory to liner shipping strategic alliances}. \emph{Maritime Policy and
  Management} 29(3):285--301,
  \urlprefix\url{http://dx.doi.org/10.1080/03088830210132632}.

\bibitem[{Steenken, Voß, \protect\BIBand{} Stahlbock(2004)}]{steenken2004a}
Steenken D, Voß S, Stahlbock R, 2004 \emph{Container terminal operation and
  operations research - a classification and literature review}. \emph{Or
  Spectrum} 26(1):3--49,
  \urlprefix\url{http://dx.doi.org/10.1007/s00291-003-0157-z}.

\bibitem[{Sun et~al.(2018)Sun, Niu, Xu, \protect\BIBand{} Ying}]{sun2018a}
Sun B, Niu B, Xu H, Ying W, 2018 \emph{Cooperative optimization for port and
  shipping line with unpredictable disturbance consideration}. \emph{Icnc-fskd
  2018 - 14th International Conference on Natural Computation, Fuzzy Systems
  and Knowledge Discovery} 8686901, 113--118,
  \urlprefix\url{http://dx.doi.org/10.1109/FSKD.2018.8686901}.

\bibitem[{{TGI}(2021)}]{tgi}
{TGI}, 2021 \emph{{TGI Maritime Software}}. \url{https://www.tgims.com/},
  accessed: 2021-09-06.

\bibitem[{UNCTAD(2019)}]{UNCTAD19}
UNCTAD, 2019 \emph{Review of Maritime Transport 2019}.
  \urlprefix\url{http://dx.doi.org/https://doi.org/https://doi.org/10.18356/17932789-en},
  (Accessed on 27.06.2020).

\bibitem[{Venturini et~al.(2017)Venturini, Iris, Kontovas, \protect\BIBand{}
  Larsen}]{venturini2017a}
Venturini G, Iris C, Kontovas CA, Larsen A, 2017 \emph{The multi-port berth
  allocation problem with speed optimization and emission considerations}.
  \emph{Transportation Research. Part D: Transport and Environment}
  54:142--159, \urlprefix\url{http://dx.doi.org/10.1016/j.trd.2017.05.002}.

\bibitem[{Wang, Wang, \protect\BIBand{} Meng(2014)}]{wang2014a}
Wang H, Wang S, Meng Q, 2014 \emph{Simultaneous optimization of schedule
  coordination and cargo allocation for liner container shipping networks}.
  \emph{Transportation Research Part E: Logistics and Transportation Review}
  70(1):261--273, \urlprefix\url{http://dx.doi.org/10.1016/j.tre.2014.07.005}.

\bibitem[{Wang, Liu, \protect\BIBand{} Qu(2015)}]{wang2015a}
Wang S, Liu Z, Qu X, 2015 \emph{Collaborative mechanisms for berth allocation}.
  \emph{Advanced Engineering Informatics} 29(3):572, 332--338,
  \urlprefix\url{http://dx.doi.org/10.1016/j.aei.2014.12.003}.

\bibitem[{Wang \protect\BIBand{} Meng(2012)}]{wang2012a}
Wang S, Meng Q, 2012 \emph{Liner ship route schedule design with sea
  contingency time and port time uncertainty}. \emph{Transportation Research
  Part B} 46(5):615--633,
  \urlprefix\url{http://dx.doi.org/10.1016/j.trb.2012.01.003,
  10.1016/j.trb.2012.01.003}.

\bibitem[{Wang, Meng, \protect\BIBand{} Liu(2013{\natexlab{a}})}]{wang2013a}
Wang S, Meng Q, Liu Z, 2013{\natexlab{a}} \emph{Bunker consumption optimization
  methods in shipping: A critical review and extensions}. \emph{Transportation
  Research Part E: Logistics and Transportation Review} 53(1):49--62,
  \urlprefix\url{http://dx.doi.org/10.1016/j.tre.2013.02.003}.

\bibitem[{Wang, Meng, \protect\BIBand{} Liu(2013{\natexlab{b}})}]{wang2013b}
Wang S, Meng Q, Liu Z, 2013{\natexlab{b}} \emph{A note on "berth allocation
  considering fuel consumption and vessel emissions"}. \emph{Transportation
  Research Part E: Logistics and Transportation Review} 49(1):48--54,
  \urlprefix\url{http://dx.doi.org/10.1016/j.tre.2012.07.002}.

\bibitem[{Wen et~al.(2019)Wen, Larsen, Røpke, Petersen, \protect\BIBand{}
  Madsen}]{wen2019a}
Wen M, Larsen R, Røpke S, Petersen HL, Madsen OB, 2019 \emph{Centralised
  horizontal cooperation and profit sharing in a shipping pool}. \emph{Journal
  of the Operational Research Society} 70(5):737--750,
  \urlprefix\url{http://dx.doi.org/10.1080/01605682.2018.1457481}.

\bibitem[{Wen, Pacino, \protect\BIBand{} Kontovas(2017)}]{wen2017a}
Wen M, Pacino D, Kontovas C, 2017 \emph{A multiple ship routing and speed
  optimization problem under time, cost and environmental objectives}.
  \emph{Transportation Research. Part D: Transport and Environment}
  52(A):303--321, \urlprefix\url{http://dx.doi.org/10.1016/j.trd.2017.03.009}.

\bibitem[{{World Shipping Council}(2015)}]{shippingcouncil2015}
{World Shipping Council}, 2015 \emph{Some observations on port congestion,
  vessel size and vessel sharing agreements}. Technical report, World Shipping
  Council,
  \urlprefix\url{https://www.worldshipping.org/industry-issues/transportation-infrastructure/Observations_on_Port_Congestion_Vessel_Size_and_VSA_May_28_2015.pdf},
  (Accessed on 07.03.2021).

\bibitem[{Zhen et~al.(2020)Zhen, Hu, Yan, Zhuge, \protect\BIBand{}
  Wang}]{zhen2020a}
Zhen L, Hu Z, Yan R, Zhuge D, Wang S, 2020 \emph{Route and speed optimization
  for liner ships under emission control policies}. \emph{Transportation
  Research Part C: Emerging Technologies} 110:330--345,
  \urlprefix\url{http://dx.doi.org/10.1016/j.trc.2019.11.004}.

\end{thebibliography}


\end{document}